\documentclass[10pt,letter]{amsart}

\usepackage[T1]{fontenc}
\usepackage[utf8x]{inputenc}
\usepackage[english]{babel}
\usepackage{amsmath,amsfonts,amssymb,amsthm,graphicx,xspace,esint}
\usepackage{lscape}
\usepackage{qtree}
\usepackage[clearempty]{titlesec}
\usepackage{lmodern}
\usepackage{url}
\usepackage{enumerate}
\usepackage{tabularx}
\usepackage{fullpage}
\usepackage{verbatim}%
\usepackage{lscape}%
\makeatletter%
\def\timenow{\@tempcnta\time
  \@tempcntb\@tempcnta
  \divide\@tempcntb60
  \ifnum10>\@tempcntb0\fi\number\@tempcntb
  \multiply\@tempcntb60
  \advance\@tempcnta-\@tempcntb
:\ifnum10>\@tempcnta0\fi\number\@tempcnta}
\makeatother

\newtheorem{defi}{Definition}[section]
\newtheorem{prop}[defi]{Proposition}
\newtheorem{thm}[defi]{Theorem}
\newtheorem{lem}[defi]{Lemma}

\newtheorem{rmq}[defi]{Remark}

\newtheorem{innercustomthm}{Theorem}
\newenvironment{customthm}[1]
  {\renewcommand\theinnercustomthm{#1}\innercustomthm}
  {\endinnercustomthm}

\newenvironment{pr}{\noindent\textbf{Proof}~:\\}{$\blacksquare$\\}
\newenvironment{prf}[2]{\noindent\textbf{Proof of #1 #2}~:\\}{$\blacksquare$\\}

\newcommand{\N}{\mathbb{N}}
\newcommand{\Z}{\mathbb{Z}}

\newcommand{\R}{\mathbb{R}}
\newcommand{\C}{\mathbb{C}}
\newcommand{\F}{\mathcal{F}}

\newcommand{\E}{\mathcal{E}}
\newcommand{\B}{\mathcal{B}}
\newcommand{\M}{\mathcal{M}}
\newcommand{\He}{\mathcal{H}}
\newcommand{\Sc}{\mathcal{S}}

\newcommand{\Res}{\mathcal{R}}
\newcommand{\m}{\underline{m}}
\newcommand{\n}{\underline{n}}
\newcommand{\p}{\underline{p}}

\newcommand{\norh}[2]{\left\Vert #1 \right\Vert _{H^{#2}}}
\newcommand{\norhh}[2]{\left\lVert #1 \right\rVert _{\dot{H}^{#2}}}
\newcommand{\norlp}[3]{\left\lVert #1 \right\rVert _{L^{#2}_{#3}}}
\newcommand{\norw}[4]{\left\lVert #1 \right\rVert _{W^{#2,#3}_{#4}}}
\newcommand{\nor}[2]{\left\lVert #1 \right\rVert _{#2}}

\newcommand{\T}{\mathcal{T}}
\newcommand{\inte}[4]{\int_{#1}^{#2}{#3 d#4}}

\newcommand{\eps}{\varepsilon}

\newcommand{\xc}{\langle x\rangle}
\newcommand{\cro}[1]{\langle #1\rangle}

\newcommand{\tld}[1]{\widetilde{#1}}
\newcommand{\sgn}{\mathrm{sign}}
\newcommand{\supp}{\mathrm{supp}}
\newcommand{\cacher}[1]{}
\newcommand{\sqnor}[1]{\mathcal{B}(#1)}

\numberwithin{equation}{section}
\usepackage{hyperref}

\title{Long-time existence and resonant approximation for the quadratic nonlinear wave equation with an anisotropic harmonic trapping}

\author{Nicolas Laillet}

\begin{document}
\maketitle%

\begin{abstract}
We establish long-time existence and uniqueness for the 2D wave equation with a harmonic potential in one direction. This proof relies on a fine study of the so-called space-time resonances of the equation. Then we derive a resonant system for this equation and we prove that it is a satisfying approximation for the original equation.
\end{abstract}

\part{Introduction}
\section{Framework}\label{framework}

\subsection{Presentation of the equation}

The goal of this paper is to study the equation
\begin{equation}\label{equation}
\left\{
\begin{array}{rcl}
\partial_{t}^2u-\Delta u+x_{2}^2 u+u&=&u^2,\\
u(0,x_1,x_2)&=&u_{0}(x_1,x_2),\\
\partial_{t}u(0,x_1,x_2)&=&u_{1}(x_1,x_2),
\end{array}\right.
\end{equation}
where $u:(t,x_1,x_2)\in\R_{+}\times\R^2\mapsto u(t,x_1,x_2)\in\R$. This equation is based on a wave equation with a harmonic potential in one direction. \\

We want to study this equation in the weakly nonlinear regime, i.e. in the small data regime. This framework, quite classical in the study of nonlinear dispersive PDEs, allows to obtain long-time existence theorems. In addition to this long-time study, we want to try to understand the long-time dynamics of the solution, approximating it by a simpler one.\\

The study of wave equations with a potential has a pretty long history: for a review of the different dispersive effects, see the work of Schlag in \cite{Sch07}, or \cite{GV03} and \cite{DS10} for specific and more recent examples. Some global existence theorems have been proven in the case of polynomially decreasing potentials: see the books of Strauss (\cite{Str89}) and of Shatah and Struwe (\cite{SS98}) for reviews or \cite{GI92} for a specific example. These results mainly rely on the fact that a localized or decreasing-at-infinity potential should be invisible for solutions far from the origin: its effect should be either neglictible or well-understood from a global point of view. This is not the case for a harmonic potential, and other methods have to be considered.
Introducing a harmonic potential (i.e. non-decaying, non-localized) in a dispersive equation has been studied in the past years, in the particular case of Schrödinger's equation: see for example \cite{Car03}, \cite{AMSW05}, or more recently in \cite{HPTV13} (which considers a toric geometry, quite close to the geometry created by the harmonic potential) or \cite{HT14}. Considering a harmonic potential forces to consider the harmonic structure of the equation and to study the frequencies interactings, i.e. the resonances, inside the nonlinearity.\\
This fine study of resonances has been introduced by Klainerman in \cite{Kla86} and developed for example in \cite{Kla93}. To be more precise, we are going to use the new version of this study of resonances, developped by Germain, Masmoudi and Shatah in \cite{GMS09}, and used for the wave equation by Pusateri and Shatah in \cite{PS13}. \\
\\
Studying the equation $(\ref{equation})$ is therefore quite new. Equation $(\ref{equation})$ is weakly nonlinear, but with a quadratic nonlinearity, which means that the resonant interactions will not be able to be compensated simply by using the weakness of the nonlinearity. Moreover, the geometry given by the harmonic potential, which is physically know as the "cigar-shaped" geometry in the case of Bose-Einstein condensates (for a cubic Schrödinger equation with a harmonic potential) gives birth to very specific resonant interactions and will force us to understand in detail the resonant zones in the frequency space.\\
This fine understanding allows us to understand better the dynamics of our equation: in particular we are able to find a resonant system for $(\ref{equation})$, in the spirit of what has been done in \cite{HPTV13} or \cite{HT14}. This system should be simpler to study and we prove that it is a good approximation of the solutions of the original equation.\\
\\
The rest of Section \ref{framework} is devoted to presenting the main results of this paper: the existence theorem \ref{thm}, the existence theorem for the resonant system \ref{thm-existence-resonant-system} and Theorem \ref{approximation theorem} establishing the validity of the approximation by the resonant system. Section \ref{strategy} gives a strategy for the proof of Theorem \ref{thm}, which is detailed in Sections \ref{strategy-HM}, \ref{I-HF}, \ref{I-HHe} and \ref{I-LFHe}. Then the three next sections focus on the resonant system: the way to obtain it (Section \ref{derivation-resonant-system}), the existence theorem (Section \ref{long-time-existence-resonant system}) and the approximation theorem (Section \ref{section-validity-approx}).

\subsection{Mathematical framework}

Since we want to prove long-time existence theorems in a high-regularity framework, we have to define the regularity spaces we are going to use: given the anisotropic structure of the operator $-\Delta+x^2_{2}$, we have to define anisotropic spaces and anisotropic transforms adapted to the differential operator.

\begin{defi}
The $n$-th Hermite function $\psi_{n}$ on $\R$ is defined as follows
\begin{align}\label{def-Hermite}
\psi_{n}(x):=(-1)^n e^{\frac{-x^2}{2}}\frac{d^n}{dx^n}\left(e^{-x^2}\right).
\end{align}
It is the $n$-th eigenfunction of the harmonic oscillator:
\begin{align}
\psi_{n}''(x)+x^2\psi_{n}=(2n+1)\psi_n.
\end{align}
We also define the interaction term between three Hermite functions $\M(m,n,p)$:
\begin{equation}\label{def-interaction-term}
\mathcal{M}(m,n,p):=\inte{\R}{}{\psi_m(x)\psi_n(x)\psi_p(x)}{x}.
\end{equation}
\end{defi}
\begin{rmq}
We recall that the family $(\psi_n)_{n\in\N}$ is a hilbertian basis of $L^2(\R)$. However, contrary to what happens for complex exponentials, the product of two Hermite functions is not a Hermite function. That is why we need to define the $\mathcal{M}(m,n,p)$: its properties are studied in Appendix \ref{Fourier-multipliers}.
\end{rmq}
\begin{defi}\label{ondes-harm-harmonic-transforms}
The Fourier transform of a function $g$ defined on $\R$ is given by
\begin{equation*}
\F(g)(\xi)=\hat{g}(\xi):=\inte{\R}{}{e^{-ix\xi}g(x)}{x}.
\end{equation*}
 The \emph{Fourier-Hermite} transform of a function $f$ defined on $\R^2$ is defined by 
 \begin{equation*}
\widetilde{f}_{p}(t,\xi)=\inte{\R}{}{  \inte{\R}{}{  f(t,x_1,x_2)e^{-ix_1 \xi}\psi_{p}(x_2)  }{x_1}  }{x_2}, 
\end{equation*}
where $\psi_{p}$ is the $p$-th Hermite function defined in $(\ref{def-Hermite})$.\\
We also define $f_p:=\F_{x_1}^{-1}(\tld{f_p})$.
\end{defi}

Given the form of the operator $(-\Delta+x_{2}^2+1)$, anisotropic regularity spaces will have to be defined. We are going to define two different kinds of "Hermite regularity spaces", depending on whether or not we give a global definition or a strongly anisotropic one: actually both will be related one to the other.\\
The isotropic point of view consists in defining the regularity with the operator, in the same fashion that, for example, $\nor{f}{H^s}=\nor{(-\Delta)^{s/2}f}{2}$:

\begin{defi}
For all integers $N$, for all $f$ in $L^{2}(\R^2)$, we define the $\tld{\He}^N(\R^2)$ norm of $f$ by
\begin{align*}
\nor{f}{\tld{\He}^{N}}=\nor{(-\Delta+ x_{2}^2+1)^{N}f}{L^{2}(\R^2)}.
\end{align*}
\end{defi}

However, given the anisotropy, it will be useful to be able to study each direction separately.
\begin{defi}
Let $M$ be an integer. The space at Hermite regularity $M$, written~$\He^M(\R)$, is defined by the following norm:
\begin{align}\label{def-HeM}
\nor{f}{\He^M}&:=\nor{(-\partial_{x}^2+ x^{2}+1)^{M}f}{L^{2}(\R)}.
\end{align}
\end{defi}
This one-dimensional definition allows us to introduce the following space of functions defined on $\R^2$.
\begin{defi}
Let $N,M$ be two integers. The space $\He_{x_2}^{M}H_{x_1}^{N}$ (written $\He^{M}H^{N}$ in the rest of the paper) is defined on functions on $\R^2$ by the following norm
\begin{align*}
\nor{f}{\He^{M}H^N}=\nor{\left(p^{M}\nor{f_{p}}{H^{N}_{x_1}}\right)_{p\in\N}}{\ell^2}.
\end{align*}
\end{defi}
\begin{rmq}
The space $\He^M H^N$ is related to $\He^M$ in the following way : Since the eigenvalues of $-\partial_{x}^2+ x^{2}+1$ are $(2n+2)_{n\in\N}$, 
the following equivalence of norms holds
\begin{align*}
\nor{f}{\He^M}&\sim \left(\sum_{p\in\N}(2p+2)^{2M}|f_{p}|^2\right)^{\frac{1}{2}}\sim \left(\sum_{p\in\N}\p^{2M}|f_{p}|^2\right)^{\frac{1}{2}},
\end{align*}
with $\p=\max(1,p)$.
\end{rmq}
\begin{rmq}
The order chosen to define $\He^{M}H^N$ is very important: in our proof we are going to start to work with a given Hermite mode, and then sum over all the Hermite modes.
\end{rmq}

Finally, the following spaces will be defined so as to have lighter and simpler notations for time-dependent functions.

\begin{defi}\label{defispace}
Let $M$ and $N$ be two integers, $t$ a non negative real number. Then the spaces $B_{t}$, $\mathcal{B}^{M}_{t}$ and $S^{M,N}_{t}$ are defined, for all function $f$ defined on $\R_{+}\times\R^2$, by the norms
\begin{align}
\nor{f(t)}{B_{t}}&:=\cro{t}^{-\frac{1}{2}}\nor{f(t)}{H^{\frac{3}{2}}(\cro{x_1})},\\
\nor{f(t)}{\mathcal{B}^{M}_{t}}&:=\cro{t}^{-\frac{1}{2}}\nor{f(t)}{\He^{M}H^{\frac{3}{2}}(\cro{x_1})},\\
\nor{f(t)}{S^{M,N}_{t}}&:= \nor{f(t)}{\tld{\He}^N}+\nor{f(t)}{\B^{M}_{t}}.
\end{align}
Finally, the $S^{M,N}_{t}$ norm of a vector is defined as follows:
\begin{align}\label{def-norm-S-R2}
\nor{\left(\begin{array}{c}
u\\
v
\end{array}\right)}{S^{M,N}_{t}}=\nor{u}{S^{M,N}_{t}}+\nor{v}{S^{M,N}_{t}}.
\end{align}
The fixed-point space in which we are going to work is then defined as follows: if $T>0$, $\Sigma^{M,N}_{T}$ is defined by 
\begin{equation}
\nor{f}{\Sigma^{M,N}_{T}}=\sup_{0\leq t\leq T} \nor{f(t)}{S^{M,N}_{t}}.
\end{equation}
\end{defi}

\begin{rmq}\label{rmqnorm}
\begin{itemize}

\item We have for all $f$
\begin{align*}
\nor{f}{\He^{\frac{N}{2}}H^{N}}\leq \nor{f}{\tld{\He}^N}\leq\nor{f}{\He^{N}H^{2N}}.
\end{align*}
\item if $f\in S_t^{N,M}$, then for all $p\in\N$, there exists $(a_{p}(t))_{p\in\N}$, such that
\begin{align*}
\norh{f_p(t)}{N}{}&=\p^{-\frac{N}{2}}a_p(t) \nor{f(t)}{S_t^{N,M}},\\
\nor{f_{p}(t)}{B_t}&=\p^{-M}a_p(t) \nor{f(t)}{S_t^{N,M}},
\end{align*}
with $\nor{(a_p(t))_{p\in\N}}{\ell^2}\leq 1$.
\end{itemize}
\end{rmq}

\subsection{Main results}
Our aim is to obtain a \emph{long-time existence} result for \emph{small initial data} at \emph{high regularity}.
\subsubsection{Existence theorem}
First of all, direct energy estimates give the following result, proven in \cite{These}:

\begin{prop}
Let $\eps>0$. Let $u_0$ be in $\tld{\He}^N$, $N\geq 0$, with $\nor{u_0}{\tld{\He}^N}\leq \eps/2$. Then Equation $(\ref{equation})$ has a unique solution in the space $L^{\infty}([0,T),\tld{\He}^N)$ with 
\begin{align*}
T=C\nor{u_0}{\tld{\He}^N}^{-1},
\end{align*}
wth $C$ independent of $u_0$, such that 
\begin{align*}
 \sup_{t\in [0,T)} \nor{u(t)}{\tld{\He}^N}\leq \eps.
 \end{align*} 
\end{prop}

The whole point of this article is to be able to have a long-time existence, that is to say an existence time of order $\eps^{-a}$ with $a>1$, if $\eps$ is the size of the initial data. We are going to prove the following Theorem:

\begin{customthm}{A}\label{thm}
Let $\delta>0$, $\eps>0$. Let
\begin{align*}
 T=C(\delta)\eps^{-a},~\text{with}~a=\frac{4}{3(1+\delta)},
 \end{align*} 
with $C(\delta)$ depending on $\delta$ only. \\
Then, given $M$ and $N$ integers satisfying
\begin{align}
\label{cond-M}&M>3,\\
\label{cond-N}&N>\frac{1}{\delta}+\frac{3}{2}+2M,
\end{align} 
if $(u_0,u_1)$ satisfies 
\begin{equation*}
\nor{u_0}{S_{0}^{N+1,M+1}}+\nor{u_1}{S_{0}^{N,M}}\leq\frac{\eps}{2},
\end{equation*} then there exists a unique solution $u$ in $\Sigma^{M+1,N+1}_{T}$ to $(\ref{equation})$ with 
\begin{equation}
\begin{aligned}
\nor{u}{\Sigma^{M+1,N+1}_{T}}&\leq\eps,\\
\nor{\partial_{t}u}{\Sigma^{M,N}_{T}}&\leq\eps.
\end{aligned}
\end{equation}
\end{customthm}

This result comes from a fine study of resonance phenomena occurring in the equation and of the dispersive properties of the Klein-Gordon operator with a harmonic potential.

\subsubsection{Notations}
\begin{itemize}
\item For all real numbers $x$, 
\begin{align*}
\cro{x}:=\sqrt{1+x^2}.
\end{align*}
\item For all $\eta\in\R$ and $m\in\N$,
\begin{align*}
\cro{\eta}_{m}:=\sqrt{\eta^2+2m+2}
\end{align*}
\item we write $D=i\partial$.
\item we write $f\lesssim g$ if there exists a universal constant $C$ such that $f\leq Cg$.
\item if $A$ and $B$ are two functions, $m$ and $n$ two integers, 
\begin{align}\label{defi-lesssim-sym}
\left(A\lesssim_{m\leftrightarrow n} B(m,n)\right)\Leftrightarrow\left(A\lesssim B(m,n)+B(n,m)\right).
\end{align}
\end{itemize}

\subsubsection{Duhamel formula}
In order to understand the resonances of $(\ref{equation})$ and find a resonant system for it, we are going to find its Duhamel formula, i.e. its integral formulation.\\
So as to establish this Duhamel formula, we are going to work with the profile of a solution instead of a solution itself.

\begin{defi}
Let $u:\R_{+}\times\R^{d}\rightarrow\C$ be a solution in $L^2$ of $\partial_{t}u+iL(D)u=N(u)$, where $N$ is a nonlinearity and $L(D)$ is a Fourier multiplier, i.e. for all function $f$, $\F(L(D)f)(\xi)=L(\xi)f(\xi)$ with $L:\R^d\rightarrow\C$ is a function of $\xi$, then the \emph{profile} of $u$ is defined by
\begin{equation}
e^{itL(D)}u(t,x):=\frac{1}{(2\pi)^d}\F_{\xi}^{-1}\left(e^{itL(\xi)}\hat{u}(t,\xi)\right).
\end{equation}
\end{defi}

\begin{rmq}
 In the linear case, i.e. $\partial_{t}u+iLu=0$ and $u_0$ the initial condition, the profile of $u$ is $u_0$. This is the way the profile should be understood: it is the solution transported backwards by the linear propagator.
\end{rmq}

Generally, a first step to write a Duhamel formulation for a PDE of order $2$ in time is to reduce it to an equation of order one: this is the reason for the next definition.

\begin{defi} The \emph{left-traveling} part of $u$ (resp. the \emph{right-traveling} part) denoted $u_{+}$ (resp. $u_{-}$) is defined by
\begin{equation*}
  u_{\pm}:=\partial_{t} u\pm i\left(-\Delta+x_2^2+1\right)^{1/2}u.
\end{equation*}
\end{defi}

\begin{rmq}\label{rmq-equivalence-udtu-upm}
It is important to remark that we have the following equivalence for all $t$:
\begin{equation}
\left((u(t),\partial_{t}u(t))\in S_{t}^{N+1,M+1}\times S_{t}^{N,M}\right)\Leftrightarrow \left(u_{\pm}\in S_{t}^{N,M}\right).
\end{equation}
\end{rmq}

We are now going to reformulate the equation as follows (the proof for this formula is in \cite{These}):
\begin{prop}
 A function $u$ is a solution of $(\ref{equation})$ if and only if the profile $f=(f_{+},f_{-})$ satisfies
 \begin{equation}\label{duh}
 f=\mathcal{A}(f),
\end{equation}
where $\mathcal{A}(f)=f_0+A(f)$ and
\begin{align*}
\widetilde{A(f)}_{p}(t,\xi)=\left(\begin{array}{c}
\widetilde{A(f)}_{+,p}(t,\xi)\\
\widetilde{A(f)}_{-,p}(t,\xi)\end{array}\right),
\end{align*}
where
\begin{align}\label{duh-defA}
\widetilde{A(f)}_{\pm,p}(t,\xi)=\sum_{m,n}\sum_{\alpha,\beta=\pm 1}\alpha\beta\mathcal{M}(n,m,p)\inte{0}{t}{  \inte{\R}{}{  e^{\mp is\phi^{\alpha,\beta}_{m,n,p}}\frac{\widetilde{f}_{\alpha,m}(\eta)}{\cro{\eta}_{m}} \frac{\widetilde{f}_{\beta,n}(\xi-\eta)}{\cro{\xi-\eta}_{n}}  }{\eta}  }{s},
\end{align}
with $\phi^{\alpha,\beta}_{m,n,p}=\cro{\xi}_{p}+\alpha\cro{\eta}_{m}+\beta\cro{\xi-\eta}_{n}$ and $\M$ is the hermite functions interaction term $(\ref{def-interaction-term})$.
\end{prop}
\begin{rmq}
Here, the frequency variables $\xi$ or $\eta$ have to be understood as $\xi_{1}$ or $\eta_{1}$, i.e. the frequency variable associated to the first space variable. We do not write this index in order to simplify notations.\\
Moreover, from now on we are going to write $\widetilde{f}_{\alpha,m}(\eta)$ instead of $\widetilde{f}_{\alpha,m}(s,\eta)$ when the dependence in $s$ is obvious.
\end{rmq}

\subsubsection{Remarks on Theorem \ref{thm} and resonant system.}
One can remark that a similar result to Theorem \ref{thm} can be proven quite easily if instead of the space $S^{M,N}_{t}$ we take the space defined by
\begin{align*}
\nor{f}{\tld{S}^{M,N}_{t}}:=\nor{f}{\tld{\He}^N}+\frac{1}{\cro{t}^{\frac{1}{4}}}\nor{f}{\He^M W^{\frac{3}{2},1}}.
\end{align*}
Proving this result relies on the dispersive estimate for the Klein-Gordon with a harmonic potential, i.e. a decay inequality for the $L^{\infty}$ norm of the solution.\\
However it does not involve any the study of resonances. It is really weaker than the result of Theorem \ref{thm}, because it does not give any estimate on $f$ in a weighted Sobolev norm. And having an estimate on weighted norms is fundamental when it comes to the study of the dynamics of the system, in particular when approximating it by a resonant system, as it is done in Part \ref{part-resonant-system}.\\
In fact, in order to study the dynamics of (\ref{equation}), it may be useful to derive a simpler system, generating the dynamics of the original equation. This study of a resonant system has been initiated in hyperbolic equations  by Klainerman and Majda in \cite{KM81} for incompressible fluids, then by Grenier \cite{Gre97} and Schochet \cite{Sch94} with the so-called "filtering method" for highly rotating fluids. \\
Using this notion of resonant system in the framework of dispersive equations is more recent: Ionescu and Pausader in \cite{Io-Pa-12-a} and \cite{Io-Pa-12-b} studied the nonlinear Schrödinger equation on $\R\times \mathbb{T}^3$; other studies have been made by Hani, Pausader, Tzvetkov, Visciglia in \cite{HPTV13} for NLS in $\R\times\mathbb{T}^d$ ($1\leq d\leq 4$), by Pausader, Tzvetkov and Wang for NLS on $\mathbb{S}^3$ (\cite{PTW14}) and more recently by Hani and Thomann in for the NLS with a harmonic trapping (\cite{HT14}).\\
In all of these articles, the main idea is that the dynamics of a system is governed by the resonant frequencies. So as to understand this idea, assume that a quadratic dispersive PDE has the following Duhamel formula:
\begin{align}\label{initial-equation}
\hat{f}(t,\xi)=\hat{f}(0,\xi)+\inte{0}{t}{  \inte{}{}{e^{is\phi(\xi,\eta)}\hat{f}(\eta)\hat{f}(\xi-\eta)}{\eta}  }{s}.
\end{align}

Assume that for all $\xi$, there is exactly one $\eta_{0}(\xi)$ such that $\partial_{\eta}\phi (\xi,\eta_{0}(x))=0$. Then, a Stationary Phase Lemma will give

\begin{align*}
\inte{}{}{e^{is\phi(\xi,\eta)}\hat{f}(\eta)\hat{f}(\xi-\eta)}{\eta}=&\frac{C}{\sqrt{s}}e^{is\phi(\xi,\eta_{0}(\xi))}\hat{f}(\eta_{0}(\xi))\hat{f}(\xi-\eta_{0}(\xi))\\
&+\text{remainder decreasing with time.}
\end{align*}
Moreover, if $\phi(\xi,\eta_{0}(\xi))=0$, the integral
\begin{align*}
\inte{0}{t}{e^{is\phi(\xi,\eta_{0}(\xi))}\hat{f}(\eta_{0}(\xi))\hat{f}(\xi-\eta_{0}(\xi))}{s}
\end{align*}
is an oscillating integral, and it is bounded as $t$ goes to infinity thanks to Riemann-Lebesgue's Lemma. Hence, if $\Xi$ is the set of $\xi$ such that $\phi(\xi,\eta_{0}(\xi))=0$ and $\mathbb{I}_{\xi\in\Xi}$ the indicator function of $\Xi$, the leading term in the Duhamel formula is 
\begin{align*}
\mathbb{I}_{\xi\in\Xi}
\inte{0}{t}{\hat{f}(\eta_{0}(\xi))\hat{f}(\xi-\eta_{0}(\xi))}{s}.
\end{align*}
We will call the equation
\begin{align}\label{resonant-equation-ex}
\hat{f}(t,\xi)=\hat{f}(0,\xi)+\mathbb{I}_{\xi\in\Xi}
\inte{0}{t}{\hat{f}(\eta_{0}(\xi))\hat{f}(\xi-\eta_{0}(\xi))}{s}
\end{align}
the \emph{resonant equation}. This equation is simpler since we restricted the original one to some \emph{resonant modes}. In the case of anisotropic models as ours (with one free direction and one direction trapped by a harmonic potential), the resonant equation keeps the same form but with trickier resonant conditions.\\
A good resonant system has to satisfy the two following properties:
\begin{enumerate}
\item it has to be a good approximation of the initial equation, i.e. if $f$ is a solution of (\ref{initial-equation}) and $g$ is a solution of (\ref{resonant-equation-ex}) with the same initial data, then $f-g$ goes to zero as $t$ goes to infinity (if we are in the lucky case of a global existence).
\item we should be able to understand its dynamics. For example, in \cite{HPTV13}, Hani, Pausader, Tzvetkov, Visciglia were able to build solutions of the resonant system with growing Sobolev norms, and consequently prove that the initial equation had solutions with growing Sobolev norms.
\end{enumerate}

In this paper we focus on the derivation of the resonant system, the existence of long-time solutions for this system and the validity of the approximation of the initial equation by this system.\\
The resonant equation associated to (\ref{equation}) is
\begin{equation}\label{resonant-equation}
\tld{f}_{\pm,p}(t,\xi)=\tld{f}_{\pm,p}(0,\xi)+\inte{0}{t}{
\sum_{\substack{m,n\in\Z\\ \alpha,\beta\in\{\pm 1\}\\ m\neq n\text{ or }\alpha\neq -\beta\\ (\ref{cond-res})\text{ satisfied}}}\frac{ \mathcal{M}(m,n,p)}{\sqrt{s|\partial^{2}_{\eta}\phi(\xi,\lambda^{\alpha,\beta}_{m,n}\xi)|}}\frac{\widetilde{f}_{\alpha,m}(\lambda^{\alpha,\beta}_{m,n}\xi)}{\cro{\lambda^{\alpha,\beta}_{m,n}\xi}_{m}} \frac{\widetilde{f}_{\beta,n}((1-\lambda^{\alpha,\beta}_{m,n})\xi)}{\cro{(1-\lambda^{\alpha,\beta}_{m,n})\xi}_{n}} 
}{s},
\end{equation}
with $\lambda^{\alpha,\beta}_{m,n}:=\frac{1}{1+\alpha\beta\sqrt{\frac{n+1}{m+1}}}$ and $(\ref{cond-res})$ is the resonant condition appearing in Theorem \ref{thm-phase}. It will be formally derived in Section \ref{derivation-resonant-system}\\
Since the system is simpler given it involves only selected interacting modes, we are able to prove a better existence and uniqueness result than for the original one.

\begin{customthm}{B}\label{thm-existence-resonant-system}
Let $\eps>0$, $T=C/\eps^2$ with $C$ a universal constant. Then, given $M$, $N$ and $\kappa$ satisfying
\begin{equation}
M>6,~N\geq\frac{3}{2},~\kappa=1\text{ or }2,
\end{equation}if $f_0=(f_{0,+},f_{0,-})$ is an initial data with $\nor{f_{0}}{\He^M H^N(\cro{x_1}^\kappa)}\leq \eps/2$, then there exists one and only one solution $f=(f_{+},f_{-})$ to the resonant system on the interval $[0,T)$, belonging to the $L^{\infty}_{t}\left([0,T), \He^M H^N(\cro{x_1}^{\kappa}) \right)$. Moreover, for all $t\in[0,T)$, $\nor{f_{\pm}}{\He^M H^N(\cro{x_1}^{\kappa})}\leq\eps$.
\end{customthm}

Moreover we are able to prove that this resonant system is a good approximation of the initial equation:

\begin{customthm}{C}\label{approximation theorem}
Let $0<\alpha<\frac{5}{3}$, $0<\omega<1-\frac{3}{5}\alpha$, $\eps>0$. Let $N$ and $M$ satisfying (\ref{cond-M}), (\ref{cond-N}) with in addition $N\geq 9-\frac{1}{4}$. Let $0\leq M_{0}< M-\frac{1}{8}$.\\
There exists a $C(\alpha,\omega)$ such that, for $\eps$ small enough, if $ T=C(\alpha,\omega)\eps^{-\frac{4}{3+\omega}}$, if 
\begin{itemize}
\item $f$ is a solution to the initial system in $\Sigma^{M,N}_{T}$ with initial data $f_0$ in the ball of center $0$ and radius $\eps/2$ of $S^{M,N}_{0}$,
\item $g$ is a solution to the resonant system in $\Sigma^{M,N}_{T}$ with the same initial data,
\end{itemize}
then we have, for all $t\leq T$, 
\begin{align*}
\nor{(f-g)(t)}{\He^{M_{0}}L^{2}}\leq \eps^{\alpha}.
\end{align*}
\end{customthm}

\begin{rmq}
We have to compare the size of $f-g$ to the variation of $f$ and $g$ during a time $\eps^{-\frac{4}{3+\omega}}$: we prove in Theorem \ref{prop-contraction} that the one for $f$ is of order 
\begin{align*}
\inte{0}{\eps^{-\frac{4}{3+\omega}}}{s^{-\frac{1}{4}}\eps^2}{s}\sim \eps^{2-\frac{3}{3+\omega}}.
\end{align*}
Similarly, the increase of $g$ is of order 
\begin{align*}
\inte{0}{\eps^{-\frac{4}{3+\omega}}}{s^{-\frac{1}{2}}\eps^2}{s}\sim\eps^{2-\frac{2}{3+3\omega}}<<\eps^{2-\frac{3}{3+\omega}}.
\end{align*}
But if $\omega$ is small enough, we have $\eps^{\alpha}<<\eps^{2-\frac{2}{3+3\omega}}$, which means that the size of $f-g$ is small compared to the variation of $f$ and $g$.
\end{rmq}

\section{Strategy}\label{strategy}

\subsection{Contraction estimates}\label{section-contractions}

The Duhamel formula allows us to write Equation $(\ref{equation})$ as a fixed-point problem:
\begin{align*}
f=\mathcal{A}(f),
\end{align*}
where $\mathcal{A}(f)=f_0+A(f)$, and $A$ is defined in $(\ref{duh-defA})$.

It suffices to prove that the operator $\mathcal{A}$ is a contraction which maps the unit ball for the $\Sigma^{M,N}_{T}$ norm into itself (for well-chosen $M$, $N$ and $T$). More precisely, let $\eps>0$. The goal is to prove that if $\nor{f_0}{S^{M,N}_{0}}\leq\frac{\eps}{2}$ then $\nor{f}{\Sigma^{M,N}_{T}}\leq\eps\Rightarrow\nor{\mathcal{A}(f)}{\Sigma^{M,N}_{T}}\leq\eps$. \\
 So as to do this, we have to prove an inequality of the form
\begin{align*}
\nor{A(f)}{\Sigma^{M,N}_{T}}\leq CT^{a}\nor{f}{\Sigma^{M,N}_{T}}^k\leq CT^{a}\eps^k,
\end{align*}
where $a<1$ and $k\in\N$, $k>1$. This will allow existence on a time 
\begin{align*}
T=\left(2C\eps^{k-1}\right)^{-\frac{1}{a}}.
\end{align*}
In order to get the existence time $T=\eps^{-\frac{4}{3(1+\delta)}}$, we are going to prove the following theorem.

\begin{thm}\label{prop-contraction}
For all $\omega>0$, for all $M$ and $N$ satisfying $(\ref{cond-M})-(\ref{cond-N})$, we have the following inequality:
\begin{align}\label{prop-contraction-ineg}
\nor{\mathcal{A}(f)}{S^{M,N}_{t}}&\lesssim \inte{0}{t}{s^{-\frac{1}{4}+\omega}\nor{f}{S^{M,N}_{s}}^2+s^{\frac{1}{2}+\omega}\nor{f}{S^{M,N}_{s}}^3}{s}+C(t),
\end{align}
with $C(t)=\cro{t}^{\omega+\frac{1}{4}}\left(\nor{f(t)}{S^{M,N}_{t}}^2+\nor{f(1)}{S^{M,N}_{1}}^2\right)$.
\end{thm}
\begin{rmq}
We can make two remarks on this property:
\begin{itemize}
\item the following inequality is a direct consequence of $(\ref{prop-contraction-ineg})$:
\begin{align}\label{ineg-contr-2}
\nor{A(f)}{\Sigma^{M,N}_{T}}\lesssim \max\left(T^{\frac{3}{4}+\omega}\nor{f}{\Sigma^{M,N}_{T}}^2,T^{\frac{3}{2}+\omega}\nor{f}{\Sigma^{M,N}_{T}}^3\right).
\end{align}
with the same notations and hypothesis as in Theorem \ref{prop-contraction}.
\item Going from Theorem \ref{prop-contraction} to Theorem \ref{thm} is then quite straightforward: the inequality $(\ref{ineg-contr-2})$ gives an existence time equal to (up to a constant):
\begin{align*}
\min\left(\eps^{\frac{1}{\frac{3}{4}+\omega}},\eps^{2.\frac{1}{\frac{3}{2}+\omega}}\right)=\min(\eps^{-\frac{4}{3+4\omega}},\eps^{-\frac{4}{3+2\omega}})=\eps^{-\frac{4}{3+4\omega}}.
\end{align*}
Then taking $\omega=\frac{3\delta}{4}$ gives the result.
\end{itemize}
\end{rmq}

For now on, our goal will be to prove Theorem \ref{prop-contraction}.

\subsection{Space-time resonances}\label{space-time-resonances}
In order to obtain a large existence time, we are going to use the method of \emph{space-time resonances}, introduced by Germain, Masmoudi and Shatah in \cite{GMS09}. Take a general Duhamel formula:
\begin{align*}
\hat{f}(t,\xi)=\hat{f}_{0}(\xi)+\inte{0}{t}{\inte{\R}{}{e^{-is\phi(\xi,\eta)}\hat{f}(\eta)\hat{f}(\xi-\eta)}{\eta}}{s},
\end{align*}
where $\phi:=L(\xi)-L(\eta)-L(\xi-\eta)$. This corresponds to the dispersive equation $\partial_{t}u+L(D)u=u^2$.\\
So as to deal with long-time existence, we will have to find a way of gaining some powers of time. This is the aim of space-time resonances: we have to study the phase $\phi$ to make some transformations.
\begin{enumerate}
\item if the phase $\phi$ does not vanish, we can write $e^{-is\phi}=\frac{i}{\phi}\partial_{s}\left(e^{-is\phi}\right)$ and perform an integration by parts. This \emph{normal forms transformation} (according to Shatah's terminology in \cite{Sha85}, refering to the classical concept of normal forms in dynamical systems) will lead to an expression of the form
\begin{align*}
 \inte{\R}{}{e^{-it}\hat{f}(t,\eta)\hat{f}(t,\xi-\eta)}{\eta}+ \inte{0}{t}{\inte{\R}{}{\frac{1}{\phi}e^{-is\phi(\xi,\eta)}\partial_{s}\hat{f}(\eta)\hat{f}(\xi-\eta)}{\eta}}{s}+\text{ symmetric or easier terms}.
 \end{align*} A short calculation shows that $\partial_{s}\hat{f}=e^{-isL}\F(u^2)$. The nonlinearity $f\partial_{s}f$ is now a cubic one, easier to deal with since we are working with small data.\\
 Having $\phi=0$ can be understood as $L(\xi)=L(\eta)+L(\xi-\eta)$: this corresponds to what is also called resonances in physics, i.e. a situation when two plane waves interact and create another plane wave.
\item if the derivative of the phase $\partial_{\eta}\phi$ does not vanish, write $e^{-is\phi}=\frac{1}{-is\partial_{\eta}\phi}\partial_{\eta}\left(e^{-is\phi}\right)$ and perform an integration by parts, to obtain
\begin{align*}
\inte{0}{t}{\inte{\R}{}{\partial_{\eta}\left(\frac{1}{-is\partial_{\eta}\phi}\right)e^{-is\phi(\xi,\eta)}\hat{f}(\eta)\hat{f}(\xi-\eta)}{\eta}}{s}&+\inte{0}{t}{\inte{\R}{}{\frac{1}{-is\partial_{\eta}\phi}e^{-is\phi(\xi,\eta)}\partial_{\eta}\hat{f}(\eta)\hat{f}(\xi-\eta)}{\eta}}{s}\\
&+\text{ symmetric or easier terms}.
 \end{align*} This allows to gain a power of $s$ and improve the long-time existence.\\
 Concretely, the case $\partial_{\eta}\phi\neq 0$ has to be seen as $L'(\eta)\neq L'(\xi-\eta)$, i.e. the group velocities of both interacting waves are different: \emph{if we assume the initial data to be localized}, this guarantees that a wave packet at frequency $\eta$ and a wave packet at frequency $\xi-\eta$ will not interact with each other after a sufficiently long time.\\
 This method has first been developped by Klainerman in \cite{Kla86}.
\item in the general case, we have zones where $\phi$ vanishes (call this zone $\T$, the time resonant set), where $\partial_{\eta}\phi$ vanishes ($\Sc$, the space resonant set) and a zone where both vanish: call it $\Res$, the space-time resonant set. This zone is problematic because none of the integration by parts presented before are feasable.\\
If the space-time resonant set is of measure $0$, it is reasonable to think that this difficulty will be handled by an \emph{adapted localization}: we introduce a function $\chi$ equal to $1$ around $\Res$: the integral $\inte{0}{t}{\inte{\R}{}{\chi(\xi,\eta)e^{-is\phi(\xi,\eta)}\hat{f}(\eta)\hat{f}(\xi-\eta)}{\eta}}{s}$ should be small, maybe if we make the localization narrower as time grows. Then one of the integration by parts (in $s$ or in $\eta$) should give estimates for $\inte{0}{t}{\inte{\R}{}{(1-\chi)e^{-is\phi(\xi,\eta)}\hat{f}(\eta)\hat{f}(\xi-\eta)}{\eta}}{s}$, maybe if the $(\xi,\eta)$ plane is cut off again in two zones, one where $\phi$ does not vanish and one where $\partial_{\eta}\phi$ does not vanish.\\
So as to estimating terms of the form 
$\inte{0}{t}{\inte{\R}{}{m(\xi,\eta)e^{-is\phi(\xi,\eta)}\hat{f}(\eta)\hat{f}(\xi-\eta)}{\eta}}{s}$, we will need some results on \emph{Fourier bilinear multiplier}: the ones needed in this article are gathered in Appendix \ref{Fourier-multipliers}.
\end{enumerate}

All the difficulty will be to deal with these different zones: we have to identify them (study of the phase), cut off carefully the frequency space and deal with each zone. This is the aim of the following sections. We start by defining the resonant sets adapted to the method just evoked.
\begin{defi}
Let $\phi(\xi,\eta)$ be a real phase.
\begin{enumerate}[(i)]
\item The \emph{time resonant set} is the set $\T:=\lbrace(\xi,\eta)\in\R^2| \phi(\xi,\eta)=0\rbrace$.
\item The \emph{space resonant set} is the set $\Sc:=\lbrace(\xi,\eta)\in\R^2|\partial_{\eta}\phi(\xi,\eta)=0\rbrace$.
\item The \emph{space coresonant set} is the set $\tld{\Sc}:=\lbrace(\xi,\eta)\in\R^2|\partial_{\xi}\phi(\xi,\eta)=0\rbrace$.
\item The \emph{space-time resonant set} is the set $\mathcal{R}:=\lbrace(\xi,\eta)\in\R^2|\phi(\xi,\eta)=\partial_{\eta}(\xi,\eta)=0\rbrace$.
\end{enumerate}
\end{defi}
\begin{rmq}
The space coresonant set is not directly involved in the study of resonances. However, since the fixed-point space is built on weighted norms, we will have to differentiate with respect to $\xi$ the integral term in the Duhamel formula $(\ref{duh-defA})$: this will make terms of the form $\partial_{\xi}\phi$ appear. Hence it can be interesting to compare $\partial_{\xi}\phi$ to $\partial_{\eta}\phi$: this explains that we also want to study the space coresonant set $\tld{\Sc}$. 
\end{rmq}

The precise study of the resonances and of the resonant sets is done in Appendix \ref{resonances-asympt}: here we only state the main theorem.

\begin{thm}\label{thm-phase}
Let $\alpha$ and $\beta$ be two elements of $\lbrace -1,1\rbrace$. Consider the phase 
\begin{equation*}
\phi(\xi,\eta)=\phi^{\alpha,\beta}_{m,n,p}(\xi,\eta)=\sqrt{\xi^2+2p+2}+\alpha\sqrt{\eta^2+2m+2}+\beta\sqrt{(\xi-\eta)^2+2n+2}.
\end{equation*}
Then the space resonant set is
\begin{equation*}
\Sc=\left\lbrace \left(\left(1+\frac{\beta}{\alpha}\sqrt{\frac{2n+2}{2m+2}}\right)\eta,\eta\right),\eta\in\R\right\rbrace.
\end{equation*}
\begin{enumerate}
\item In the case $(\alpha,\beta)=(1,1)$, there are no time resonances: $\T=\emptyset$.
\item Otherwise,
  \begin{enumerate}
    \item If $\alpha\beta p+\beta m<0$ or $\alpha\beta p+\beta n<0$, there are no time resonances.
    \item If $\alpha\beta p+\beta m\geq 0$ and $\alpha\beta p+\beta n\geq 0$, there are space-time resonances if and only if the following condition is satisfied.
    \begin{equation}\label{cond-res}\tag{C}
	\left\lbrace\begin{array}{c}
		\alpha\beta p+\beta m+\alpha n\geq 0,\\
		m^2+n^2+p^2-2mn-2pm-2pn-2m-2n-2p-3=0.
	\end{array}\right.
	\end{equation}
    In that case, the space coresonant set is equal to the space resonant set: $\tld{\Sc}=\Sc$.
  \end{enumerate}
\end{enumerate}
\end{thm}

Hence we are in a rather new situation, compared to the previously studied situations. For example, in \cite{GMS12-2D}, Germain, Masmoudi and Shatah study the nonlinear Schrödinger equation with a quadratic nonlinearity which corresponds to the following phases:
\begin{itemize}
\item if the nonlinearity is equal to $\bar{u}^2$, then the phase is $\phi=\xi^2+\eta^2+(\xi-\eta)^2$. The time resonant set is $\{\xi=\eta=0\}$, so is the space-time resonant one.
\item if the nonlinearity is equal to $u^2$, then the then the phase is $\phi=\xi^2-\eta^2-(\xi-\eta)^2$. Hence the space resonant set is $\{\xi=2\eta\}$ and the time resonant set is $\{\eta.(\xi-\eta)=0\}$: this leads to a space-time resonant set equal to $\{\xi=\eta=0\}$.
\item if the nonlinearity is equal to $u\bar{u}$, then the phase is $\phi=\xi^2-\eta^2+(\xi-\eta)^2$. Hence the space resonant set is $\{\xi=0\}$, the time resonant set is $\{\eta.\xi=0\}$ and the space-time resonant set equal to $\{\xi=0\}$.
\end{itemize}

These considerations explain why Germain, Masmoudi and Shatah focused on a nonlinearity of the form $\bar{u}^2$ or $u^2$: the space-time resonant set being a point, it must be easier to deal with than the space-time resonant set of the case $|u|^2$ (where blow-up is expected).\\
\\
In our situation, like in the nonlinear Schrödinger equation with a $|u|^2$ nonlinearity, the space-time resonant set is a line; one more difficulty is that this line depends on the input and output Hermite modes. This kind of problem is a really new situation which will require a fine adaptation of the Germain-Masmoudi-Shatah method.
\begin{rmq}
Here is how the integers satisfying the condition $C_{\alpha,\beta}$ distribute (in the case $\alpha=\beta=1$: they all are on a surface of degree $3$ but, more interesting, they seem to be uniformly distributed.
\begin{tabularx}{\textwidth}{XX}
   \includegraphics[scale=0.42]{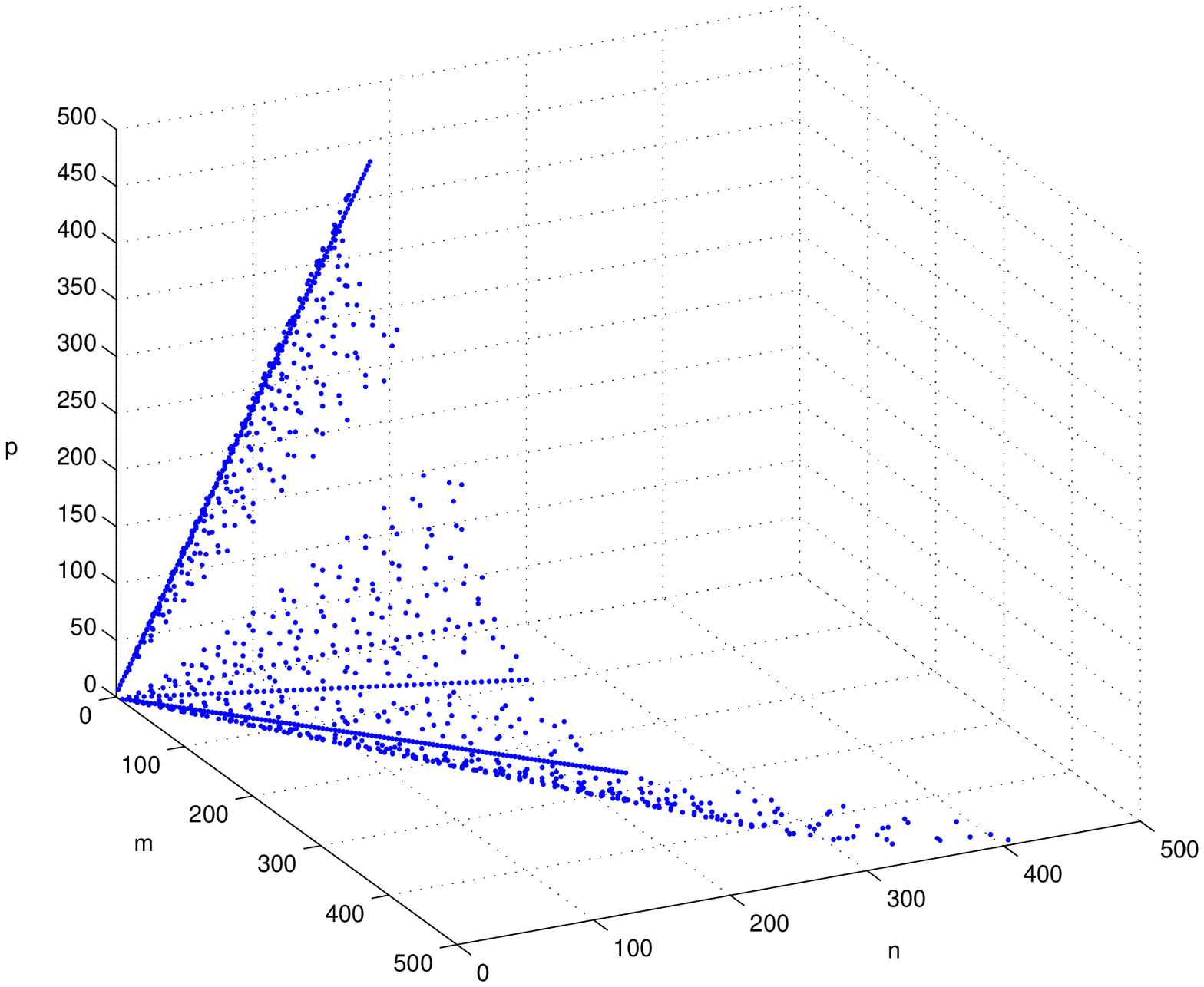} &
   \includegraphics[scale=0.44]{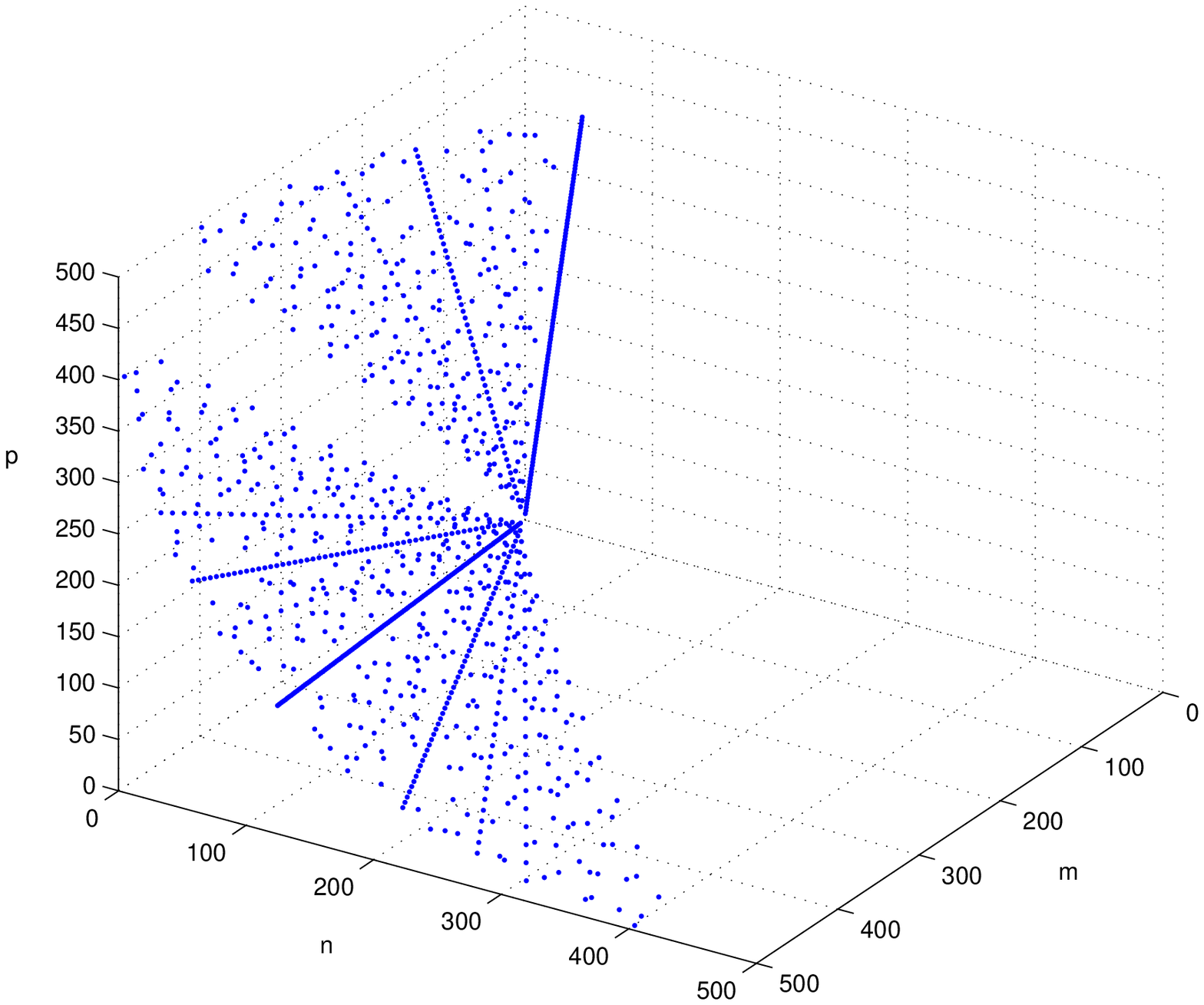} \\
\end{tabularx}
\end{rmq}
\begin{rmq}\label{rmqmass}
This is also the reason why we add a mass to the harmonic potential : if we were to study $\partial_{t}^2u-\Delta u+x_{2}^2 u=u^2$, the resonant structure would be much simpler, since the algebraic resonance condition would be false for all integers $m,n,p$.
\end{rmq}
In order to study separately space and time resonant sets, a precise understanding of the phase and its derivatives is necessary: this work is done in Appendix \ref{resonances-asympt}.\\
This study being done, we introduce the following functions, used in order to explicitly split the $(\xi,\eta)$-frequency space into several zones:
\begin{defi}\label{cutoffs}
We define $\chi$ to be a smooth function, homogeneous of degree $0$ such that $\chi(\eta,\xi-\eta)=0$ for $|\xi-\eta|\leq 2|\eta|$. \\
The function $\theta$ is defined as a smooth function supported on $[0,2]$, equal to $1$ on $[0,1]$. For all $R>0$ we define $\theta_{R}(x)=\theta(\cro{R}x)$. \\
Moreover we impose $\chi$ and $\theta$ to satisfy Coifman-Meyer's theorem's hypotheses.
\end{defi}

\begin{rmq}
The function $\chi$ will allow us to make a \emph{paraproduct decomposition} adapted to the convolution. The function $\theta_{R}$ will be introduced in order deal separately with high and low frequencies, and the parameter $R$ may be time-dependent.
\end{rmq}

\part{Long-time existence and uniqueness}

First of all, we are going to skip the proof of the $\tld{H}^N$ estimate: 

\begin{prop}\label{contrHN}
The operator $\mathcal{B}$ defined as
\begin{equation}\label{HN-def-B}
\mathcal{B}(f):=\inte{0}{t}{e^{is\sqrt{-\Delta+x_{2}^2+1}}\sum_{\alpha,\beta=\pm 1}\frac{e^{-\alpha is\sqrt{-\Delta+x_{2}^2+1}}f_{\alpha}}{\sqrt{-\Delta+x_{2}^2+1}}\frac{e^{-\beta is\sqrt{-\Delta+x_{2}^2+1}}f_{\beta}}{\sqrt{-\Delta+x_{2}^2+1}}}{s}.
\end{equation}
satisfies, for all $f$, for all $M$ satisfiying he condition $(\ref{cond-M})$ and $N$ satisfying $(\ref{cond-N})$.
\begin{equation*}
\nor{\mathcal{B}(f)}{\tld{H}^N}\lesssim \inte{0}{t}{s^{-\frac{1}{4}}\nor{f}{S^{M,N}_{s}}^2}{s}.
\end{equation*}
\end{prop}

The proof can be found in \cite{These} and relies on a product law for the harmonic oscillator (proven in \cite{DG09}) of the form $\nor{fg}{\tld{\He}^{N}}\lesssim \nor{f}{\tld{\He}^{N}}\norlp{g}{\infty}{}+\norlp{f}{\infty}{}\nor{g}{\tld{\He}^{N}}$.\\

\section{High regularity results -- strategy for the weighted norm}\label{strategy-HM}

We now focus on the weighted norm and we want to prove the following proposition:

\begin{prop}\label{prop-H-3/2}
If $U$ is the bilinear operator defined as 
\begin{align*}
\F(U_{m,n}(a,b)):= \inte{0}{t}{  \inte{}{}{  \partial_{\xi}\left(e^{-is\phi}\frac{\hat{a}(\eta)}{\langle \eta\rangle_{m}}\frac{\hat{b}(\xi-\eta)}{\langle \xi-\eta\rangle_{n}}\right)  }{\eta}  }{s}
\end{align*}
then for all $\delta>0$, for all $M$, $N$ satisfying the conditions $(\ref{cond-M})-(\ref{cond-N})$, there exists a constant $K(\delta)$ such that
\begin{equation}\label{ineg-H-3/2}
\begin{aligned}
&\frac{1}{\sqrt{t}}p^M \sum_{m,n} \mathcal{M}(m,n,p)\nor{ U_{m,n}(f_n,f_m)}{H^{3/2}} \\
&\leq  K(\delta)\bigg(\inte{0}{t}{\left(s^{4\delta}s^{-\frac{1}{4}}\nor{f}{S^{M,N}_{s}}^2 a_p(s) + s^{3\delta}s^{\frac{1}{2}}\nor{f}{S^{M,N}_{s}}^3 b_{p}(s)\right)}{s}+C(t)c_{p}(t)\bigg),
\end{aligned}
\end{equation}
with $C(t)=\cro{t}^{\frac{5\delta}{2}+\frac{1}{4}}\left(\nor{f(t)}{S^{M,N}_{s}}^2+\nor{f(1)}{S^{M,N}_{1}}^2\right)$ and where $(a_{p}(s))_{p\in\N}$, $(b_{p}(s))_{p\in\N}$ and $(c_{p}(t))_{p\in\N}$ are $\ell^2$ sequences of norm bounded by $1$.
\end{prop}
\begin{rmq}
Taking the $\ell^2$ norm (in $p$) of (\ref{ineg-H-3/2}) allows us to write
\begin{align*}
\nor{ \sum_{m,n} \mathcal{M}(m,n,p)\frac{1}{\sqrt{t}} U_{m,n}(f_n,f_m)}{\mathcal{B}^M_{t}}  \leq K(\delta)\left(\inte{0}{t}{\left(s^{4\delta}s^{-\frac{1}{4}}\nor{f}{S^{M,N}_{s}}^2  + s^{3\delta}s^{\frac{1}{2}}\nor{f}{S^{M,N}_{s}}^3\right) }{s}+C(t)\right),
\end{align*}
which corresponds to the inequality in Theorem \ref{prop-contraction}.
\end{rmq}
\begin{prf}{Proposition}{\ref{prop-H-3/2}}
We know that
\begin{align*}
 \frac{1}{\sqrt{t}}\norhh{U(f_n,f_m)}{3/2} &=  \frac{1}{\sqrt{t}}\norlp{|\xi|^{\frac{3}{2}}\inte{0}{t}{  \inte{}{}{  \partial_{\xi}\left(e^{-is\phi}\frac{\hat{f}_{m}(\eta)}{\langle \eta\rangle_{m}}\frac{\hat{f}_{n}(\xi-\eta)}{\langle \xi-\eta\rangle_{n}}\right)  }{\eta}  }{s}}{2}{\xi}.
\end{align*}
By the Leibniz rule the integral can be rewritten as follows:
\begin{equation*}
 |\xi|^{\frac{3}{2}}\inte{0}{t}{  \inte{}{}{  \partial_{\xi}\left(e^{-is\phi}\frac{\hat{f}_{m}(\eta)}{\langle \eta\rangle_{m}}\frac{\hat{f}_{n}(\xi-\eta)}{\langle \xi-\eta\rangle_{n}}\right)  }{\eta}  }{s}=I_{m,n}+J_{m,n}+K_{m,n},
\end{equation*}
where
\begin{align*}
& I_{m,n}:=|\xi|^{\frac{3}{2}}\inte{0}{t}{  \inte{}{}{ -is\partial_{\xi}\phi e^{-is\phi}\frac{\hat{f}_{m}(\eta)}{\langle \eta\rangle_{m}}\frac{\hat{f}_{n}(\xi-\eta)}{\langle \xi-\eta\rangle_{n}}  }{\eta}  }{s},\\
& J_{m,n}:=|\xi|^{\frac{3}{2}}\inte{0}{t}{  \inte{}{}{ e^{-is\phi}\frac{\hat{f}_{m}(\eta)}{\langle \eta\rangle_{m}}\frac{\partial_{\xi}\hat{f}_{n}(\xi-\eta)}{\langle \xi-\eta\rangle_{n}}  }{\eta}  }{s},\\
& K_{m,n}:=|\xi|^{\frac{3}{2}}\inte{0}{t}{  \inte{}{}{ e^{-is\phi}\frac{\hat{f}_{m}(\eta)}{\langle \eta\rangle_{m}}\frac{\hat{f}_{n}(\xi-\eta)(-\xi)}{\left((\xi-\eta)^2+2n+2\right)^{3/2}}  }{\eta}  }{s}.
\end{align*}
The integral term $J_{m,n}$ will be treated in Section \ref{J} (Proposition \ref{contrJ}) ; $K_{m,n}$ will be dealt with in Section \ref{K} (Proposition \ref{contrK}). Estimating the integral term $I_{m,n}$ will be harder and explained in Section \ref{I}.\\
Although the estimates for $J_{m,n}$ and $K_{m,n}$ might not seem sharp in view of their proof, they fit in the ones for the $I_{m,n}$ term, which is the harder one to deal with.
\subsection{Estimates for $J_{m,n}$}\label{J}
We are going to prove the following inequality:
\begin{prop}\label{contrJ}
There exists $(a_{p}(s))_{p\in\N}$ in the unit ball of $\ell^2$ such that
\begin{align*}
\frac{1}{\sqrt{t}}p^M\sum_{m,n\in\N}\mathcal{M}(m,n,p)\norlp{J_{m,n}}{2}{}\lesssim \inte{0}{t}{\cro{s}^{-\frac{1}{4}+\frac{3\delta}{2}}\nor{f}{S^{M,N}_{s}}^2 a_{p}}{s},
\end{align*}
for all $N$ and $M$ satisfying (\ref{cond-M})-(\ref{cond-N}).
\end{prop}
We are going to proceed in two steps: first, we establish $n$ and $m$-dependent bounds on $\norlp{J_{m,n}}{2}{}$ ; then it remains to sum these bounds.

\subsubsection{Bounds for $\norlp{J_{m,n}}{2}{}$}
\begin{lem}\label{contrJ-lemma}
The following bound holds.
\begin{align}\label{contrJ-lemma-ineg}
\frac{1}{\sqrt{t}}\norlp{J_{m,n}}{2}{}\lesssim_{m\leftrightarrow n} \inte{0}{t}{\cro{s}^{-\frac{1}{4}+\frac{3\delta}{2}}\frac{\max(\m,\n)^{\frac{1}{4}}}{\sqrt{\m\n}}\Big(\sqnor{f_m}(s)+\nor{f_m(s)}{H^N}\Big)\nor{f_n(s)}{B_s}}{s},
\end{align}
with $\lesssim_{m\leftrightarrow n}$ defined in (\ref{defi-lesssim-sym}) and $\sqnor{f_m}(s)$ is the following quantity:
\begin{equation}\label{def-sqnor}
\sqnor{f}(s):=\sqrt{\nor{f(s)}{H^N}\nor{f(s)}{B_{s}}}.
\end{equation}
\end{lem}

\begin{pr}
Here we will perform two cutoffs:
\begin{enumerate}
\item the paraproduct decomposition, using the function $\chi$, equal to $0$ for $|\xi-\eta|\leq 2|\eta|$.
\item high-low frequency cut-off: $\theta_{s^\delta}(|\eta|)$ is the smooth function localizing in the zone $|\eta|\leq s^{\delta}$.
\end{enumerate}
Write the $L^2$ norm of $J$ as follows.
\begin{align*}
\norlp{J_{m,n}}{2}{}&\leq
 \norlp{|\xi|^{\frac{3}{2}}\inte{0}{t}{  \inte{}{}{ (1-\chi(\eta,\xi-\eta))e^{-is\phi}\frac{\hat{f}_{m}(\eta)}{\langle \eta\rangle_{m}}\frac{\partial_{\xi}\hat{f}_{n}(\xi-\eta)}{\langle \xi-\eta\rangle_{n}}  }{\eta}  }{s}}{2}{\xi}\\
 &+\norlp{|\xi|^{\frac{3}{2}}\inte{0}{t}{  \inte{}{}{ \chi(\eta,\xi-\eta)\theta_{s^\delta}(|\eta^2|)e^{-is\phi}\frac{\hat{f}_{m}(\eta)}{\langle \eta\rangle_{m}}\frac{\partial_{\xi}\hat{f}_{n}(\xi-\eta)}{\langle \xi-\eta\rangle_{n}}  }{\eta}  }{s}}{2}{\xi}\\
 &+\norlp{|\xi|^{\frac{3}{2}}\inte{0}{t}{  \inte{}{}{ \chi(\eta,\xi-\eta)(1-\theta_{s^\delta}(|\eta^2|))e^{-is\phi}\frac{\hat{f}_{m}(\eta)}{\langle \eta\rangle_{m}}\frac{\partial_{\xi}\hat{f}_{n}(\xi-\eta)}{\langle \xi-\eta\rangle_{n}}  }{\eta}  }{s}}{2}{\xi}\\
 &=\norlp{J_{1-\chi}}{2}{}+\norlp{J_{\chi,l}}{2}{}+\norlp{J_{\chi,h}}{2}{}.
\end{align*}
We are going to estimate separately each of those three terms.
\begin{itemize}
\item \textbf{Term $J_{1-\chi}$}\\
We can rewrite $J_{1-\chi}$ as
\begin{align*}
J_{1-\chi}= \inte{0}{t}{  \inte{}{}{\left(\frac{|\xi|}{|\xi-\eta|}\right)^{\frac{3}{2}} (1-\chi(\eta,\xi-\eta))e^{-is\phi}\frac{\hat{f}_{m}(\eta)}{\langle \eta\rangle_{m}}\frac{|\xi-\eta|^{\frac{3}{2}}\partial_{\xi}\hat{f}_{n}(\xi-\eta)}{\langle \xi-\eta\rangle_{n}}  }{\eta}  }{s},
\end{align*}
which allows us to write it as a bilinear Fourier operator:
\begin{align*}
\F^{-1}\left(J_{1-\chi}\right)=\inte{0}{t}{  e^{is\cro{D}_{p}} T_{m^{J}_{1-\chi}}\left(e^{-is\cro{D}_{m}}\frac{f_{m}}{\cro{D}_{m}},|D|^{\frac{3}{2}}e^{-is\cro{D}_{n}}\frac{x_{1}f_{n}}{\cro{D}_{n}}\right)  }{s},
\end{align*}
with $T_{m^{J}_{1-\chi}}$ the bilinear Fourier operator (defined in (\ref{def-bilin-mult}) ) associated to the multiplier 
\begin{equation*}
m^{J}_{1-\chi}(\eta,\zeta):=\left|\frac{\eta +\zeta}{\zeta}\right|^{\frac{3}{2}} (1-\chi(\eta,\zeta)),
\end{equation*}
 which satisfies H\"{o}lder-like inequalities thanks to Proposition \ref{CMex}.\\
The differential operator $e^{is\cro{D}}$ is continuous from $L^2$ to $L^2$. Then 
\begin{align*}
\norlp{J_{1-\chi}}{2}{}&\lesssim \inte{0}{t}{  \norlp{ T_{m^{J}_{1-\chi}}\left(e^{-is\cro{D}_{m}}\frac{f_{m}}{\cro{D}_{m}},|D|^{\frac{3}{2}}e^{-is\cro{D}_{n}}\frac{x_{1}f_{n}}{\cro{D}_{n}}\right)}{2}{}  }{s}.
\end{align*}
Proposition \ref{CMex} implies
\begin{align}\label{J1-chi-ineg}
\norlp{J_{1-\chi}}{2}{}&\lesssim \inte{0}{t}{  \norlp{ e^{-is\cro{D}_{m}}\frac{f_{m}}{\cro{D}_{m}}}{\infty}{}\norlp{|D|^{\frac{3}{2}}e^{-is\cro{D}_{n}}\frac{x_{1}f_{n}}{\cro{D}_{n}}}{2}{}  }{s}.
\end{align}
Proposition \ref{disp+mult} implies that 
\begin{align}\label{J1-chi-inegint-1}
\norlp{ e^{-is\cro{D}_{m}}\frac{f_{m}}{\cro{D}_{m}}}{\infty}{}\lesssim m^{-\frac{1}{4}}\cro{s}^{-\frac{1}{4}}\sqrt{\nor{f_m}{H^N}\nor{f_m}{B_{s}}}.
\end{align}
Thus Proposition \ref{propmultgeneral}, inequality (\ref{propmult}), gives
\begin{align}\label{J1-chi-inegint-2}
 \norlp{|D|^{\frac{3}{2}}e^{-is\cro{D}_{n}}\frac{x_{1}f_{n}}{\cro{D}_{n}}}{2}{}\lesssim n^{-\frac{1}{2}}\norlp{|D|^{\frac{3}{2}}e^{-is\cro{D}_{n}}x_{1}f_{n}}{2}{}\lesssim n^{-\frac{1}{2}}\nor{f_n}{B_{s}}.
 \end{align}
The two inequalities $(\ref{J1-chi-inegint-1})$ and $(\ref{J1-chi-inegint-2})$ allow us to rewrite $(\ref{J1-chi-ineg})$ as
\begin{align}\label{ineg-J-1-chi}
\norlp{J_{1-\chi}}{2}{}&\lesssim m^{-\frac{1}{4}}n^{-\frac{1}{2}}\inte{0}{t}{  \cro{s}^{-\frac{1}{4}}\sqrt{\nor{f_m}{H^N}\nor{f_m}{B_{s}}}\nor{f_n}{B_{s}}  }{s}.
\end{align}
\item \textbf{Term $J_{\chi,l}$}\\
First of all since we are in the zone $\{\sqrt{|\xi|^2+|\eta|^2}\leq \cro{s}^{\delta}\}\cap \{|\xi-\eta|\leq|\eta|\}$, we can bound $|\xi|^{\frac{3}{2}}$ by $\cro{s}^{\frac{3\delta}{2}}$ (up to a constant). This gives
\begin{align*}
\norlp{J_{\chi,l}}{2}{}&\lesssim \norlp{\inte{0}{t}{  \inte{}{}{ \cro{s}^{\frac{3\delta}{2}}\chi(\eta,\xi-\eta)\theta_{s^\delta}(|\eta^2|)e^{-is\phi}\frac{\hat{f}_{m}(\eta)}{\langle \eta\rangle_{m}}\frac{\partial_{\xi}\hat{f}_{n}(\xi-\eta)}{\langle \xi-\eta\rangle_{n}}  }{\eta}  }{s}}{2}{\xi}\\
&\lesssim \inte{0}{t}{  \cro{s}^{\frac{3\delta}{2}}\norlp{ T_{m^{J}_{\chi,l}}\left(e^{-is\cro{D}_{m}}\frac{f_{m}}{\cro{D}_{m}},e^{-is\cro{D}_{n}}\frac{x_{1}f_{n}}{\cro{D}_{n}}\right)}{2}{}  }{s},
\end{align*}
where $m^{J}_{\chi,l}(\xi,\eta):=\chi(\eta,\xi-\eta)\theta_{s^\delta}(|\eta^2|)$. Thanks to Coifman-Meyer estimates (Theorem \ref{CM}) we can write
\begin{align*}
\norlp{J_{\chi,l}}{2}{}&\lesssim \inte{0}{t}{  \cro{s}^{\frac{3\delta}{2}}\norlp{ e^{-is\cro{D}_{m}}\frac{f_{m}}{\cro{D}_{m}}}{\infty}{}\norlp{e^{-is\cro{D}_{n}}\frac{x_{1}f_{n}}{\cro{D}_{n}}}{2}{}  }{s}.
\end{align*}
The conclusion is the same as previously: we find
\begin{align}\label{ineg-J-chi-l}
\norlp{J_{\chi,l}}{2}{} &\lesssim \inte{0}{t}{\cro{s}^{\frac{3\delta}{2}-\frac{1}{4}}m^{-\frac{1}{4}}n^{-\frac{1}{2}}\nor{f_m}{B_{s}}\sqrt{\nor{f_n}{H^N}\nor{f_n}{B_{s}}}}{s}.
\end{align}

\item \textbf{Term $J_{\chi,h}$}\\
Inequality (\ref{hf}) will be crucial to deal with high frequencies. Firstly, $J_{\chi,h}$ can be rewritten as
\begin{align*}
\norlp{J_{\chi,h}}{2}{}&=\norlp{\inte{1}{t}{  \inte{}{}{ \left|\frac{\xi}{\eta}\right|^{\frac{3}{2}}\chi(\eta,\xi-\eta)(1-\theta_{s^\delta}(|\eta^2|))|\eta|^{\frac{3}{2}}e^{-is\phi}\frac{\hat{f}_{m}(\eta)}{\langle \eta\rangle_{m}}\frac{\partial_{\xi}\hat{f}_{n}(\xi-\eta)}{\langle \xi-\eta\rangle_{n}}  }{\eta}  }{s}}{2}{\xi}\\
&\lesssim  \inte{0}{t}{  \norlp{ T_{m^{J}_{\chi,h}}\left((1-\theta_{s^\delta}(|D|))|D|^{\frac{3}{2}}e^{-is\cro{D}_{m}}\frac{f_{m}}{\cro{D}_{m}},e^{-is\cro{D}_{n}}\frac{x_{1}f_{n}}{\cro{D}_{n}}\right)}{2}{}  }{s},
\end{align*}
with $m^{J}_{\chi,h}(\eta,\zeta):=\left|\frac{\xi}{\eta}\right|^{\frac{3}{2}}\chi(\eta,\xi-\eta)(1-\theta_{s^\delta}(|\eta^2|))$. Then Theorem \ref{CM} combined with Proposition \ref{CMex} give
\begin{align}\label{ineg-J-chi-h}
\norlp{J_{\chi,h}}{2}{}&\lesssim  \inte{0}{t}{  \norlp{(1-\theta_{s^\delta}(|D|))|D|^{\frac{3}{2}}e^{-is\cro{D}_{m}}\frac{f_{m}}{\cro{D}_{m}}}{\infty}{}\norlp{e^{-is\cro{D}_{n}}\frac{x_{1}f_{n}}{\cro{D}_{n}}}{2}{}  }{s}.
\end{align}
First, by the multiplier estimate (\ref{propmult}),
\begin{equation}\label{Jchi-inegint-1}
\norlp{e^{-is\cro{D}_{n}}\frac{x_{1}f_{n}}{\cro{D}_{n}}}{2}{}\lesssim n^{-\frac{1}{2}}\nor{f_n(s)}{B_{s}}.
\end{equation}
Then, the Sobolev embedding $\norlp{u}{\infty}{}\lesssim\norh{u}{N}$ for $N>1/2$ -- actually it is enough to use ${\norlp{u}{\infty}{}\lesssim\norh{u}{1}}$ -- allows us to deal with the other factor in (\ref{ineg-J-chi-h}):
\begin{align*}
\norlp{(1-\theta_{s^\delta}(|D|))|D|^{\frac{3}{2}}e^{-is\cro{D}_{m}}\frac{f_{m}}{\cro{D}_{m}}}{\infty}{}&\lesssim \norlp{(1-\theta_{s^\delta}(|D|))|D|^{\frac{5}{2}}e^{-is\cro{D}_{m}}\frac{f_{m}}{\cro{D}_{m}}}{2}{}\\
&\lesssim m^{-\frac{1}{2}}\norlp{(1-\theta_R(D))|D|^{\frac{5}{2}}f_{m}}{2}{}.
\end{align*}
Thanks to Proposition \ref{propmultgeneral}, inequality (\ref{hf}),
\begin{align}\label{Jchi-inegint-2}
\norlp{(1-\theta_{s^\delta}(|D|))|D|^{\frac{3}{2}}e^{-is\cro{D}_{m}}\frac{f_{m}}{\cro{D}_{m}}}{\infty}{}&\lesssim \frac{1}{\sqrt{m}}\cro{s}^{-\left(N-\frac{5}{2}\right)\delta}\nor{f_m}{H^N}.
\end{align}
Finally, using $(\ref{Jchi-inegint-1})$ and $(\ref{Jchi-inegint-2})$ in $(\ref{ineg-J-chi-h})$ gives
\begin{align*}
\frac{1}{\sqrt{t}}\norlp{J_{\chi,h}}{2}{}&\lesssim \inte{0}{t}{\frac{1}{\sqrt{\m\n}}\frac{1}{\cro{s}^{\delta\left(N-\frac{3}{2}\right)}} \norh{f_m}{N}{}\nor{f_n}{B_{s}}}{s}\\
&\lesssim \frac{1}{\sqrt{\m\n}}\inte{0}{t}{\cro{s}^{-\frac{1}{4}} \norh{f_m}{N}{}\nor{f_n}{B_{s}}}{s},
\end{align*}
as soon as $N\geq \frac{3}{2}+\frac{1}{4\delta}$ (which is true because of Hypothesis (\ref{cond-N})).
\end{itemize}
 Combining Inequalities (\ref{ineg-J-1-chi}), (\ref{ineg-J-chi-l}) and (\ref{ineg-J-chi-h}) concludes the proof of Lemma \ref{contrJ-lemma}.\end{pr}

\subsubsection{Summation}\label{resumJ}

Now that Lemma \ref{contrJ-lemma} is proven, going back to Proposition \ref{contrJ} reduces to summing over the indices $m$ and $n$, i.e. to proving the following result:

\begin{align}\label{resumJ-ineq}
\p^M\sum_{m,n}\mathcal{M}(m,n,p)\frac{\max(\m,\n)^{\frac{1}{4}}}{\sqrt{\m\n}}\Big(\sqnor{f_m}(s)+\nor{f_m(s)}{H^N}\Big)\nor{f_n(s)}{B_s} \lesssim  \nor{f(s)}{S^{M,N}_{s}}^2 a_{p}(s) ,
\end{align}
with $\nor{(a_{p}(s))_{p\in\N}}{\ell^2}\leq 1$.\\ \\
\begin{pr}
First of all, for all integers $p$,
\begin{align*}
\nor{f_{p}(s)}{B_{s}}=\nor{f(s)}{S^{M,N}_{s}} \p^{-M}a_{p}(s),
\end{align*}
with $(a_{p}(s))_{p\in\N}$ in the unit ball of $\ell^2$. Similarly, there exists $(b_{p}(s))_{p\in\N}$ in the unit ball of $\ell^2$ such that for all $p\in\N$,
\begin{align*}
\norh{f_p(s)}{N}{}&=\p^{-\frac{N}{2}}\nor{f(s)}{S^{M,N}_{s}}b_{p}(s)\\
&\leq \p^{-M}\nor{f(s)}{S^{M,N}_{s}}b_{p}(s),
\end{align*}
since $N\geq 2M$ thanks to condition (\ref{cond-N}). Finally, using that if $(a_{p}(s))_{p\in\N},(b_{p}(s))_{p\in\N}\in\ell^2$, then $(\sqrt{a_{p}b_{p}(s)})_{p\in\N}$ is in $\ell^2$, we can bound each of the factors
\begin{align*}
 &\sqnor{f_m}(s)\nor{f_n(s)}{B_s},~\nor{f_m(s)}{H^N}\nor{f_n(s)}{B_s},\\\ &\sqnor{f_n}(s)\nor{f_m(s)}{B_s},~\nor{f_m(s)}{B_s}\nor{f_n(s)}{H^N},
 \end{align*} 
by
\begin{align*}
 \m^{-M}\mu_m(s) \n^{-M}\eta_n(s) \nor{f}{S^{M,N}_{s}}^2,
 \end{align*} 
 with $(\mu_m(s))_{m\in\N}$ and $(\eta_n(s))_{n\in\N}$ in the unit ball of $\ell^2$. Then we can bound 
\begin{align*}
\Big(\sqnor{f_m}(s)+\nor{f_m(s)}{H^N}\Big)\nor{f_n(s)}{B_s}\lesssim \m^{-M}\alpha_m(s) \n^{-M}\beta_n(s) \nor{f(s)}{S^{M,N}_{s}}^2,
\end{align*}
with $(\alpha_m(s))_{m\in\N}$ and $(\beta_n(s))_{n\in\N}$ in the unit ball of $\ell^2$. This inequality implies 
\begin{align*}
&\p^M\sum_{m,n}\mathcal{M}(m,n,p)\frac{\max(\m,\n)^{\frac{1}{4}}}{\sqrt{\m\n}}\Big(\sqnor{f_m}(s)+\nor{f_m(s)}{H^N}\Big)\nor{f_n(s)}{B_s} \\
&\lesssim \nor{f(s)}{S^{N,M}_{s}}^2\p^M\sum_{m,n}\mathcal{M}(m,n,p)\frac{\max(\m,\n)^{\frac{1}{4}}}{\sqrt{\m\n}}\m^{-M}\alpha_m(s) \n^{-N}\beta_n(s).
\end{align*}
Then we are in the framework of the resummation theorem \ref{resummation}-\ref{classicalresummation}, with $a=\frac{1}{4}$, hence 
\begin{equation}\label{ineg-resum-J}
\begin{aligned}
\p^M\sum_{m,n}\mathcal{M}(m,n,p)\frac{\max(\m,\n)^{\frac{1}{4}}}{\sqrt{\m\n}}\Big(\sqnor{f_m}(s)+\nor{f_m(s)}{H^N}\Big)\nor{f_n(s)}{B_s} 
&\lesssim  C_{\gamma}\nor{f}{S^{N,M}_{s}}^2 \p^{-\frac{1}{2}+\gamma}a_{p}(s)\\
&\lesssim \nor{f}{S^{M,N}_{s}}^{2} a_{p}(s),
\end{aligned}
\end{equation}

with $\nor{(a_{p}(s))_{p\in\N}}{\ell^2}\leq 1$. Combining $(\ref{ineg-resum-J})$ and the integral inequality $(\ref{contrJ-lemma-ineg})$ ends the proof of Proposition \ref{contrJ}.
\end{pr}

\subsection{Estimates for $K_{m,n}$}\label{K}
Here the property we need to prove is very similar to the one for $J$: this is the reason why the proof of the following property will be skipped (the details are in \cite{These}).

\begin{prop}\label{contrK}
There exists $(a_{p}(s))_{p\in\N}$ in $\ell^2$ such that
\begin{align*}
\frac{1}{\sqrt{t}}\p^M\sum_{m,n\in\N}\mathcal{M}(m,n,p)\norlp{K_{m,n}}{2}{}\lesssim \inte{0}{t}{s^{-\frac{1}{4}}\nor{f}{S^{M,N}_{s}}^2 u_{p}(s)}{s}.
\end{align*}
\end{prop}

\subsection{Estimates for $I_{m,n}$}\label{I}

The integral term $I_{m,n}$ concentrates the main difficulties of the proof: it corresponds to the case when the differentiation in $\xi$ hits the complex exponential $e^{is\phi}$ and appears to make long-time estimates impossible given the additional power of $s$ given by $\partial_{\xi}e^{is\phi}$. This is the reason why we are going to try to find a way to get additional decay in time. We write
\begin{align*}
\sum_{m,n}\mathcal{M}(m,n,p)I_{m,n}&=I^{hf}+I^{hm}+I^{lf,lm},
\end{align*}
where
\begin{enumerate}
\item $I^{hf}$ corresponds to the \emph{high frequency term}:
\begin{align*}
I^{hf}:=\sum_{m,n}\mathcal{M}(m,n,p)I^{hf}_{m,n},
\end{align*}
with
\begin{align*}
I^{hf}_{m,n}:= -|\xi|^{\frac{3}{2}}\inte{0}{t}{  \inte{}{}{ \left(1-\theta_{s^\delta}(|\eta|)\theta_{s^\delta}(|\xi-\eta|)\right)is\partial_{\xi}\phi e^{-is\phi}\frac{\hat{f}_{m}(\eta)}{\langle \eta\rangle_{m}}\frac{\hat{f}_{n}(\xi-\eta)}{\langle \xi-\eta\rangle_{n}}  }{\eta}  }{s},
\end{align*}
the function $\theta_{s^\delta}(|\zeta|)$ being smooth and localizing in the zone $|\zeta|\leq s^{\delta}$.\\
High frequencies are quite easy to deal with: in fact, since we are in a high-regularity framework, this means that the high frequencies have a small amplitude. This can be understood with the high-frequency inequality in Proposition \ref{propmultgeneral}, inequality (\ref{hf}).
\item $I^{hm}$ corresponds to the \emph{high Hermite modes}:
\begin{align*}
I^{hm}:= \sum_{m,n\geq \cro{t}^{\delta}}\mathcal{M}(m,n,p)I^{hm}_{m,n},
\end{align*}
with
\begin{align*}
I^{hm}_{m,n}:= -|\xi|^{\frac{3}{2}}\inte{0}{t}{  \inte{}{}{ \theta_{s^\delta}(|\eta|)\theta_{s^\delta}(|\xi-\eta|)is\partial_{\xi}\phi e^{-is\phi}\frac{\hat{f}_{m}(\eta)}{\langle \eta\rangle_{m}}\frac{\hat{f}_{n}(\xi-\eta)}{\langle \xi-\eta\rangle_{n}}  }{\eta}  }{s}.
\end{align*}
The idea of a high regularity leading to good estimates for high frequencies applies also in the framework of Hermite modes.
\item $I^{lf,lm}$ is the remaining term, corresponding to low frequencies and low Hermite modes:
\begin{align*}
I^{lf,lm}:= \sum_{m\leq \cro{t}^{\delta},n\in\N}\mathcal{M}(m,n,p)I^{lf}_{m,n}+\sum_{m\in\N,n\leq \cro{t}^{\delta}}\mathcal{M}(m,n,p)I^{lf}_{m,n}.
\end{align*}
This last sum will be treated thanks to the space-time resonances method: in particular we are going to distinguish when there are space-times resonances (condition (\ref{cond-res}), page \pageref{cond-res}, satisfied) or when there are not space-time resonances (condition (\ref{cond-res}) not satisfied).
\end{enumerate}
The situation is summed up in the following tree.\\
\\
\begin{center}
\Tree[.{$\displaystyle  \sum_{m,n}\mathcal{M}(m,n,p) I_{m,n}:=\sum_{m,n}\mathcal{M}(m,n,p)|\xi|^{\frac{3}{2}}\inte{0}{t}{  \inte{}{}{ -is\partial_{\xi}\phi e^{-is\phi}\frac{\hat{f}_{m}(\eta)}{\langle \eta\rangle_{m}}\frac{\hat{f}_{n}(\xi-\eta)}{\langle \xi-\eta\rangle_{n}}  }{\eta}  }{s}$} {High Fourier modes\\ term $I^{hf}$\\ section \ref{I-HF}} {High Hermite modes\\ term $I^{hm}$\\ section \ref{I-HHe}} [.{Low Fourier and\\Hermite modes\\ term $I^{lf,lm}$\\ section \ref{I-LFHe}} {Small times\\ \ref{I-LFHE-st}} [.{Large times} {Condition $C$\\satisfied\\\ref{I-LFHE-Lt-C}} {Condition $C$\\ not satisfied\\ \ref{I-LFHE-Lt-notC}} ] ] ] 
\end{center}

\section{High-frequency estimates}\label{I-HF}
Since we are in a high regularity framework, dealing with high frequencies should not be problematic. We need to prove the following proposition.

\begin{prop}\label{contrIhf}
There exists $(a_p(s))_{p\in\N}$ in $\ell^2$, $\nor{(a_p(s))_{p\in\N}}{\ell^2}\leq 1$ such that for all $N,M$ satisfying Conditions $(\ref{cond-M})-(\ref{cond-N})$,
\begin{align*}
\frac{1}{\sqrt{\cro{t}}}\p^M\sum_{m,n\in\N}\mathcal{M}(m,n,p)\norlp{I^{hf}_{m,n}}{2}{}\lesssim \inte{0}{t}{s^{-\frac{1}{4}}\nor{f}{S^{M,N}_{s}}^2 a_{p}(s)}{s}.
\end{align*}
\end{prop}
\begin{pr}

First of all, we are going to prove the following result:
\begin{lem}\label{contrIhf-lemma}
For all integers $m$ and $n$,
\begin{equation*}
 \frac{1}{\sqrt{\cro{t}}}\norlp{I^{hf}_{m,n}}{2}{}\lesssim_{m\leftrightarrow n} \inte{0}{t}{ \frac{\max(\m,\n)^{\frac{1}{4}}}{\sqrt{\m\n}}s^{-\frac{1}{4}}\norh{f_{m}(s)}{N} \sqnor{f_n}(s)}{s},
\end{equation*}
whenever $N$ satisfies $(\ref{cond-N})$.
\end{lem}
The sum over $m$ and $n$ will be skipped: we are in the framework of the summation theorem \ref{resummation}, with the same parameters as in Section \ref{resumJ}: Proposition \ref{contrIhf} is deduced in the same way.\\\\
\begin{prf}{Lemma}{\ref{contrIhf-lemma}}
Recall the expression of $I^{hf}_{m,n}$:
\begin{align*}
I^{hf}_{m,n}:= -|\xi|^{\frac{3}{2}}\inte{0}{t}{  \inte{}{}{ \left(1-\theta_{s^\delta}(|\eta|)\theta_{s^\delta}(|\xi-\eta|)\right)is\partial_{\xi}\phi e^{-is\phi}\frac{\hat{f}_{m}(\eta)}{\langle \eta\rangle_{m}}\frac{\hat{f}_{n}(\xi-\eta)}{\langle \xi-\eta\rangle_{n}}  }{\eta}  }{s}.
\end{align*}
Our idea is to say that if $|\xi-\eta|+|\eta|$ is large, it means that either $|\xi-\eta|$ is large or $|\eta|$ is large. In order to do so, we can remark that
\begin{align*}
1-\theta_{s^\delta}(|\eta|)\theta_{s^\delta}(|\xi-\eta|)=
\left(1-\theta_{s^\delta}(|\eta|)\right)+\theta_{s^\delta}(|\eta|)\left(1-\theta_{s^\delta}(|\xi-\eta|)\right).
\end{align*}

Using the paraproduct decomposition (with the function $\chi$, equal to $0$ for $|\xi-\eta|\leq 2|\eta|$) leads to the following splitting of $I^{hf}_{m,n}$:
\begin{align*}
I^{hf}_{m,n}=I^1 + I^2 + I^3 + I^4,
\end{align*}
where
\begin{align*}
I^1:=& |\xi|^{\frac{3}{2}}\inte{0}{t}{  \inte{}{}{ s\partial_{\xi}\phi\chi(\eta,\xi-\eta) \left(1-\theta_{s^\delta}(|\eta|)\right)e^{-is\phi}\frac{\hat{f}_{m}(\eta)}{\langle \eta\rangle_{m}}\frac{\hat{f}_{n}(\xi-\eta)}{\langle \xi-\eta\rangle_{n}}  }{\eta}  }{s},\\
I^2:=& |\xi|^{\frac{3}{2}}\inte{0}{t}{  \inte{}{}{ s\partial_{\xi}\phi\chi(\eta,\xi-\eta) \theta_{s^\delta}(|\eta|)\left(1-\theta_{s^\delta}(|\xi-\eta|)\right)e^{-is\phi}\frac{\hat{f}_{m}(\eta)}{\langle \eta\rangle_{m}}\frac{\hat{f}_{n}(\xi-\eta)}{\langle \xi-\eta\rangle_{n}}  }{\eta}  }{s},\\
I^{3}:=&|\xi|^{\frac{3}{2}}\inte{0}{t}{  \inte{}{}{s\partial_{\xi}\phi(1-\chi(\eta,\xi-\eta)) \left(1-\theta_{s^\delta}(|\eta|)\right)e^{-is\phi}\frac{\hat{f}_{m}(\eta)}{\langle \eta\rangle_{m}}\frac{\hat{f}_{n}(\xi-\eta)}{\langle \xi-\eta\rangle_{n}}  }{\eta}  }{s},\\ 
I^{4}:=&|\xi|^{\frac{3}{2}}\inte{0}{t}{  \inte{}{}{s\partial_{\xi}\phi(1-\chi(\eta,\xi-\eta)) \theta_{s^\delta}(|\eta|)\left(1-\theta_{s^\delta}(|\xi-\eta|)\right)e^{-is\phi}\frac{\hat{f}_{m}(\eta)}{\langle \eta\rangle_{m}}\frac{\hat{f}_{n}(\xi-\eta)}{\langle \xi-\eta\rangle_{n}}  }{\eta}  }{s}. 
\end{align*}
We shall only deal with the integral term $I^1$, the others can be dealt with in a similar way.\\

Using the expression of $\partial_{\xi}\phi=\frac{\xi}{\cro{\xi}_{p}}-\frac{\xi-\eta}{\cro{\xi-\eta}_{n}}$, we get
\begin{align*}
 I^{1}=&\inte{0}{t}{  \inte{}{}{ \frac{|\xi|^{\frac{3}{2}}}{|\eta|^{\frac{3}{2}}} s\partial_{\xi}\phi\chi(\eta,\xi-\eta) (1-\theta_{s^\delta}(|\eta|))e^{-is\phi}|\eta|^{\frac{3}{2}}\frac{\hat{f}_{m}(\eta)}{\langle \eta\rangle_{m}}\frac{\hat{f}_{n}(\xi-\eta)}{\langle \xi-\eta\rangle_{n}}  }{\eta}  }{s}\\
 =&\inte{0}{t}{  \inte{}{}{ \frac{|\xi|^{\frac{3}{2}}}{|\eta|^{\frac{3}{2}}} s\frac{\xi}{\cro{\xi}_{p}} (1-\theta_{s^\delta}(|\eta|))e^{-is\phi}|\eta|^{\frac{3}{2}}\frac{\hat{f}_{m}(\eta)}{\langle \eta\rangle_{m}}\frac{\hat{f}_{n}(\xi-\eta)}{\langle \xi-\eta\rangle_{n}}  }{\eta}  }{s}\\
 &-\inte{0}{t}{  \inte{}{}{ \frac{|\xi|^{\frac{3}{2}}}{|\eta|^{\frac{3}{2}}} s\frac{\xi-\eta}{\cro{\xi-\eta}_{n}} (1-\theta_{s^\delta}(|\eta|))e^{-is\phi}|\eta|^{\frac{3}{2}}\frac{\hat{f}_{m}(\eta)}{\langle \eta\rangle_{m}}\frac{\hat{f}_{n}(\xi-\eta)}{\langle \xi-\eta\rangle_{n}}  }{\eta}  }{s}\\
 =&I^{1,1}+I^{1,2}.
\end{align*}
Let us focus only on $I^{1,1}$, $I^{1,2}$ being very similar.\\
The bilinear multiplier associated to the symbol $(\eta,\zeta)\mapsto \frac{|\eta+\zeta|^{\frac{3}{2}}}{|\eta|^{\frac{3}{2}}}\chi(\eta,\zeta)$ satisfies H\"{o}lder-type estimates thanks to Lemma \ref{CMex}:
\begin{align*}
 \norlp{I^{1,1}}{2}{}&\lesssim \frac{1}{\sqrt{\m\n}} \inte{0}{t}{ s\norlp{\frac{D}{\cro{D}_{p}}|D|^{3/2}(1-\theta_{s^\delta}(|D|))f_{m}}{2}{} \norlp{e^{-is\langle iD\rangle_m}f_{n}}{\infty}{}}{s}.
\end{align*}
Thanks to Proposition \ref{propmultgeneral}, inequality (\ref{propmult2}), we know that the multiplier $\frac{|D|}{\cro{D}_{p}}$ is bounded in~$L^2$. Hence we can write:
\begin{align*}
 \norlp{I^{1,1}}{2}{}&\lesssim \frac{1}{\sqrt{\m\n}} \inte{0}{t}{ s\norlp{|D|^{3/2}(1-\theta_{s^\delta}(|D|))f_{m}}{2}{} \norlp{e^{-is\langle iD\rangle_n}f_{n}}{\infty}{}}{s}.
\end{align*}
Then Proposition \ref{disp+mult} implies
\begin{align*}
\norlp{e^{-is\langle iD\rangle_n}f_{n}}{\infty}{} \lesssim \n^{\frac{1}{4}}\cro{s}^{-\frac{1}{4}}\sqrt{\nor{f_n}{H^{N}}\nor{f_n}{B_{s}}}.
\end{align*}
Thus
\begin{equation*}
 \norlp{I^{1,1}}{2}{}\lesssim \frac{\n^{\frac{1}{4}}}{\sqrt{\m\n}}\inte{0}{t}{ \cro{s}^{\frac{3}{4}} \norlp{|D|^{3/2}(1-\theta_{s^\delta}(|D|))f_m}{2}{}\sqrt{\nor{f_n}{H^{N}}\nor{f_n}{B_{s}}} }{s}.
\end{equation*}
Then use the high frequencies proposition \ref{hf} to write
\begin{equation*}
 \norlp{I^{1,1}}{2}{}\lesssim  \frac{\n^{\frac{1}{4}}}{\sqrt{\m\n}}\inte{0}{t}{\cro{s}^{\frac{3}{4}}\frac{1}{\cro{s}^{\delta(N-\frac{3}{2})}}\norh{f_{m}}{N} \sqrt{\nor{f_n}{H^{N}}\nor{f_n}{B_{s}}}}{s}
\end{equation*}
Finally, in order to estimate the norm in the space $B_s$, we have to divide by $\sqrt{\cro{t}}$. Then
\begin{equation*}
 \frac{1}{\sqrt{\cro{t}}}\norlp{I^{1,1}}{2}{}\lesssim \inte{0}{t}{ \frac{\max(\m,\n)^{\frac{1}{4}}}{\sqrt{\m\n}}\cro{s}^{-\frac{1}{4}}\norh{f_{m}}{N} \sqrt{\nor{f_n}{H^{N}}\nor{f_n}{B_{s}}}}{s},
\end{equation*}
whenever $N>\frac{1}{\delta}+\frac{3}{2}$, which is true thanks to Condition $(\ref{cond-N})$. This ends the proof of Lemma \ref{contrIhf-lemma}.
\end{prf}
This also ends the proof of Proposition \ref{contrIhf}.
\end{pr}
\section{High Hermite modes estimates}\label{I-HHe}
\begin{prop}\label{contrIhm}
For all $M,N$ satisfying (\ref{cond-M})-(\ref{cond-N}), there exists $(a_p(s))_{p\in\N}$ in $\ell^2$, $\nor{(a_p(s))_{p\in\N}}{\ell^2}\leq 1$ such that
\begin{align*}
&\frac{1}{\sqrt{\cro{t}}}\p^M\sum_{m\in\N,n\geq \cro{t}^{\delta}}\mathcal{M}(m,n,p)\norlp{I^{lf}_{m,n}}{2}{} + \frac{1}{\sqrt{\cro{t}}}\p^M\sum_{m\geq \cro{t}^{\delta},n\in\N}\mathcal{M}(m,n,p)\norlp{I^{lf}_{m,n}}{2}{}\\
& \lesssim \inte{0}{t}{\cro{s}^{-\frac{1}{4}}\nor{f}{S^{M,N}_{s}}^2 a_{p}(s)}{s} .
\end{align*}
\end{prop}
We are going to skip the proof of this result, the idea being similar to the one of the previous section: high regularity leads to a decay of high frequencies (details to be found in \cite{These}).

\section{Low frequencies and low Hermite modes estimates}\label{I-LFHe}
Our aim is to prove the following proposition:
\begin{prop}\label{contrIlfm}
If $M$ and $N$ satisfy Conditions (\ref{cond-M})-(\ref{cond-N}), there exists $(a_p(s))_{p\in\N}$, $(b_p(s))_{p\in\N}$ and $(c_{p}(s))$ in $\ell^2$, with norm less than or equal to $1$, such that
\begin{align*}
\p^M\sum_{m,n\leq \cro{t}^{\delta}}\frac{1}{\sqrt{\cro{t}}}\norlp{I^{lf}_{m,n}}{2}{} \lesssim \inte{0}{t}{\Big(Q(s) a_p(s) + C(s) b_p(s)\Big)}{s} + c_{p}(t)R(t),
\end{align*}
with 
\begin{align*}
Q(s)&:=\cro{s}^{4\delta-\frac{1}{4}}\nor{f}{S^{M,N}_{s}}^2,\\
C(s)&:=\cro{s}^{3\delta+\frac{1}{2}}\nor{f}{S^{M,N}_{s}}^3,\\
R(t)&:=\cro{t}^{\frac{5\delta}{2}+\frac{1}{4}}\left(\nor{f(t)}{S^{M,N}_{s}}^2+\nor{f(1)}{S^{M,N}_{1}}^2\right).
\end{align*}

\end{prop}

\subsection{Intermediate result and summation}

In this section we are going to prove the following proposition:
\begin{prop}\label{intermediate-result-I}
For all $m$ and $n$ integers, and $N\geq \frac{3}{2}$,
\begin{align*}
\frac{1}{\sqrt{t}}\norlp{I^{lf}_{m,n}}{2}{}\lesssim_{m\leftrightarrow n} &\int_{0}^{t}\bigg(\frac{\max(\m,\n)^{3+\frac{3}{4}}}{\sqrt{\m\n}}\cro{s}^{3\delta}\cro{s}^{\frac{1}{4}}\nor{f_m(s)}{H^N}\sqnor{f_n}(s)\\
&+\frac{\max(n,m)^{\frac{1}{2}}}{\sqrt{\m\n}}\cro{s}^{3\delta}\cro{s}^{\frac{3}{2}}\nor{f_m(s)}{H^N}\sqnor{f_m}(s)\sqnor{f_n}(s)\bigg) ds\\
&+(\sqrt{n+1}+\sqrt{m+1})^{2} \frac{\max(n,m)^{\frac{3}{4}}}{\sqrt{\m\n}}t^{\frac{5\delta}{2}+\frac{1}{4}}\left(A(t)+A(1)\right),
\end{align*}
with $A(t)=\norh{f_m(t)}{N}{}\sqnor{f_n}(s)$.
\end{prop}
Going from Proposition \ref{intermediate-result-I} to Proposition \ref{contrIlfm} is quite easy: by Remark \ref{rmqnorm}, we know that 

\begin{enumerate}
\item the quadratic term can be bounded as follows:
\begin{align*}
\nor{f_m(s)}{H^N}\sqnor{f_n}(s)\lesssim \m^{-2M}a_m(s)^2 \n^{-M}b_n(s) \nor{f}{S^{M,N}_{s}}^2,
\end{align*}
with $(a_m(s))_{m\in\N}$ and $(b_n(s))_{n\in\N}$ in the unit ball of $\ell^2$.
\item the cubic term leads to a better bound:
\begin{align*}
\nor{f_m(s)}{H^N}\sqnor{f_m}(s)\sqnor{f_n}(s)\lesssim \m^{-M}a_m(s) \n^{-M}b_n(s) \nor{f}{S^{M,N}_{s}}^3,
\end{align*}
with $(a_m(s))_{m\in\N}$ and $(b_n(s))_{n\in\N}$ in the unit ball of $\ell^2$. Hence
\begin{align*}
\nor{f_m(s)}{H^N}\sqnor{f_m}(s)\sqnor{f_n}(s)\lesssim \m^{-2M}a_m(s) \n^{-M}b_n(s) \nor{f}{S^{M,N}_{s}}^3,
\end{align*}
with $(a_m(s))_{m\in\N}$ and $(b_n(s))_{n\in\N}$ in  the unit ball of $\ell^2$.
\item a same bound can be found for the remaining term: 
\begin{equation*}
A(t)\lesssim \m^{-2M}a_m(s)^2 \n^{-M}b_n(s) \nor{f}{S^{M,N}_{s}}^2,
\end{equation*} with $(a_m(s))_{m\in\N}$ and $(b_n(s))_{n\in\N}$ in the unit ball of $\ell^2$.
\end{enumerate}

We fit in the framework of the bounded resummation theorem \ref{resummation}-(\ref{bddresummation}), and Proposition \ref{contrIlfm} is proven.\\
\\
The proof of Proposition \ref{intermediate-result-I} can be summed up as follows:
\begin{enumerate}
\item first of all we say a word about small times (section \ref{I-LFHE-st}).
\item then, for large times, we have to split the space in two zones: around the space resonant set and far from it. But the way of dealing with those zones will depend on whether Condition (\ref{cond-res}) (page \pageref{cond-res}) is satisfied or not.
\begin{enumerate}
\item if this condition is satisfied,
\begin{enumerate}
\item near the space resonant set, we are going to take advantage of the narrowness of the zone (Section \ref{condC-around}).
\item outside the space resonant set, we are going to perform an integration by parts in $\eta$ and gain some powers of time (Section \ref{condC-outside}).
\end{enumerate}
\item if the condition is not satisfied,
\begin{enumerate}
\item near the space resonant set, we are going to take advantage of the non cancellation of the phase and perform an integration by parts in time (Section \ref{notcondC-around}).
\item outside the space resonant set, we are going to perform an integration by parts in $\eta$ and gain some powers of time (Section \ref{notcondC-outside}).
\end{enumerate}
\end{enumerate}
\end{enumerate}
\subsection{Small times}\label{I-LFHE-st}
Establishing contraction estimates for small times is not a big matter when we study weakly nonlinear dispersive equations. In our situation, we have the following theorem:
\begin{prop}\label{smalltimes}
For all $0<t<1$,
\begin{equation}\label{smalltimes-ineg}
\begin{aligned}
\frac{1}{\sqrt{\cro{t}}}\norlp{I^{lf}_{m,n}(t)}{2}{}\lesssim_{m\leftrightarrow n} \frac{\max(\m,\n)^{\frac{1}{4}}}{\sqrt{\m\n}}\int_{0}^{t}  s^{\frac{3}{4}}  \norh{f_m}{N}{} \sqnor{f_n}(s) ds.
\end{aligned}
\end{equation}
\end{prop}
Let us remark that since we are considering times smaller than $1$, we have $s^{\frac{3}{4}}\leq \cro{s}^{-\frac{1}{4}}$: this is the reason why Proposition \ref{smalltimes} will imply Proposition \ref{intermediate-result-I}. The elementary proof of this result will be skipped, but can be found in \cite{These}.
\subsection{Large times. Estimates for $I^{lf}_{m,n}$ with condition (\ref{cond-res}) satisfied}\label{I-LFHE-Lt-C}
We are going to prove the following result:

\begin{lem}\label{contrIlflm-cp}
For all $m$ and $n$ integers such that condition (\ref{cond-res}) is satisfied, for all $M$, $N$ satisfying (\ref{cond-M})-(\ref{cond-N}) and for all $t\geq 1$,
\begin{align}
\frac{1}{\sqrt{t}}\norlp{I^{lf}_{m,n}}{2}{}\lesssim_{m\leftrightarrow n} \frac{\max(n,m)^{2}}{\sqrt{\m\n}}\inte{1}{t}{\cro{s}^{3\delta}\cro{s}^{-\frac{1}{4}}\nor{f_m}{H^{N}}\mathcal{B}(f_n)(s)}{s}.
\end{align}
\end{lem}
Recall that the study of resonances in Appendix \ref{resonances-asympt} (Lemma \ref{dxi<deta}) implies that if condition (\ref{cond-res}) is satisfied, then the space-resonant set is the space-time-resonant set, and is the straight line $\left\{\xi=\Lambda_{m,n}\eta\right\}$, where 
\begin{equation}
\Lambda_{m,n}=1+\sqrt{\frac{n+1}{m+1}}.
\end{equation} Hence it seems natural to distinguish two zones: close to this set and far away from it.\\
Let $\varphi^{s}(\xi,\eta):=\theta(\sqrt{\cro{s}}\partial_{\eta}\phi(\xi,\eta))$ where $\theta$ is equal to $1$ around $0$: $\varphi^{s}$ localizes in the zone 
\begin{equation*}
-\frac{1}{\sqrt{\cro{s}}}\leq \partial_{\eta}\phi(\xi,\eta)\leq \frac{1}{\sqrt{\cro{s}}}.
\end{equation*} Let us now write $I^{lf}_{m,n}=I^{lf,r}_{m,n}+I^{lf,nr}_{m,n}$, where
\begin{itemize}
 \item $I^{lf,r}_{m,n}$ is the low-frequency, resonant term.
\begin{equation}\label{formule-Ilfr}
I^{lf,r}_{m,n}:=|\xi|^{\frac{3}{2}}\inte{1}{t}{  \inte{}{}{  s\partial_{\xi}\phi \theta_{s^\delta}(|\eta|)\theta_{s^\delta}(|\xi-\eta|)\varphi^{s}(\xi,\eta) e^{-is\phi}\frac{\hat{f}_{m}(\eta)}{\langle \eta\rangle_{m}}\frac{\hat{f}_{n}(\xi-\eta)}{\langle \xi-\eta\rangle_{n}}  }{\eta}  }{s} 
\end{equation}
 \item $I^{lf,nr}_{m,n}$ is the low-frequency, non resonant term.
 
\begin{equation}\label{formule-Ilfnr}
 I^{lf,nr}_{m,n}:=|\xi|^{\frac{3}{2}}\inte{1}{t}{  \inte{}{}{  s\partial_{\xi}\phi \theta_{s^\delta}(|\eta|)\theta_{s^\delta}(|\xi-\eta|) (1-\varphi^{s}(\xi,\eta))e^{-is\phi}\frac{\hat{f}_{m}(\eta)}{\langle \eta\rangle_{m}}\frac{\hat{f}_{n}(\xi-\eta)}{\langle \xi-\eta\rangle_{n}}  }{\eta}  }{s} 
\end{equation}
\end{itemize}

\subsubsection{Around the space resonant set,}\label{condC-around} 
We are going to use the narrowness of the zone where we are localizing in order to prove the following result (which implies Lemma \ref{contrIlflm-cp}).
\begin{lem}\label{lem-condC-around}
For all $m$ and $n$ integers, for all $M$ satisfying (\ref{cond-M}) and $N$ satisfying (\ref{cond-N}), 
\begin{align}
\frac{1}{\sqrt{\cro{t}}}\norlp{I^{lf,r}_{m,n}}{2}{}\lesssim_{m\leftrightarrow n} \inte{1}{t}{\cro{s}^{3\delta}\frac{\max(n,m)^{\frac{3}{4}}}{\sqrt{\m\n}}\cro{s}^{-\frac{1}{4}}\nor{f_m}{H^{N}}\mathcal{B}(f_n)(s)}{s}.
\end{align}
\end{lem}

\begin{pr}
First of all, we use the fact that in the zone where $\theta_{s^\delta}(|\eta|)\theta_{s^\delta}(|\xi-\eta|)\neq 0$, we have $|\eta|\lesssim t^{\delta}$ and $|\xi-\eta|\lesssim t^\delta$, so $|\xi|\lesssim t^\delta$, hence:
\begin{align*}
\norlp{I^{lf,r}_{m,n}}{2}{}\lesssim t^{\frac{3\delta}{2}}\norlp{\inte{1}{t}{  \inte{}{}{ s\partial_{\xi}\phi \theta_{s^\delta}(|\eta|)\theta_{s^\delta}(|\xi-\eta|)\varphi^{s} (\xi,\eta)e^{-is\phi}\frac{\hat{f}_{m}(\eta)}{\langle \eta\rangle_{m}}\frac{\hat{f}_{n}(\xi-\eta)}{\langle \xi-\eta\rangle_{n}}  }{\eta}  }{s}}{2}{}.
\end{align*}

Let us write 
\begin{equation}
S(\xi,\eta):=\sqrt{s}\partial_{\xi}\phi \theta_{s^\delta}(|\eta|)\theta_{s^\delta}(|\xi-\eta|)\varphi^{s} (\xi,\eta).
\end{equation}
Our aim is to get nice Hölder-like estimates for the symbol $S$, which is a bilinear multiplier localizing in a narrow zone, around a curve: this is the point of the paper of Bernicot and Germain \cite{BG2011}, and in particular Theorem \ref{thmger} and its refined version \ref{estbilinband}.\\

\textbf{1. Study of the symbol $S$}\\
Here we are interested in size and width of the support of $S$ and in derivative estimates. We have the following lemma:
\begin{lem}\label{Ill-lt-around-lemS}
The symbol $S$ satisfies the following properties:
\begin{itemize}
\item $S$ is supported in a ball of radius $\rho=s^{\delta}$,
\item $S$ is supported in a band of width $\omega=\frac{s^{\frac{3\delta}{2}}}{(2m+2)\sqrt{s}}$,
\item the derivatives of $S$ satisfy the following inequalities:\begin{align*}
\left|\partial^{a}_{\xi}\partial^{b}_{\eta}S\right|\lesssim \sqrt{s}^{a+b}.
\end{align*}
\end{itemize}
 Hence the symbol $S$ satisfies the hypotheses of Theorem \ref{estbilinband}, with $\rho=s^{\delta}$, $\omega=\frac{s^{3\delta}\sqrt{\max(\m,\n)}}{\sqrt{s}}$, $\mu=\frac{1}{\sqrt{s}}$.
 \end{lem}
\begin{prf}{Lemma}{\ref{Ill-lt-around-lemS}}
\begin{itemize}
\item It is straightforward that $S$ is supported in a ball of radius $s^{\delta}$.
\item Then, we have to determine the width of the support of $S$, that is to say the width of the zone $|\partial_{\eta}\phi|\leq\frac{1}{\sqrt{s}}$. So as to do this, use the asymptotics for $\partial_{\eta}\phi$ computed in Appendix \ref{resonances-asympt}.
\begin{enumerate}
\item In the zone $|\eta|\ll \sqrt{\m}$, (\ref{width-eta-small}) applies, and the witdh of this zone is bounded by
\begin{align*}
\frac{\sqrt{\min(m,n)}}{\sqrt{s}}.
\end{align*}
\item In the zone $|\eta| \geq c  \sqrt{\m}$ ($c\in\R$), since $\eta^2\lesssim s^{\delta}\ll \sqrt{s}$, we are in the asymptotics of (\ref{width-rho-small}). 
Hence the width of the band $|\partial_{\eta}\phi|\leq\frac{1}{\sqrt{s}}$ is less than $\frac{s^{\frac{3\delta}{2}}}{(2m+2)\sqrt{s}}$.
\end{enumerate}
This completes the proof.
\item Finally, we have to estimate the derivatives of $S$.
Thanks to Lemma \ref{dxi<deta}, we know that on the band $|\partial_{\eta}\phi|\leq\frac{1}{\sqrt{s}}$, $|\partial_{\xi}\phi|\leq\frac{1}{\sqrt{t}}$. Hence the inequality is satisfied for $a=b=0$.\\
Then we have to study $\partial^{a}_{\xi}\partial^{b}_{\eta}\left(\sqrt{s}\partial_{\xi}\phi \theta_{s^\delta}(|\eta|)\theta_{s^\delta}(|\xi-\eta|)\varphi^{s}\right)$:
\begin{itemize}
\item any derivative of $\xi$ and $\eta$ of $\partial_{\eta}\phi$ or $\partial_{\xi}\phi$ is a sum of fractions of negative order in $\xi$, in $\eta$, and in $m$, $n$ and $p$. As an example, we have $\partial_{\eta}(\partial_{\eta}\phi)=\frac{2m+2}{\cro{\eta}_{m}^3}+\frac{2n+2}{\cro{\xi-\eta}_{n}^3}$, $\partial_{\xi}(\partial_{\eta}\phi)=-\frac{2n+2}{\cro{\xi-\eta}_{n}^3}$. Then 
\begin{align*}
\left|\partial^{a}_{\xi}\partial^{b}_{\eta}\left(\partial_{\eta}\phi\right)\right|\lesssim 1,\\
\left|\partial^{a}_{\eta}\partial^{b}_{\eta}\left(\partial_{\eta}\phi\right)\right|\lesssim 1.
\end{align*}
\item then $|\partial_{\xi}^{a}\partial_{\eta}^{b}\theta_{s^\delta}(|\eta|)\theta_{s^\delta}(|\xi-\eta|)| \lesssim\sqrt{s}^{a+b}$.
\item finally, 
\begin{align*}
\partial_{\eta}\varphi^{s}=\sqrt{t}\partial_{\eta}(\partial_{\eta}\phi)\theta'\left(\sqrt{s}\partial_{\eta}\phi\right)\\
\partial_{\xi}\varphi^{s}=\sqrt{s}\partial_{\xi}(\partial_{\xi}\phi)\theta'\left(\sqrt{t}\partial_{\xi}\phi\right).
\end{align*}
We just proved that the different derivatives of $\partial_{\eta}\phi$ were bounded by a universal constant. This leads to $|\partial_{\xi}^{a}\partial_{\eta}^{b}\varphi^{s}|\lesssim\sqrt{s}^{a+b}$.
\end{itemize}
Leibniz' rule concludes the proof.
\end{itemize}
Lemma \ref{Ill-lt-around-lemS} is now proved.
\end{prf}

\textbf{2. Estimates.}

The term to estimate is 
\begin{align*}
&\norlp{\inte{1}{t}{  \inte{}{}{ s\partial_{\xi}\phi \theta_{s^\delta}(|\eta|)\theta_{s^\delta}(|\xi-\eta|)\varphi^{s} (\xi,\eta)e^{-is\phi}\frac{\hat{f}_{m}(\eta)}{\langle \eta\rangle_{m}}\frac{\hat{f}_{n}(\xi-\eta)}{\langle \xi-\eta\rangle_{n}}  }{\eta}  }{s}}{2}{}\\
&\lesssim 
\inte{1}{t}{  \norlp{\inte{}{}{ e^{is\cro{\xi}_{p}}s\partial_{\xi}\phi \theta_{s^\delta}(|\eta|)\theta_{s^\delta}(|\xi-\eta|)\varphi^{s} (\xi,\eta)e^{is\cro{\eta}_{m}}\frac{\hat{f}_{m}(\eta)}{\langle \eta\rangle_{m}}e^{is\cro{\xi-\eta}_{n}}\frac{\hat{f}_{n}(\xi-\eta)}{\langle \xi-\eta\rangle_{n}}  }{\eta}  }{2}{}}{s}\\
&\lesssim\inte{1}{t}{\sqrt{s}\norlp{T_{S}\left(e^{-is\cro{D}_{m}}\frac{f_m}{\cro{D}_{m}},e^{-is\cro{D}_{n}}\frac{f_n}{\cro{D}_{n}}\right)}{2}{}}{s},
\end{align*}
with $T_S$ the bilinear Fourier multiplier associated to $S$ as defined in (\ref{def-bilin-mult}). Lemma \ref{Ill-lt-around-lemS} and Theorem \ref{estbilinband} lead to the following estimate.

\begin{align*}
\norlp{T_{S}\left(e^{-is\cro{D}_{m}}\frac{f_m}{\cro{D}_{m}},e^{-is\cro{D}_{n}}\frac{f_n}{\cro{D}_{n}}\right)}{2}{}
&\lesssim \sqrt{\max(\m,\n)}s^{\frac{3\delta}{2}}\norlp{e^{-is\cro{D}_{m}}\frac{f_m}{\cro{D}_{m}}}{2}{}\norlp{e^{-is\cro{D}_{n}}\frac{f_n}{\cro{D}_{n}}}{\infty}{}\!\!\!\!\!\!\!.
\end{align*}
Now by the modified dispersion proposition \ref{disp+mult}, we get the following inequality.
\begin{align*}
\norlp{T_{S}\left(e^{-is\cro{D}_{m}}\frac{f_m}{\cro{D}_{m}},e^{-is\cro{D}_{n}}\frac{f_n}{\cro{D}_{n}}\right)}{2}{}&\lesssim \sqrt{\max(\m,\n)}s^{\frac{3\delta}{2}}\norlp{\frac{f_m}{\cro{D}_{m}}}{2}{}\n^{-\frac{1}{4}}s^{-\frac{1}{4}}\mathcal{B}(f_n)(s).
\end{align*}
Then, thanks to the linear multiplier inequality (\ref{propmult}), we obtain the final inequality.
\begin{align*}
\norlp{T_{S}\left(e^{-is\cro{D}_{m}}\frac{f_m}{\cro{D}_{m}},e^{-is\cro{D}_{n}}\frac{f_n}{\cro{D}_{n}}\right)}{2}{}&\lesssim  \frac{\sqrt{\max(\m,\n)}}{\sqrt{\m}\n^{1/4}}s^{\frac{3\delta}{2}-\frac{1}{4}}\nor{f_m}{H^{N}}\mathcal{B}(f_n)(s).
\end{align*}

Now it remains to integrate over $s$, and divide by $\sqrt{t}$ to get the $B$ norm: we obtain
\begin{align*}
\frac{1}{\sqrt{t}}\norlp{I^{lf,r}_{m,n}}{2}{}
&\lesssim \frac{1}{\sqrt{t}}\inte{1}{t}{\sqrt{s}\sqrt{\max(\m,\n)}\frac{1}{\sqrt{\m}\n^{1/4}}s^{\frac{3\delta}{2}-\frac{1}{4}}\nor{f_m}{H^{N}}\mathcal{B}(f_n)(s)}{s},
\end{align*}
which proves Lemma \ref{lem-condC-around}.\end{pr}

\subsubsection{Outside the space resonant set,}\label{condC-outside}
We have to take advantage fo the non-cancellation of $\partial_{\eta}\phi$: we are going to prove the following result.
\begin{lem}\label{lem-condC-outside}
For all $m$ and $n$ integers, $N\geq 3/2$, $t\geq 1$,
\begin{align}\label{ineq-condC-outside}
\frac{1}{\sqrt{\cro{t}}}\norlp{I^{lf,nr}_{m,n}}{2}{}\lesssim_{m\leftrightarrow n} \frac{\max(n,m)^{2}}{\sqrt{\m\n}}\inte{1}{t}{\cro{s}^{3\delta}\cro{s}^{-\frac{1}{4}}\nor{f_m}{H^{N}}\mathcal{B}(f_n)(s)}{s}.
\end{align}
\end{lem}
\begin{pr}
\paragraph{\textbf{1. Space resonances method.}} Write 
\begin{equation*}
 e^{-is\phi}=\frac{i}{s}\frac{1}{\partial_{\eta}\phi}\partial_{\eta}\left(e^{-is\phi}\right).
\end{equation*}

The term $I^{lf,nr}_{m,n}$ can be rewritten as follows, for $t>1$.

\begin{align*}
I^{lf,nr}_{m,n}=
|\xi|^{\frac{3}{2}}\inte{1}{t}{  \inte{}{}{  s\partial_{\xi}\phi \theta_{s^\delta}(|\eta|)\theta_{s^\delta}(|\xi-\eta|) (1-\varphi^{s})\frac{i}{s}\frac{1}{\partial_{\eta}\phi}\partial_{\eta}\left(e^{-is\phi}\right)\frac{\hat{f}_{m}(\eta)}{\langle \eta\rangle_{m}}\frac{\hat{f}_{n}(\xi-\eta)}{\langle \xi-\eta\rangle_{n}}  }{\eta}  }{s}.
\end{align*}
An integration by parts in $\eta$ leads to
\begin{equation*}
 I^{lf,nr}_{m,n} = \mathcal{I}^1+ \mathcal{I}^2+ \mathcal{I}^3+ \mathcal{I}^4,
\end{equation*}
where
\begin{align*}
 \mathcal{I}^1&:=|\xi|^{\frac{3}{2}}\inte{1}{t}{  \inte{}{}{   \theta_{s^\delta}(|\eta|)\theta_{s^\delta}(|\xi-\eta|) (1-\varphi^{s})i\frac{\partial_{\xi}\phi}{\partial_{\eta}\phi}\frac{\partial_{\eta}^2\phi}{\partial_{\eta}\phi}e^{-is\phi}\frac{\hat{f}_{m}(\eta)}{\langle \eta\rangle_{m}}\frac{\hat{f}_{n}(\xi-\eta)}{\langle \xi-\eta\rangle_{n}}  }{\eta}  }{s}, \\
 \mathcal{I}^2&:= -|\xi|^{\frac{3}{2}}\inte{1}{t}{  \inte{}{}{   \theta_{s^\delta}(|\eta|)\theta_{s^\delta}(|\xi-\eta|) (1-\varphi^{s})i\frac{\partial_{\xi}\phi}{\partial_{\eta}\phi}e^{-is\phi}\partial_{\eta}\left(\frac{\hat{f}_{m}(\eta)}{\langle \eta\rangle_{m}}\right)\frac{\hat{f}_{n}(\xi-\eta)}{\langle \xi-\eta\rangle_{n}}  }{\eta}  }{s}, \\
 \mathcal{I}^3&:=-|\xi|^{\frac{3}{2}}\inte{1}{t}{  \inte{}{}{   \theta_{s^\delta}(|\eta|)\theta_{s^\delta}(|\xi-\eta|) (1-\varphi^{s})i\frac{\partial_{\xi}\phi}{\partial_{\eta}\phi}e^{-is\phi}\frac{\hat{f}_{m}(\eta)}{\langle \eta\rangle_{m}}\partial_{\eta}\left(\frac{\hat{f}_{n}(\xi-\eta)}{\langle \xi-\eta\rangle_{n}}\right)  }{\eta}  }{s}, \\
 \mathcal{I}^4&:=-|\xi|^{\frac{3}{2}}\inte{1}{t}{  \inte{}{}{   \partial_{\eta}\Big( \theta_{s^\delta}(|\eta|)\theta_{s^\delta}(|\xi-\eta|) (1-\varphi^{s})\partial_{\xi}\phi\Big)i\frac{1}{\partial_{\eta}\phi}e^{-is\phi}\frac{\hat{f}_{m}(\eta)}{\langle \eta\rangle_{m}}\frac{\hat{f}_{n}(\xi-\eta)}{\langle \xi-\eta\rangle_{n}}  }{\eta}  }{s}. 
\end{align*}
We will only prove
\begin{align*}
\frac{1}{\sqrt{\cro{t}}}\norlp{\mathcal{I}^1}{2}{}\lesssim \frac{\max(n,m)^{2}}{\sqrt{\m\n}}\inte{1}{t}{\cro{s}^{3\delta}\cro{s}^{-\frac{1}{4}}\nor{f_m}{H^{N}}\mathcal{B}(f_n)(s)}{s},
\end{align*}
the inequalities for $\mathcal{I}^j$, $j=2,3,4$ being treated similarly. The term $\mathcal{I}^{1}$ is actually the harder to deal with since the fraction $\frac{\partial_{\eta}^{2}\phi}{\partial_{\eta}\phi}$ gets big close to the space-time resonant zone. The terms of the form $\partial_{\eta}f$ appearing in the terms $\mathcal{I}^j$, $j=2,3,4$ are not really problematic: the difficulty coming from them is compensated by the absence of $\frac{\partial_{\eta}^{2}\phi}{\partial_{\eta}\phi}$.\\
\paragraph{\textbf{2. Estimates for $\mathcal{I}^1$.}} The main problem arising here is to be able to find a bilinear estimate for the symbol
\begin{align}\label{symbol-condC-outside}
S:=\theta_{s^\delta}(|\eta|)\theta_{s^\delta}(|\xi-\eta|)(1-\varphi^{s})\frac{\partial_{\xi}\phi}{\partial_{\eta}\phi}\frac{\partial_{\eta}^2\phi}{\partial_{\eta}\phi}.
\end{align}

This symbol does not enter directly in the framework of the Bernicot-Germain theorem \ref{estbilinband}. In order to understand better the behaviour of this multiplier, we split the frequency space along the level lines of $\partial_{\eta}\phi$. So as to do this cutoff in a smooth way, let us define the following functions.
\begin{defi}
Let $\omega$ be a real function supported in $\left[\frac{1}{2},1\right]$ such that
\begin{align*}
\forall x\neq 0,~\sum_{j\in\Z}(\omega(2^{j}x)+\omega(-2^{j}x))=1.
\end{align*} Define the following functions
\begin{align*}
\mathbb{I}_{a\sim b}&=\omega\left(\frac{a}{b}\right),\\
\mathbb{I}_{a\leq m}&=\sum_{2^{k}\leq m} \omega(2^{k}a).
\end{align*}
\end{defi} 
Now we can write
\begin{align}\label{decomposition-phi}
(1-\varphi^{s})\frac{\partial_{\xi}\phi}{\partial_{\eta}\phi}\frac{\partial_{\eta}^2\phi}{\partial_{\eta}\phi}= \sum_{1/2\leq 2^{j}\leq \sqrt{s}}S^{+}_{j}+S^{-}_{j},
\end{align}
where
\begin{align}\label{def-Sj}
S^{\pm}_{j}&=\theta_{s^\delta}(|\eta|)\theta_{s^\delta}(|\xi-\eta|) \mathbb{I}_{\partial_{\eta}\phi\sim \pm 2^{-j}}(1-\varphi^{s})\frac{\partial_{\xi}\phi}{\partial_{\eta}\phi}\frac{\partial_{\eta}^2\phi}{\partial_{\eta}\phi}.
\end{align}

Then, we split dyadically the frequency space: the asymptotics of $\phi$ and its derivatives strongly depend on the comparison between the size of the frequencies and the size of $|\partial_{\eta}\phi|$: in Lemma \ref{asymptphi} we obtain three different asymptotical regimes, depending on a parameter $\varrho(m,j,\eta)$ defined by
\begin{equation}\label{def-varrho}
\varrho(m,j,\eta):=\frac{\eta^2}{2^{j}m}.
\end{equation} To deal with this, we are going to use the smooth functions $\mathbb{I}_{\sqrt{|\xi|^2+|\eta|^2}\sim 2^{k}}$ and $\mathbb{I}_{|\xi|^2+|\eta|^2\leq m}$.
This is why we need to define the following symbol:
\begin{align*}
S^{\pm}_{j,k}=\mathbb{I}_{\partial_{\eta}\phi\sim \pm 2^{-j}}\mathbb{I}_{\sqrt{|\xi|^2+|\eta|^2}\sim 2^{k}}(1-\varphi^{s})\frac{\partial_{\xi}\phi}{\partial_{\eta}\phi}\frac{\partial_{\eta}^2\phi}{\partial_{\eta}\phi}.
\end{align*}
Let us finally rewrite the symbol $S$ defined in $(\ref{symbol-condC-outside})$ in order to take into account the different asymptotics for $\partial_{\eta}\phi$. Write
\begin{align*}
S=M^{1}+M^{2}+M^{3},
\end{align*}
where $M^{1}$, $M^{2}$ and $M^{3}$ are defined as follows.
\begin{enumerate}
\item The symbol $M^{1}$ corresponds to small values of $|\xi|,|\eta|$.
\begin{align*}
M^{1}&= \mathbb{I}_{|\xi|^2+|\eta|^2\leq \frac{\m}{2}}\left(\sum_{1/2\leq 2^{j}\leq \sqrt{s}}S^{+}_{j}+S^{-}_{j}\right).
\end{align*}
\item The symbol $M^{2}$ corresponds to small values of the parameter $\varrho(m,j,2^k)$ defined in (\ref{def-varrho}).
\begin{align*}
M^{2}=(1-\mathbb{I}_{|\xi|^2+|\eta|^2\leq \frac{\m}{2}})\sum_{1/2\leq 2^{j}\leq \sqrt{s}}~~\sum_{k|\varrho(m,j,2^k)\ll 1}S^{\pm}_{j,k}.
\end{align*}
\item Finally the symbol $M^{3}$ corresponds to the remaining terms, i.e. large values of $\varrho(m,j,\eta)$.
\begin{align*}
M^{3}=(1-\mathbb{I}_{|\xi|^2+|\eta|^2\leq \frac{\m}{2}})\sum_{1/2\leq 2^{j}\leq \sqrt{s}}~~\sum_{k|\varrho(m,j,2^k)\gtrsim 1}S^{\pm}_{j,k}.
\end{align*}
\end{enumerate}
If $q\in\{1,2,3\}$, we write $J^{q}$ for
\begin{align*}
J^{q}:&=\inte{1}{t}{  \inte{}{}{  M^{q}(\xi,\eta)e^{-is\phi}\frac{\hat{f}_{m}(\eta)}{\langle \eta\rangle_{m}}\frac{\hat{f}_{n}(\xi-\eta)}{\langle \xi-\eta\rangle_{n}}  }{\eta}  }{s}.
\end{align*}
With these notations, we remark that 
\begin{align*}
\mathcal{I}^{1}=|\xi|^{\frac{3}{2}}(J^{1}+J^{2}+J^{3}).
\end{align*}

\textbf{1. Estimates for $M^{1}$.}

\begin{lem}\label{lem-condC-M1}(Bilinear estimate for low frequencies)
We have the following estimate.
\begin{align}\label{ineq-M1}
\frac{1}{\sqrt{\cro{t}}}\norlp{|\xi|^{\frac{3}{2}}J^{1}}{2}{}&\lesssim_{m\leftrightarrow n} \frac{\n^{\frac{1}{4}}}{\sqrt{\m\n}}\min(\sqrt{\m},\sqrt{\n})\inte{1}{t}{s^{\frac{3}{2}\delta}s^{-\frac{1}{4}}\nor{f_m}{H^{N}}\mathcal{B}(f_n)(s)}{s}.
\end{align}
\end{lem}
\begin{rmq}
Remark that the inequality $(\ref{ineq-M1})$ is stronger than the inequality~$(\ref{ineq-condC-outside})$ in Lemma~\ref{lem-condC-outside}.
\end{rmq}
\begin{pr}
Let $J^{1}_{j}$ be the following quantity.
\begin{align*}
J^{1,\pm}_{j}:=\inte{}{}{  \mathbb{I}_{|\xi|^2+|\eta|^2\leq \frac{\m}{2}}S_{j}^{\pm}e^{-is\phi}\frac{\hat{f}_{m}(\eta)}{\langle \eta\rangle_{m}}\frac{\hat{f}_{n}(\xi-\eta)}{\langle \xi-\eta\rangle_{n}}  }{\eta},
\end{align*}
where $S^{\pm}_{j}$ is defined in $(\ref{def-Sj})$.\\

First we will establish a bilinear estimate for the symbol $\mathbb{I}_{|\xi|^2+|\eta|^2\leq \frac{\m}{2}}S_{j}^{-}$ which will be denoted $\tld{S}_{j}$ for the sake of simplicity. Given the central symmetry for the level lines of $\partial_{\eta}\phi$, estimates for $\tld{S}_{j}^{-}$ will also be valid for $\tld{S}_{j}^{+}$.\\

Let us adopt the following strategy: first we study completely the symbol $\tld{S}_{j}$, then we rescale it to fit in the Bernicot-Germain theorem's hypotheses.\\
\begin{lem}\label{Ill-lt-outside-lemtldSj}
The symbol $\tld{S}_{j}$ satisfies the following properties:
\begin{itemize}
\item $\tld{S}_{j}$ is supported in a ball of radius $\rho=\sqrt{\frac{\m}{2}}$,
\item $\tld{S}_{j}$ is supported in a band of width $\omega=2^{-j}\min(\sqrt{\m},\sqrt{\n})$,
\item the derivatives of $\tld{S}_{j}$ satisfy the following inequalities: for all $a,b$ integers,\begin{align*}
\left|\partial_{\xi}^{a}\partial_{\eta}^{b}\tld{S}_{j}\right|\lesssim 2^j(2^{j})^{a+b}.
\end{align*}
\end{itemize}
 Hence the symbol $2^{-j}\tld{S}_{j}$ satisfies the hypotheses of Theorem \ref{estbilinband}, with 
${\rho=\sqrt{\frac{\m}{2}}}$, ${\omega=2^{-j}\min(\sqrt{\m},\sqrt{\n})}$, ${\mu=2^{-j}}$.
 \end{lem}
The proof is skipped here and can be found in \cite{These}.\\

 Rewriting $J^{1,\pm}_{j}$ as a bilinear operator gives
 \begin{align}\label{M1-prf-0}
 \norlp{J^{1}_{j}}{2}{}=2^j\norlp{T_{2^{-j}\tld{S}_{j}}\left(e^{-is\cro{D}_{m}}\frac{f_m}{\cro{D}_{m}},e^{-is\cro{D}_{n}}\frac{f_n}{\cro{D}_{n}}\right)}{2}{\xi}.
 \end{align}
 Using Theorem \ref{estbilinband}, we obtain
 \begin{align*}
& \norlp{T_{2^{-j}\tld{S}_{j}}\left(e^{-is\cro{D}_{m}}\frac{f_m}{\cro{D}_{m}},e^{-is\cro{D}_{n}}\frac{f_n}{\cro{D}_{n}}\right)}{2}{\xi}\\
 &\lesssim  \max\left(1,\frac{\omega}{\mu}\right)(\rho\omega)^{\frac{1}{2}+\frac{1}{2}-1} \norlp{e^{-is\cro{D}_{m}}\frac{f_m}{\cro{D}_{m}}}{2}{}\norlp{e^{-is\cro{D}_{n}}\frac{f_n}{\cro{D}_{n}}}{\infty}{}\\
  &\lesssim \min(\sqrt{\m},\sqrt{\n})\norlp{e^{-is\cro{D}_{m}}\frac{f_m}{\cro{D}_{m}}}{2}{}\norlp{e^{-is\cro{D}_{n}}\frac{f_n}{\cro{D}_{n}}}{\infty}{}.
 \end{align*}
  Then use the Fourier multiplier Proposition \ref{propmult}, to get:
 \begin{align}\label{M1-prf-1}
\norlp{e^{-is\cro{D}_{m}}\frac{f_m}{\cro{D}_{m}}}{2}{}&\lesssim \frac{1}{\sqrt{\m}}\norlp{f_m}{2}{}\\
&\lesssim \frac{1}{\sqrt{\m}}\nor{f_m}{H^{N}}.
 \end{align}
 Similarly, Proposition \ref{disp+mult} implies
 \begin{align}\label{M1-prf-2}
 \norlp{e^{-is\cro{D}_{n}}\frac{f_n}{\cro{D}_{n}}}{\infty}{}\lesssim \n^{-\frac{1}{4}}s^{-\frac{1}{4}}\mathcal{B}(f_n)(s).
 \end{align}
 By (\ref{M1-prf-1}) and (\ref{M1-prf-2}), 
 \begin{equation}\label{M1-prf-3}
 \norlp{T_{2^{-j}\tld{S}_{j}}\left(e^{-is\cro{D}_{m}}\frac{f_m}{\cro{D}_{m}},e^{-is\cro{D}_{n}}\frac{f_n}{\cro{D}_{n}}\right)}{2}{\xi} \lesssim \frac{\min(\m,\n)}{\n^{\frac{1}{4}}\sqrt{\m}}\nor{f_m}{H^{N}}\mathcal{B}(f_n)(s).
 \end{equation}
 Now by (\ref{M1-prf-0}), (\ref{M1-prf-3}) and since $|\xi|\leq s^{\delta}$, we get
 \begin{align}\label{M1-prf-4}
 \norlp{|\xi|^{\frac{3}{2}}J^1_{j}}{2}{}\lesssim 2^{j}\frac{\n^{\frac{1}{4}}}{\sqrt{\m\n}}\min(\sqrt{\m},\sqrt{\n})s^{\frac{3}{2}\delta}s^{\frac{1}{4}}\frac{1}{\sqrt{s}}\nor{f_m}{H^{N}}\mathcal{B}(f_n)(s).
 \end{align}
Then, since
\begin{align*}
\frac{1}{\sqrt{t}}\norlp{|\xi|^{\frac{3}{2}}\inte{1}{t}{\sum_{\frac{1}{2}\leq 2^j\leq\sqrt{s}}J_{j}^{1}(s)}{s}}{2}{}\lesssim &\frac{1}{\sqrt{t}}\inte{1}{t}{\sum_{\frac{1}{2}\leq 2^j\leq\sqrt{s}}\norlp{J_{j}^{1}(s)}{2}{}}{s},
\end{align*}
Inequality (\ref{M1-prf-4}) gives
\begin{align*}
&\frac{1}{\sqrt{t}}\norlp{|\xi|^{\frac{3}{2}}\inte{1}{t}{\sum_{\frac{1}{2}\leq 2^j\leq\sqrt{s}}J_{j}^{1}(s)}{s}}{2}{}\\
&\lesssim \frac{\n^{\frac{1}{4}}}{\sqrt{\m\n}}\min(\sqrt{\m},\sqrt{\n}) \inte{1}{t}{\nor{f_m}{H^{N}}\mathcal{B}(f_n)(s)\frac{1}{\sqrt{s}}s^{\frac{3}{2}\delta}s^{\frac{1}{4}}\frac{1}{\sqrt{s}}\sum_{\frac{1}{2}\leq 2^j\leq\sqrt{s}}2^{j}}{s}\\
&\lesssim \frac{\n^{\frac{1}{4}}}{\sqrt{\m\n}}\min(\sqrt{\m},\sqrt{\n})\inte{1}{t}{s^{\frac{3}{2}\delta}s^{-\frac{1}{4}}\nor{f_m}{H^{N}}\mathcal{B}(f_n)(s)}{s},
\end{align*}
which concludes the proof of Lemma \ref{lem-condC-M1}.
\end{pr}

\textbf{2. Estimates for $M^{2}$.}

\begin{lem}\label{lemM2}(Bilinear estimate for high frequencies and small values of $\varrho(m,j,\eta)$)\\
We have the following inequality.
\begin{align}\label{ineq-M2}
\frac{1}{\sqrt{\cro{t}}}\norlp{|\xi|^{\frac{3}{2}}J^{2}}{2}{}&\lesssim \frac{\max(\m,\n)^{\frac{1}{4}}}{\m\sqrt{\m\n}}\inte{1}{t}{s^{\frac{3}{2}\delta}s^{-\frac{1}{4}}\nor{f_m}{H^{N}}\mathcal{B}(f_n)(s)}{s}.
\end{align}
\end{lem}
\begin{rmq}
The inequality $(\ref{ineq-M2})$ is stronger than the inequality~$(\ref{ineq-condC-outside})$ in Lemma~\ref{lem-condC-outside}.
\end{rmq}

\begin{pr}

Recall that:
\begin{equation*}
J^{2}=\inte{1}{t}{ \sum_{\frac{2^{k}}{\sqrt{\m}2^{j}}\ll 1,~2^k\gtrsim\sqrt{\m}}J^{\pm}_{j,k}(s)  }{s}=\inte{1}{t}{ \sum_{\frac{2^{k}}{\sqrt{\m}2^{j}}\ll 1,~2^k\gtrsim\sqrt{\m}} \inte{}{}{S^{\pm}_{j,k}e^{-is\phi}\frac{\hat{f}_{m}(\eta)}{\langle \eta\rangle_{m}}\frac{\hat{f}_{n}(\xi-\eta)}{\langle \xi-\eta\rangle_{n}}}{\eta}}{s}.
\end{equation*}
We start by stating multilinear estimates for the symbol $S_{j,k}=S^{-}_{j,k}$ (the case $S^{+}_{j,k}$ is similar). We skip the proof of the following result, very similar to Lemma \ref{Ill-lt-outside-lemtldSj}.\\
\begin{lem}\label{Ill-lt-outside-J2-lemSjk}
The symbol $S_{j,k}$ satisfies the following properties:
\begin{itemize}
\item $S_{j,k}$ is supported in a ball of radius $\rho=2^{k}$,
\item $S_{j,k}$ is supported in a band of width $\omega=\frac{2^{3k}2^{-j}}{2m+2}$,
\item the derivatives of $S_{j,k}$ satisfy the following inequalities: for all $a,b$ integers,\begin{align*}
\left|\partial_{\xi}^{a}\partial_{\eta}^{b}S_{j,k}\right|\lesssim 2^j(2^{j})^{a+b}.
\end{align*}
\end{itemize}
 Hence the symbol $2^{-j}S_{j,k}$ satisfies the hypotheses of Theorem \ref{estbilinband}, with ${\rho=2^{k}}$, ${\omega=2^{-j}\frac{2^{3k}}{2m+2}}$, ${\mu=2^{-j}}$.
 \end{lem}
 If we rewrite $J_{j,k}$ as a bilinear operator,

 \begin{align*}
 \norlp{J_{j,k}}{2}{}=2^j\norlp{T_{2^{-j}S_{j,k}}\left(e^{-is\cro{D}_{m}}\frac{f_m}{\cro{D}_{m}},e^{-is\cro{D}_{n}}\frac{f_n}{\cro{D}_{n}}\right)}{2}{\xi},
 \end{align*}
we can apply Theorem \ref{estbilinband} to obtain
 \begin{align*}
 &\norlp{T_{2^{-j}S_{j}}\left(e^{-is\cro{D}_{m}}\frac{f_m}{\cro{D}_{m}},e^{-is\cro{D}_{n}}\frac{f_n}{\cro{D}_{n}}\right)}{2}{\xi}\\
 &\lesssim \frac{2^{3k}}{2m+2}\norlp{e^{-is\cro{D}_{m}}\frac{f_m}{\cro{D}_{m}}}{2}{}\norlp{e^{-is\cro{D}_{n}}\frac{f_n}{\cro{D}_{n}}}{\infty}{},
 \end{align*}
 by the dilation lemma \ref{multdilat}. \\
 Then, by the dispersive estimate and Proposition \ref{propmultgeneral}, inequality (\ref{propmult}),
 \begin{align*}
 \norlp{|\xi|^{\frac{3}{2}}J_{j,k}}{2}{}\lesssim 2^{j}\frac{\n^{\frac{1}{4}}}{\sqrt{\m\n}}\frac{2^{3k}}{2m+2}s^{\frac{3}{2}\delta}s^{\frac{1}{4}}\frac{1}{\sqrt{s}}\nor{f_m}{H^{N}}\mathcal{B}(f_n)(s).
 \end{align*}
Now we sum over $k$:
\begin{align*}
\sum_{2^k\leq s^\delta} \norlp{|\xi|^{\frac{3}{2}}J_{j,k}}{2}{}
\lesssim 2^{j}\frac{\n^{\frac{1}{4}}}{m\sqrt{\m\n}}s^{\frac{3}{2}\delta}s^{\frac{1}{4}}\frac{1}{\sqrt{s}}\nor{f_m}{H^{N}}\mathcal{B}(f_n)(s).
\end{align*}
Finally by sum over $j$ and integrating in time,
\begin{align*}
\frac{1}{\sqrt{t}}\norlp{\inte{1}{t}{\sum_{2^k\leq s^\delta}\sum_{\frac{1}{2}\leq 2^j\leq\sqrt{s}}J_{j,k}(s)}{s}}{2}{}\lesssim \frac{\max(\m,\n)^{\frac{1}{4}}}{\m\sqrt{\m\n}}\inte{1}{t}{s^{\frac{3}{2}\delta}s^{-\frac{1}{4}}\nor{f_m}{H^{N}}\mathcal{B}(f_n)(s)}{s}.
\end{align*}
This ends the proof of Lemma \ref{lemM2}.
\end{pr}

\textbf{4. Estimates for $M^{3}$.}
\begin{lem}\label{lemM3}(Bilinear estimate for high frequencies, and small values of $j$)
We have the following inequality.
\begin{align}\label{ineq-M3}
\frac{1}{\sqrt{\cro{t}}}\norlp{|\xi|^{\frac{3}{2}}J^{3}}{2}{}&\lesssim \frac{\max(\m,\n)^{\frac{1}{4}}}{\sqrt{\m\n}}\max\left(\sqrt{\frac{n}{m}},1\right)\inte{1}{t}{s^{3\delta}s^{-\frac{1}{4}}\nor{f_m}{H^{N}}\mathcal{B}(f_n)(s)}{s}.
\end{align}
\begin{rmq}
Inequality $(\ref{ineq-M3})$ is stronger than the inequality $(\ref{ineq-condC-outside})$ in Lemma~\ref{lem-condC-outside}.
\end{rmq}
\end{lem}
\begin{pr}
Recall that:
\begin{equation*}
J^{3}=\inte{1}{t}{ \sum_{\frac{2^{k}}{\sqrt{\m}2^{j}}\gtrsim 1,~2^k\gtrsim\sqrt{\m}}J^{\pm}_{j,k}(s)  }{s}=\inte{1}{t}{ \sum_{\frac{2^{k}}{\sqrt{\m}2^{j}}\gtrsim 1,~2^k\gtrsim\sqrt{\m}} \inte{}{}{S^{\pm}_{j,k}e^{-is\phi}\frac{\hat{f}_{m}(\eta)}{\langle \eta\rangle_{m}}\frac{\hat{f}_{n}(\xi-\eta)}{\langle \xi-\eta\rangle_{n}}}{\eta}}{s}.
\end{equation*}
First we establish multilinear estimates for the symbol $S_{j,k}=S^{-}_{j,k}$ (as previously, the case $S^{+}_{j,k}$ is similar): we are not going to give the proof of the following result but it depends on the asymptotics found in (\ref{etainf}) and (\ref{width-rho-large}).
\begin{lem}\label{Ill-lt-outside-J3-lemSjk}
The symbol $S_{j,k}$ satisfies the following properties:
\begin{itemize}
\item $S_{j,k}$ is supported in a ball of radius $\rho=2^{k}$,
\item $S_{j,k}$ is supported in a band of width $\omega=2^{\frac{j}{2}}\sqrt{2n+2}\lesssim 2^{\frac{j}{2}}\max(\m,\n)^{\frac{1}{2}}$,
\item the derivatives of $S_{j,k}$ satisfy the following inequalities: for all $a,b$ integers,\begin{align*}
\left|\partial_{\xi}^{a}\partial_{\eta}^{b}S_{j,k}\right|\lesssim 2^j(2^{j})^{a+b}.
\end{align*}
\end{itemize}
 Hence the symbol $2^{-j}S_{j,k}$ satisfies the hypotheses of Theorem \ref{estbilinband}, with ${\rho=2^{k}}$, ${\omega=2^{\frac{j}{2}}\sqrt{\max(\m,\n)}}$, ${\mu=2^{-j}}$.
 \end{lem}
  This leads to

 \begin{align*}
 \norlp{|\xi|^{\frac{3}{2}}J^2_{j,k}}{2}{}\lesssim 2^{j}\frac{\max(\m,\n)^{\frac{1}{4}}}{\sqrt{\m\n}}2^{\frac{3j}{2}}\max(\sqrt{\n},\sqrt{\m})\inte{1}{t}{s^{\frac{3}{2}\delta}s^{\frac{1}{4}}\frac{1}{\sqrt{s}}\nor{f_m}{H^{N}}\mathcal{B}(f_n)(s)}{s}.
 \end{align*}
 Here we are in the regime where $2^{j}\lesssim\frac{2^k}{\sqrt{\m}}\lesssim\frac{s^\delta}{\sqrt{\m}}$. Hence the inequality rewrites as follows.
  \begin{align*}
 \norlp{|\xi|^{\frac{3}{2}}J^2_{j,k}}{2}{}\lesssim 2^{j}\frac{\max(\m,\n)^{\frac{1}{4}}}{\sqrt{\m\n}}\max\left(\sqrt{\frac{n}{m}},1\right)\inte{1}{t}{s^{3\delta}s^{\frac{1}{4}}\frac{1}{\sqrt{s}}\nor{f_m}{H^{N}}\mathcal{B}(f_n)(s)}{s}.
 \end{align*}
 Then summing over $j$ and $k$ and integrating leads to
 \begin{align*}
&\frac{1}{\sqrt{t}}\norlp{\inte{1}{t}{\sum_{2^k\leq s^\delta}~~\sum_{\frac{1}{2}\leq 2^j\leq\sqrt{s}}J_{j}^{3}(s)}{s}}{2}{}\lesssim \frac{\max(\m,\n)^{\frac{1}{4}}}{\sqrt{\m\n}}\max\left(\sqrt{\frac{n}{m}},1\right)\inte{1}{t}{s^{3\delta}s^{-\frac{1}{4}}\nor{f_m}{H^{N}}\mathcal{B}(f_n)(s)}{s}.
\end{align*}
 Then Lemma \ref{lemM3} is proved.
\end{pr}
\\
Then gathering Lemma \ref{lem-condC-M1}, Lemma \ref{lemM2} and Lemma \ref{lemM3} leads to Lemma \ref{lem-condC-outside}.\end{pr}

\subsection{Case where $p>m$, $p>n$ but condition (\ref{cond-res}) is not satisfied}\label{I-LFHE-Lt-notC}
We are going to prove the following lemma:
\begin{lem}\label{contrIlflm-notcp}
For all $m$ and $n$ integers satisfying (\ref{cond-M})-(\ref{cond-N}) and such that condition (\ref{cond-res}) is not satisfied, 
\begin{align*}
\frac{1}{\sqrt{t}}\norlp{I^{lf}_{m,n}}{2}{}\lesssim&\frac{\max(\m,\n)^{3+\frac{3}{4}}}{\sqrt{\m\n}}\inte{1}{t}{s^{3\delta}s^{-\frac{1}{4}}\nor{f_m}{H^{N}}\sqnor{f_n}(s)}{s}\\
&+\frac{\max(n,m)^{\frac{1}{2}}}{\sqrt{\m\n}}\inte{1}{t}{s^{3\delta}s^{\frac{3}{2}}\nor{f_m}{H^{N}}\sqnor{f_m}(s)\sqnor{f_n}(s)}{s}\\
&+(\sqrt{n+1}+\sqrt{m+1})^{2} \frac{\max(n,m)^{\frac{3}{4}}}{\sqrt{\m\n}}t^{\frac{5\delta}{2}+\frac{1}{4}}\left(A(t)+A(1)\right),
\end{align*}
with $A(t)=\norh{f_m(t)}{N}{}\sqnor{f_n}(t)$.
\end{lem}

Recall the situation: in the case where $p>m$, $p>n$ but condition (\ref{cond-res}) is not satisfied, there are no space-time resonances. When we are close to the space-resonant straight line, a normal form transformation should help. \\

However, one of the main problems is that we do not have $|\partial_{\xi}\phi|\leq|\partial_{\eta}\phi|$. So as to deal with this new configuration, we are going to loosen the constraint on the narrowness of the zone close to $\Sc$: this will make the estimates outside this zone easier. Inside it, we will be able to use the time-resonances method since $\phi$ does not vanish.\\
We are performing two different cutoffs: 
\begin{itemize}
\item $\theta$ is a compactly supported $\mathcal{C}^{\infty}$ function equal to $1$ on $[-1,1]$
\item $\theta_{s^\delta}(|\eta|)=\theta\left(\frac{|\eta|}{s^{\delta}}\right)$
\item We have to choose a new function $\psi^{s}$ localizing around the space resonant set. Our idea is to take the widest zone which does not meet the space-resonant set. Proposition \ref{timres-width} will be very useful: if $\psi$ localizes in a neigborhood of size $\frac{c}{(\sqrt{n+1}+\sqrt{m+1})^2}\frac{1}{R}$ of $\Sc$, we can be sure that we will not meet the time-resonant set.\\
If we adapt the proof of Lemma \ref{Ill-lt-around-lemS}, we know that the zone $|\partial_{\eta}\phi|\leq d$ is of width $\sqrt{\max(\m,\n)}s^{3\delta}d$. Consequently the function $\psi^{s}$ can be chosen equal to
\begin{align}\label{def-psi-notcondC}
\psi^{s}(\xi,\eta)=\theta\left(c'\sqrt{\max(\m,\n)}s^{3\delta}(\sqrt{n+1}+\sqrt{m+1})^2 s^\delta|\partial_{\eta}\phi|\right).
\end{align}
\end{itemize}
Then write
\begin{equation*}
I^{lf}_{m,n}=I^{lf,r}_{m,n}+I^{lf,nr}_{m,n},
\end{equation*}
where
\begin{itemize}
 \item $I^{lf,r}_{m,n}$ is the space-resonant term.
\begin{equation*}
 I^{lf,r}_{m,n}:=|\xi|^{\frac{3}{2}}\inte{1}{t}{  \inte{}{}{  s\partial_{\xi}\phi \theta_{s^\delta}(|\eta|)\theta_{s^\delta}(|\xi-\eta|)\psi^{s} e^{-is\phi}\frac{\hat{f}_{m}(\eta)}{\langle \eta\rangle_{m}}\frac{\hat{f}_{n}(\xi-\eta)}{\langle \xi-\eta\rangle_{n}}  }{\eta}  }{s} ,
\end{equation*}
 \item $I^{lf,nr}_{m,n}$ is the non space-resonant term: 
\begin{equation*}
I^{lf,nr}_{m,n}:=|\xi|^{\frac{3}{2}}\inte{1}{t}{  \inte{}{}{  s\partial_{\xi}\phi \theta_{s^\delta}(|\eta|)\theta_{s^\delta}(|\xi-\eta|) (1-\psi^{s})e^{-is\phi}\frac{\hat{f}_{m}(\eta)}{\langle \eta\rangle_{m}}\frac{\hat{f}_{n}(\xi-\eta)}{\langle \xi-\eta\rangle_{n}}  }{\eta}  }{s} ,
\end{equation*}
\end{itemize}

\subsubsection{Around the space resonant set ($I^{lf,r}_{m,n}$).}\label{notcondC-around}

\begin{lem}\label{contrnotCp-res}
For all $m$ and $n$ integers, $M$ and $N$ integers satisfying (\ref{cond-M})-(\ref{cond-N}), $t\geq 1$,
\begin{align*}
\frac{1}{\sqrt{t}}\norlp{I^{lf,r}_{m,n}}{2}{}\lesssim \mathcal{B}+ \mathcal{Q}+\mathcal{C},
\end{align*}
where $\mathcal{B}$ is the boundary term:
\begin{align*}
\mathcal{B}:=(\sqrt{n+1}+\sqrt{m+1})^{2} \frac{\max(n,m)^{\frac{3}{4}}}{\sqrt{\m\n}}t^{\frac{5\delta}{2}+\frac{1}{4}}\left(A(t)+A(1)\right),
\end{align*}
with
\begin{align*}
A(t):=\norh{f_m(t)}{N}{}\sqnor{f_n}(t) ,
\end{align*}
$\mathcal{Q}$ is the quadratic term:
\begin{align*}
\mathcal{Q}:=\frac{\max(\m,\n)^{3+\frac{3}{4}}}{\sqrt{\m\n}}\inte{1}{t}{s^{3\delta}s^{-\frac{1}{4}}\nor{f_m}{H^{N}}\sqnor{f_n}(s)}{s},
\end{align*}
and $\mathcal{C}$ is the cubic one:
\begin{align*}
\mathcal{C}:=\frac{\max(n,m)^{\frac{1}{2}}}{\sqrt{\m\n}}\inte{1}{t}{s^{3\delta}s^{\frac{3}{2}}\nor{f_m}{H^{N}}\sqnor{f_m}(s)\sqnor{f_n}(s)}{s}.
\end{align*}
\end{lem}

\begin{pr}
First of all, use the boundedness in the frequency space.
\begin{align*}
\norlp{I^{lf,r}_{m,n}}{2}{}&\lesssim t^{\frac{3\delta}{2}}\norlp{\inte{1}{t}{  \inte{}{}{ s\partial_{\xi}\phi \theta_{s^\delta}(|\eta|)\theta_{s^\delta}(|\xi-\eta|)\psi^{s} (\xi,\eta)e^{-is\phi}\frac{\hat{f}_{m}(\eta)}{\langle \eta\rangle_{m}}\frac{\hat{f}_{n}(\xi-\eta)}{\langle \xi-\eta\rangle_{n}}  }{\eta}  }{s}}{2}{}.
\end{align*}
\paragraph{\textbf{A. Integration by parts in $t$.}}\label{time-resonances-method}
Now we are going to use that there are no time resonances on the support of $\psi^{s}$, i.e. that $\phi$ does not vanish. This will allow us to write the following equality:
\begin{align*}
e^{-is\phi}=\frac{1}{-i\phi}\partial_{s}\left(e^{-is\phi}\right).
\end{align*}

Then write that 
\begin{align*} &t^{\frac{3\delta}{2}}\inte{1}{t}{  \inte{}{}{ s\partial_{\xi}\phi \theta_{s^\delta}(|\eta|)\theta_{s^\delta}(|\xi-\eta|)\psi^{s} (\xi,\eta)e^{-is\phi}\frac{\hat{f}_{m}(\eta)}{\langle \eta\rangle_{m}}\frac{\hat{f}_{n}(\xi-\eta)}{\langle \xi-\eta\rangle_{n}}  }{\eta}  }{s} \\
&=  t^{\frac{3\delta}{2}}\norlp{\inte{1}{t}{  \inte{}{}{ s\partial_{\xi}\phi \theta_{s^\delta}(|\eta|)\theta_{s^\delta}(|\xi-\eta|)\psi^{s} (\xi,\eta)\frac{1}{-i\phi}\partial_{s}\left(e^{-is\phi}\right)\frac{\hat{f}_{m}(\eta)}{\langle \eta\rangle_{m}}\frac{\hat{f}_{n}(\xi-\eta)}{\langle \xi-\eta\rangle_{n}}  }{\eta}  }{s}}{2}{},
\end{align*}
and perform an integration by parts. This operation leads to

\begin{align*}
t^{\frac{3}{2}\delta}\inte{1}{t}{\inte{}{}{ s\partial_{\xi}\phi \theta_{s^\delta}(|\eta|)\theta_{s^\delta}(|\xi-\eta|)\psi^{s} (\xi,\eta)\frac{1}{-i\phi}\partial_{s}\left(e^{-is\phi}\right)\frac{\hat{f}_{m}(\eta)}{\langle \eta\rangle_{m}}\frac{\hat{f}_{n}(\xi-\eta)}{\langle \xi-\eta\rangle_{n}}  }{\eta}  }{s}=\sum_{i=0}^{4}I_{l,r}^{i},
\end{align*}
where 
\begin{align*}
I_{l,r}^{0}&:=t^{\frac{3}{2}\delta}\left[\inte{}{}{ s\partial_{\xi}\phi \theta_{s^\delta}(|\eta|)\theta_{s^\delta}(|\xi-\eta|)\psi^{s} (\xi,\eta)\frac{1}{-i\phi}e^{-is\phi}\frac{\hat{f}_{m}(\eta)}{\langle \eta\rangle_{m}}\frac{\hat{f}_{n}(\xi-\eta)}{\langle \xi-\eta\rangle_{n}}  }{\eta}\right]_{1}^{t},\\
I_{l,r}^1&:=t^{\frac{3}{2}\delta}\inte{1}{t}{\inte{}{}{ \partial_{\xi}\phi \theta_{s^\delta}(|\eta|)\theta_{s^\delta}(|\xi-\eta|)\psi^{s} (\xi,\eta)\frac{1}{-i\phi}e^{-is\phi}\frac{\hat{f}_{m}(\eta)}{\langle \eta\rangle_{m}}\frac{\hat{f}_{n}(\xi-\eta)}{\langle \xi-\eta\rangle_{n}}  }{\eta}  }{s},\\
I_{l,r}^{2}&:=t^{\frac{3}{2}\delta}\inte{1}{t}{\inte{}{}{ s\partial_{\xi}\phi \theta_{s^\delta}(|\eta|)\theta_{s^\delta}(|\xi-\eta|)\psi^{s} (\xi,\eta)\frac{1}{-i\phi}e^{-is\phi}\frac{\partial_{s}\hat{f}_{m}(\eta)}{\langle \eta\rangle_{m}}\frac{\hat{f}_{n}(\xi-\eta)}{\langle \xi-\eta\rangle_{n}}  }{\eta}  }{s},\\
I_{l,r}^{3}&:=t^{\frac{3}{2}\delta}\inte{1}{t}{\inte{}{}{ s\partial_{\xi}\phi \theta_{s^\delta}(|\eta|)\theta_{s^\delta}(|\xi-\eta|)\psi^{s} (\xi,\eta)\frac{1}{-i\phi}e^{-is\phi}\frac{\hat{f}_{m}(\eta)}{\langle \eta\rangle_{m}}\frac{\partial_{s}\hat{f}_{n}(\xi-\eta)}{\langle \xi-\eta\rangle_{n}}  }{\eta}  }{s},\\
I_{l,r}^{4}&:=t^{\frac{3}{2}\delta}\inte{1}{t}{\inte{}{}{ s\partial_{\xi}\phi \partial_{s}\left[\theta_{s^\delta}(|\eta|)\theta_{s^\delta}(|\xi-\eta|)\psi^{s} (\xi,\eta)\right]\frac{1}{-i\phi}e^{-is\phi}\frac{\hat{f}_{m}(\eta)}{\langle \eta\rangle_{m}}\frac{\hat{f}_{n}(\xi-\eta)}{\langle \xi-\eta\rangle_{n}}  }{\eta}  }{s}.
\end{align*}

\paragraph{\textbf{B. A preliminary result.}}

The terms $I^{j}_{l,r}$ can be written as bilinear multipliers associated with the same symbol. We are going to take advantage of it and prove a general result about the multiplier associated to the following symbol:
\begin{align*}
S(\xi,\eta):=\frac{1}{(\sqrt{n+1}+\sqrt{m+1})^{2}t^\delta}\partial_{\xi}\phi(\xi,\eta) \theta_{s^\delta}(|\eta|)\theta_{s^\delta}(|\xi-\eta|)\psi^{s} (\xi,\eta)\frac{1}{\phi(\xi,\eta)}.
\end{align*}
\begin{lem}\label{bilintr}
For all $g$ and $h$ we have the following inequality.
\begin{align*}
\norlp{T_{S}\left(e^{-is\cro{D}_{m}}\frac{g}{\cro{D}_{m}},e^{-is\cro{D}_{n}}\frac{h}{\cro{D}_{n}}\right)}{2}{}&\lesssim  \frac{\max(n,m)^{\frac{3}{4}}}{\sqrt{\m\n}}s^{-\frac{1}{4}}\norlp{g}{2}{}\sqrt{\nor{h}{H^{N}}\nor{h}{B_{s}}},
\end{align*}
with $T_{S}$ the bilinear operator associated to $S$ as defined in (\ref{def-bilin-mult}).
\end{lem}
\begin{pr}
We want to apply Theorem \ref{estbilinband}: to do this, we have to estimate the size of the support of $S$ and its behavior with derivation operators:
\begin{lem}\label{Ill-lt-notC-S}
The symbol $S$ satisfies the hypotheses of Theorem \ref{estbilinband}, with 
 \begin{align*}
 \rho=s^{\delta},~\omega=\frac{1}{(\sqrt{n+1}+\sqrt{m+1})^{2}s^\delta},~\mu=\left(\sqrt{\max(\m,\n)}s^{3\delta}(\sqrt{n+1}+\sqrt{m+1})^{2}s^\delta\right)^{-1}.
 \end{align*}
 \end{lem}
\begin{prf}{Lemma}{\ref{Ill-lt-notC-S}}
We are only going to detemrine a value for $\mu$, i.e. we are going to prove that 
\begin{align*}
\left|\partial_{\xi}^{a}\partial_{\eta}^{b}S(\xi,\eta)\right|\lesssim \left(\sqrt{\max(\m,\n)}s^{3\delta}(\sqrt{n+1}+\sqrt{m+1})^{2}s^\delta\right)^{a+b}.
\end{align*}
So as to prove this inequality, we need to understand the effect of differentiation on each factor in $S$.
\begin{itemize}
\item $\left|\partial_{\xi}\phi\right|\leq 2$ and $\left|\partial^{a}_{\xi}\partial^{b}_{\eta}\phi\right|\leq 2$ for all $a,b$ such that $a+b\geq 1$.
\item $\left|\partial^{a}_{\xi}\partial^{b}_{\eta}\theta\right|\lesssim \frac{1}{s^{(a+b)\delta}}\lesssim 1$ for $s>1$.
\item $\phi_{nr}=\theta\left(\sqrt{\max(\m,\n)}s^{3\delta}(\sqrt{n+1}+\sqrt{m+1})^2 s^{\delta}|\partial_{\eta}\phi|\right)$ hence, given the boundedness of derivatives of $\partial_{\eta}\phi$,
\begin{align*}
\left|\partial_{\xi}^{a}\partial_{\eta}^{b}\psi^{s}\right|\lesssim \left(\sqrt{\max(\m,\n)}s^{3\delta}(\sqrt{n+1}+\sqrt{m+1})^{2}s^\delta\right)^{a+b}.
\end{align*}
\item Finally, by Proposition \ref{timres-width}, 
\begin{align*}
\frac{1}{|\phi|}\lesssim (\sqrt{n+1}+\sqrt{m+1})^{2}s^\delta
\end{align*} and 
\begin{align*}
\left|\partial_{\eta}\left(\frac{1}{\phi}\right)\right|&=\left|\frac{\partial_{\eta}\phi}{\phi^2}\right|\\
&\lesssim \frac{1}{|\phi|^2}\\
&\lesssim \left((\sqrt{n+1}+\sqrt{m+1})^{2}s^\delta\right)^2.
\end{align*}
\end{itemize}
This ends the proof of Lemma \ref{Ill-lt-notC-S}.
\end{prf}\\

Then applying Theorem \ref{estbilinband} leads to
\begin{align*}
\norlp{T_{S}\left(e^{-is\cro{D}_{m}}\frac{g}{\cro{D}_{m}},e^{-is\cro{D}_{n}}\frac{h}{\cro{D}_{n}}\right)}{2}{}\lesssim \sqrt{\max(\m,\n)}s^{3\delta}\norlp{\frac{g}{\cro{D}_{m}}}{2}{}\frac{\n^{\frac{1}{4}}}{\sqrt{s}}\norw{\frac{h}{\cro{D}_{n}}}{\frac{3}{2}}{1}{}.
\end{align*}
Thanks to the linear multiplier estimate (\ref{propmult}), we finally obtain
\begin{align*}
\norlp{T_{S}\left(e^{-is\cro{D}_{m}}\frac{g}{\cro{D}_{m}},e^{-is\cro{D}_{n}}\frac{h}{\cro{D}_{n}}\right)}{2}{}&\lesssim\sqrt{\max(\m,\n)}s^{3\delta} \frac{\max(n,m)^{\frac{1}{4}}}{\sqrt{\m\n}}\frac{1}{\sqrt{s}}\norlp{g}{2}{}\norw{h}{\frac{3}{2}}{1}{}\\
&\lesssim  \sqrt{\max(\m,\n)}s^{3\delta}\frac{\max(n,m)^{\frac{1}{4}}}{\sqrt{\m\n}}s^{-\frac{1}{4}}\norlp{g(s)}{2}{}\sqnor{h}(s).
\end{align*}
The result is the same when exchanging the roles of $g$ and $h$. This ends the proof of Lemma \ref{bilintr}.
\end{pr}

\paragraph{\textbf{C. Application to $I^{0}_{l,r}$.}} The integral $I^{0}_{l,r}$ can be rewritten as follows.
\begin{align*}
I_{l,r}^{0}=&t^{\frac{3}{2}\delta}\inte{}{}{ t\partial_{\xi}\phi \theta_{s^\delta}(|\eta|)\theta_{s^\delta}(|\xi-\eta|)\psi^{t} (\xi,\eta)\frac{1}{-i\phi}e^{-it\phi}\frac{\hat{f}_{m}(t,\eta)}{\langle \eta\rangle_{m}}\frac{\hat{f}_{n}(t,\xi-\eta)}{\langle \xi-\eta\rangle_{n}}  }{\eta}\\
&- \inte{}{}{ \partial_{\xi}\phi \theta_{s^\delta}(|\eta|)\theta_{s^\delta}(|\xi-\eta|)\psi^{1} (\xi,\eta)\frac{1}{-i\phi}e^{-i\phi}\frac{\hat{f}_{m}(1,\eta)}{\langle \eta\rangle_{m}}\frac{\hat{f}_{n}(1,\xi-\eta)}{\langle \xi-\eta\rangle_{n}}  }{\eta}\\
=&t^{\frac{5}{2}\delta+1} (\sqrt{n+1}+\sqrt{m+1})^{2} T_{S}\left(e^{-it\cro{D}_{m}}\frac{f_m(t)}{\cro{D}_{m}},e^{-it\cro{D}_{n}}\frac{f_n(t)}{\cro{D}_{n}}\right)\\
&- (\sqrt{n+1}+\sqrt{m+1})^{2}T_{S}\left(e^{-i\cro{D}_{m}}\frac{f_m(1)}{\cro{D}_{m}},e^{-i\cro{D}_{n}}\frac{f_n(1)}{\cro{D}_{n}}\right).
\end{align*}
Then thanks to Lemma \ref{bilintr} we have the following inequality:
\begin{align}\label{notcondC-ineg1}
\frac{1}{\sqrt{t}}\norlp{I_{l,r}^{0}}{2}{}&\lesssim (\sqrt{n+1}+\sqrt{m+1})^{2}t^{\frac{5\delta}{2}+\frac{1}{4}} \frac{\max(n,m)^{\frac{3}{4}}}{\sqrt{\m\n}}\left(A(t)+A(1)\right),
\end{align}
with $A(t)=\norh{f_m(t)}{N}{}\sqrt{\nor{f_n(t)}{H^{N}}\nor{f_n(t)}{B_{t}}}$.
\paragraph{\textbf{D. Estimates for $I^{1}_{l,r}$.}}

First we give the following formula for $I^{1}_{l,r}$:
\begin{align*}
I^{1}_{l,r}=t^{\frac{5}{2}\delta}(\sqrt{n+1}+\sqrt{m+1})^{2}\inte{1}{t}{T_{S}\left(e^{-it\cro{D}_{m}}\frac{f_m}{\cro{D}_{m}},e^{-is\cro{D}_{n}}\frac{f_n}{\cro{D}_{n}}\right)}{s}.
\end{align*}
Then dividing by $\sqrt{t}$ and using Lemma \ref{bilintr} gives
\begin{align*}
\frac{1}{\sqrt{t}}\norlp{I^{1}_{l,r}}{2}{}&\lesssim \frac{1}{\sqrt{t}}t^{\frac{5}{2}\delta}(\sqrt{n+1}+\sqrt{m+1})^{2}\inte{1}{t}{\norlp{T_{S}\left(e^{-it\cro{D}_{m}}\frac{f_m}{\cro{D}_{m}},e^{-is\cro{D}_{n}}\frac{f_n}{\cro{D}_{n}}\right)}{2}{}}{s}\\
&\lesssim t^{\frac{5}{2}\delta}(\sqrt{n+1}+\sqrt{m+1})^{2}\inte{1}{t}{s^{-\frac{3\delta}{2}}\frac{\max(n,m)^{\frac{3}{4}}}{\sqrt{\m\n}}\norlp{f_m}{2}{}s^{-\frac{1}{4}}\mathcal{B}(f_n)(s)}{s}.
\end{align*}
Then we have the following inequality:
\begin{align}\label{notcondC-ineg2}
\frac{1}{\sqrt{t}}\norlp{I^{1}_{l,r}}{2}{}&\lesssim (\sqrt{n+1}+\sqrt{m+1})^{2}\frac{\max(n,m)^{\frac{3}{4}}}{\sqrt{\m\n}} \inte{1}{t}{s^\delta s^{-\frac{1}{4}}\nor{f_m}{H^{N}}\mathcal{B}(f_n)(s)}{s}.
\end{align}

\paragraph{\textbf{E. Estimates for $I^{2}_{l,r}$.}}
We proceed as for $I^{1}_{l,r}$ and get the following.
\begin{align*}
\frac{1}{\sqrt{t}}\norlp{I^{1}_{l,r}}{2}{}&\lesssim t^{\frac{5}{2}\delta}(\sqrt{n+1}+\sqrt{m+1})^{2}\inte{1}{t}{s\frac{\max(n,m)^{\frac{3}{4}}}{\sqrt{\m\n}}s^{-\frac{1}{4}}\norlp{\partial_{s}f_m}{2}{}\mathcal{B}(f_n)(s)}{s}.
\end{align*}
Now remark that 
\begin{align*}
\partial_{s}f_{+,m}&=e^{-is\cro{D}_{m}}u_m^{2}\\
&=e^{-is\cro{D}_{m}}\left(\frac{u_{+,m}-u_{-,m}}{\cro{D}_{m}}\right)^{2},
\end{align*}
where $u_{\pm,m}=e^{\mp is\cro{D}_m}f_{\pm,m}$. Now we use the space-time resonances method. For the sake of simplicity we write the following inequality:
\begin{align*}
\norlp{\partial_{s}f_m}{2}{}&=\norlp{e^{-is\cro{D}_{m}}u_m^{2}}{2}{}\\
&=\norlp{u_m^{2}}{2}{}\\
&\lesssim \norlp{u_m}{2}{}\norlp{u_m}{\infty}{}.
\end{align*}
Then given the expression of $u_m$ and the dispersion inequality $(\ref{dispequ})$, this inequality holds.
\begin{align*}
\norlp{\partial_{s}f_m}{2}{}&\lesssim \frac{\m^{\frac{1}{4}}}{\sqrt{s}}\norlp{f_m}{2}{}\norw{f_m}{\frac{3}{2}}{1}{}\\
&\lesssim \frac{\m^{\frac{1}{4}}}{s^{\frac{1}{4}}}\sqrt{\nor{f_m}{H^{N}}^3\nor{f_m}{B_{s}}}.
\end{align*}
This leads to
\begin{equation}\label{notcondC-ineg3}
\begin{aligned}
\frac{1}{\sqrt{t}}\norlp{I^{1}_{l,r}}{2}{}\lesssim & t^{\frac{5}{2}\delta}(\sqrt{n+1}+\sqrt{m+1})^{2}\\
&\times \inte{1}{t}{s\frac{\max(n,m)^{\frac{3}{4}}}{\sqrt{\m\n}}s^{-\frac{1}{4}}\frac{\m^{\frac{1}{4}}}{s^{\frac{1}{4}}}\sqrt{\nor{f_m}{H^{N}}^3\nor{f_m}{B_{s}}}\mathcal{B}(f_n)(s)}{s}.
\end{aligned}
\end{equation}

\paragraph{\textbf{F. Estimates for $I^{3}_{l,r}$ and $I^{4}_{l,r}$}} will be skipped since they can be treated as $I^{2}_{l,r}$, even if we differentiate the function $\theta_{s}$.\\
\\
We finally gather Inequalities $(\ref{notcondC-ineg1})$, $(\ref{notcondC-ineg2})$ and $(\ref{notcondC-ineg3})$: Lemma \ref{contrnotCp-res} is now proved.\end{pr}
\subsubsection{Outside the resonant set.}\label{notcondC-outside}
The lemma to prove is the following one:
\begin{lem}\label{contrnotCp-nonres}
For all $m$ and $n$ integers, $N\geq 3/2$, $t\geq 0$,
\begin{align*}
\frac{1}{\sqrt{t}}\norlp{I^{lf,nr}_{m,n}}{2}{}\lesssim\frac{\max(\m,\n)^{3+\frac{3}{4}}}{\sqrt{\m\n}}\inte{1}{t}{s^{3\delta}s^{\frac{1}{4}}\nor{f_m}{H^{N}}\mathcal{B}(f_n)(s)}{s},
\end{align*}
where $I^{lf,nr}_{m,n}$ is defined in (\ref{formule-Ilfnr}), page \pageref{formule-Ilfnr}.
\end{lem}

Here the estimates are almost a copy-paste of the method developped page \pageref{condC-outside}, with these two changes.
\begin{enumerate}
\item The term $|\partial_{\xi}\phi|$ is no longer smaller than $|\partial_{\eta}\phi|$. Then the quantity $\left|\frac{\partial_{\xi}\phi}{\partial_{\eta}\phi}\right|$ is bounded by $2^{j}$ in the zone $\partial_{\eta}\phi\sim-2^{-j}$.
\item We also define the symbols $S_{j}^{\pm}$, $S_{j,k}^{\pm}$, $M^1$, $M^2$, $M^3$, etc. However, given the change of localization around the space resonant set (cf. the definition of $\psi^{s}$ (\ref{def-psi-notcondC})), the equality (\ref{decomposition-phi}) becomes 
\begin{align}
(1-\psi^{s})\frac{\partial_{\xi}\phi}{\partial_{\eta}\phi}\frac{\partial_{\eta}^2\phi}{\partial_{\eta}\phi}= \sum_{1/2\leq 2^{j}\leq \sqrt{\max(\m,\n)}t^{4\delta}(\sqrt{n+1}+\sqrt{m+1})^2}S^{+}_{j}+S^{-}_{j},
\end{align}
\end{enumerate}
Since
\begin{align*}
\frac{1}{\sqrt{s}}\sum_{2^{j}\leq \sqrt{\max(\m,\n)}t^{4\delta}(\sqrt{n+1}+\sqrt{m+1})^2}2^{2j}\lesssim \max(\m,\n)^{3}\frac{t^{8\delta}}{\sqrt{s}},
\end{align*}
we have the following estimate:

\begin{align*}
\frac{1}{\sqrt{t}}\norlp{I^{lf,nr}_{m,n}}{2}{}\lesssim \frac{\max(\m,\n)^{3+\frac{3}{4}}}{\sqrt{\m\n}}\inte{1}{t}{s^{\frac{3\delta}{2}}s^{-\frac{1}{4}}\nor{f_m}{H^{N}}\mathcal{B}(f_n)(s)}{s}.
\end{align*}

Lemmas \ref{contrnotCp-res} and \ref{contrnotCp-nonres} give Lemma \ref{contrIlflm-notcp}.\\
\\
Then, Proposition \ref{intermediate-result-I} is proven by combining Propositions \ref{smalltimes}, \ref{contrIlflm-cp} and \ref{contrIlflm-notcp}.
\begin{rmq}
The case where $p\leq m$ or $p\leq n$ should also be treated separately, but estimating this term is very similar to what we did in the case "$p>m$, $p>n$ but condition \ref{cond-res} not satisfied" (Section \ref{I-LFHE-Lt-notC}): actually it is even easier since we do not have time resonances.
\end{rmq}

Finally, Proposition \ref{contrIhf}, Proposition \ref{contrIhm} and Proposition \ref{contrIlfm} lead to Proposition \ref{prop-H-3/2}. \\
Propositions \ref{prop-H-3/2} and \ref{contrHN} give Theorem \ref{prop-contraction}, and consequently Theorem \ref{thm} is proven. \end{prf}

\part{Resonant system}\label{part-resonant-system}
\section{Derivation of the resonant system}\label{derivation-resonant-system}
If we want to study the dynamics of (\ref{equation}), a good way to proceed is to study the \emph{resonant system} associated to this equation. The Duhamel formula for (\ref{equation}) is the following

 \begin{equation*}
\widetilde{f}_{\pm,p}(t,\xi)=\widetilde{f}_{\pm,p}(0,\xi)+D_{1}(f,f),
\end{equation*}
where 
\begin{align*}
D_{1}(f,f)=\sum_{m,n}\sum_{\alpha,\beta=\pm 1}\alpha\beta\mathcal{M}(n,m,p)\inte{0}{t}{  \inte{\R}{}{  e^{\mp is\phi^{\alpha,\beta}_{m,n,p}}\frac{\widetilde{f}_{\alpha,m}(\eta)}{\cro{\eta}_{m}} \frac{\widetilde{f}_{\beta,n}(\xi-\eta)}{\cro{\xi-\eta}_{n}}  }{\eta}  }{s}.
\end{align*}
It describes the interaction between Fourier modes (input frequency $\eta$ and $\xi-\eta$, output frequency $\xi$) or between Hermite modes (input modes $m$ and $n$, output mode $p$). We saw in the study of existence and uniqueness for the original equation that some modes are \emph{resonant}: more particularly, they can be \emph{space resonant} when $\partial_{\eta}\phi^{\alpha,\beta}_{m,n,p}(\xi,\eta)=0$ or \emph{time resonant} when $\phi^{\alpha,\beta}_{m,n,p}(\xi,\eta)=0$. These resonant interactions must be the ones governing the dynamics of the whole equation.\\
In this section we are going to determine these resonant interactions and the corresponding resonant equation.

\paragraph{\textbf{Space resonant interactions.}} The first step is to remove the interactions which are non resonant in space, i.e. the interactions such that $\partial_{\eta}\phi \neq 0$. In Appendix \ref{resonances-asympt}, the cancellation of $\partial_{\eta}\phi$ is studied in detail :
\begin{enumerate}
\item if $\alpha=-\beta$ and $m=n$, then $\partial_{\eta}\phi$ is identically zero for $\xi=0$, and does not vanishes for $\xi\neq 0$. Moreover $\phi$ never vanishes.
\item otherwise, for all $\xi$ there exists one and only one $\eta_{0}(\xi)$ such that $\partial_{\eta}\phi(\xi,\eta_{0}(\xi))=0$:
\begin{align*}
\eta_{0}(\xi):=\lambda^{\alpha,\beta}_{m,n} \xi,~\text{where}~\lambda^{\alpha,\beta}_{m,n}=\frac{1}{1+\alpha\beta\sqrt{\frac{n+1}{m+1}}}.
\end{align*}
\end{enumerate}

Hence it is natural to approximate the Duhamel formula as follows
\begin{enumerate}
\item first, if $m=n$ and $\alpha=-\beta$, the non-cancellation of $\partial_{\eta}\phi$ allows us to perform an integration by parts in $\eta$, and consequently gain a decay in $s$. That is why we are allowed to remove those Hermite modes
 \begin{equation*}
D_{1}(f,f)\sim \sum_{\substack{m,n\in\N\\\alpha,\beta\in\{\pm 1\} \\m\neq n\text{ or }\alpha\neq - \beta}}\alpha\beta\mathcal{M}(n,m,p)\inte{0}{t}{  \inte{\R}{}{  e^{\mp is\phi^{\alpha,\beta}_{m,n,p}}\frac{\widetilde{f}_{\alpha,m}(\eta)}{\cro{\eta}_{m}} \frac{\widetilde{f}_{\beta,n}(\xi-\eta)}{\cro{\xi-\eta}_{n}}  }{\eta}  }{s}.
\end{equation*}
\item then if $m\neq n\text{ or }\alpha\neq \beta$, we are in the framework of stationary phase Lemma: the behavior of the oscillating integral in $\eta$
\begin{align*}
\inte{\R}{}{  e^{\mp is\phi^{\alpha,\beta}_{m,n,p}}\frac{\widetilde{f}_{\alpha,m}(\eta)}{\cro{\eta}_{m}} \frac{\widetilde{f}_{\beta,n}(\xi-\eta)}{\cro{\xi-\eta}_{n}}  }{\eta}.
\end{align*}
is governed by the frequencies on which $\partial_{\eta}\phi$ vanish, i.e.
\begin{align*}
&\inte{\R}{}{  e^{\mp is\phi^{\alpha,\beta}_{m,n,p}}\frac{\widetilde{f}_{\alpha,m}(\eta)}{\cro{\eta}_{m}} \frac{\widetilde{f}_{\beta,n}(\xi-\eta)}{\cro{\xi-\eta}_{n}}  }{\eta} \sim \frac{C}{\sqrt{s|\partial^{2}_{\eta}\phi(\xi,\eta_{0})|}}e^{\mp is\phi^{\alpha,\beta}_{m,n,p}(\xi,\eta_{0}(\xi))}\frac{\widetilde{f}_{\alpha,m}(\eta_{0}(\xi))}{\cro{\eta_{0}(\xi)}_{m}} \frac{\widetilde{f}_{\beta,n}(\xi-\eta_{0}(\xi))}{\cro{\xi-\eta_{0}(\xi)}_{n}},
\end{align*}
for some constant $C$. Hence the following approximation for $D_{1}(f,f)$:
\begin{align*}
D_{1}(f,f)\sim \!\!\!\!\!\!\!\!\! \sum_{\substack{m,n\in\N\\\alpha,\beta\in\{\pm 1\} \\m\neq n\text{ or }\alpha\neq - \beta}}\!\!\!\!\!\!\! \alpha\beta\mathcal{M}(n,m,p)\inte{0}{t}{\frac{Ce^{\mp is\phi^{\alpha,\beta}_{m,n,p}(\xi,\eta_{0}(\xi))}}{\sqrt{s|\partial^{2}_{\eta}\phi(\xi,\eta_{0})|}}\frac{\widetilde{f}_{\alpha,m}(\eta_{0}(\xi))}{\cro{\eta_{0}(\xi)}_{m}} \frac{\widetilde{f}_{\beta,n}(\xi-\eta_{0}(\xi))}{\cro{\xi-\eta_{0}(\xi)}_{n}}}{s}.
\end{align*}
\end{enumerate}

\paragraph{\textbf{Time resonances.}} Now the space nonresonant interactions have been removed, the approximate formula we obtained is a sum of oscillating integrals of the form 
\begin{align*}
\inte{0}{t}{e^{\mp is\phi^{\alpha,\beta}_{m,n,p}(\xi,\eta_{0}(\xi))}F^{\alpha,\beta}_{m,n,p}(s,\xi)}{s}.
\end{align*}
The integrals for which $\phi^{\alpha,\beta}_{m,n,p}$ is different from $0$ are more likely to be neglectable compared to the ones where it is. But $\phi^{\alpha,\beta}_{m,n,p}$ vanishes on the zero set of $\partial_{\eta}\phi^{\alpha,\beta}_{m,n,p}$ if and only if the condition (\ref{cond-res}) is satisfied :
\begin{equation}\tag{C}
m^2+n^2+p^2-2mn-2pm-2pn-2m-2n-2p-3=0.
\end{equation}
Hence if (\ref{cond-res}) is not satisfied, the integral 

\begin{align*}
\inte{0}{t}{\frac{Ce^{\mp is\phi^{\alpha,\beta}_{m,n,p}(\xi,\eta_{0}(\xi))}}{\sqrt{s|\partial^{2}_{\eta}\phi(\xi,\eta_{0})|}}\frac{\widetilde{f}_{\alpha,m}(\eta_{0}(\xi))}{\cro{\eta_{0}(\xi)}_{m}} \frac{\widetilde{f}_{\beta,n}(\xi-\eta_{0}(\xi))}{\cro{\xi-\eta_{0}(\xi)}_{n}}}{s}
\end{align*}
is an oscillating integral. Hence the approximation
\begin{align*}
D_{1}(f,f)\sim \inte{0}{t}{
\sum_{\substack{m,n\in\Z\\ \alpha,\beta\in\{\pm 1\}\\ m\neq n\text{ or }\alpha\neq -\beta\\ (\ref{cond-res})\text{ satisfied}}} \frac{C}{\sqrt{s|\partial^{2}_{\eta}\phi(\xi,\eta_{0})|}}\frac{\widetilde{f}_{\alpha,m}(\lambda^{\alpha,\beta}_{m,n}\xi)}{\cro{\lambda^{\alpha,\beta}_{m,n}\xi}_{m}} \frac{\widetilde{f}_{\beta,n}((1-\lambda^{\alpha,\beta}_{m,n})\xi)}{\cro{(1-\lambda^{\alpha,\beta}_{m,n})\xi}_{n}} 
}{s}.
\end{align*}

The \emph{resonant equation} is the following one
\begin{align}\label{duh-resonant}
\widetilde{f}_{\pm,p}(t,\xi)=&\widetilde{f}_{\pm,p}(0,\xi)+ \inte{0}{t}{
\sum_{\substack{m,n\in\Z\\ \alpha,\beta\in\{\pm 1\}\\ m\neq n\text{ or }\alpha\neq -\beta\\ (\ref{cond-res})\text{ satisfied}}} \frac{C}{\sqrt{t|\partial^{2}_{\eta}\phi(\xi,\eta_{0})|}}\frac{\widetilde{f}_{\alpha,m}(\lambda^{\alpha,\beta}_{m,n}\xi)}{\cro{\lambda^{\alpha,\beta}_{m,n}\xi}_{m}} \frac{\widetilde{f}_{\beta,n}((1-\lambda^{\alpha,\beta}_{m,n})\xi)}{\cro{(1-\lambda^{\alpha,\beta}_{m,n})\xi}_{n}} 
}{s}.
\end{align}
Given the heuristic approach we made, it is natural that this equation approximates correctly the dynamics of the original equation. So as to prove this, we are going to proceed in two steps:
\begin{enumerate}
\item first of all, we have to show that the solutions of the resonant system exist for a time at least as long as the solutions of the original system (we are even going to get a longer existence time).
\item then we have to prove that the solutions of the resonant system are a good approximation (in a sense explained in section \ref{section-validity-approx}) of the solutions of the original one, as long as the former exists.
\end{enumerate}
\section{Long-time existence for the resonant system}\label{long-time-existence-resonant system}
The resonant system being simpler than the original equation, it seems reasonable to find an existence time at least as good as $1/\eps^{\frac{4}{3}}$ where $\eps$ is the size of the initial data.\\

Our target is the same as in the proof of Theorem \ref{thm} : we are going to prove contraction estimates for the operator $\tld{\mathrm{Res}}_{\pm,p}(f,f)$. In order not to be too redundant, we are going to gather all the technical details in one proposition.\\
Since the Duhamel formula is symmetric in $n$ and $m$, we are going to prove a contraction estimate for this ``half Duhamel formula'':
\begin{equation}\label{half-duh-res}
\tld{f}_{\pm,p}(t,\xi)=\tld{f}_{\pm,p}(0,\xi)+\inte{0}{t}{
\sum_{\substack{m,n\in\Z\\ \alpha,\beta\in\{\pm 1\}\\m\leq n\\ m\neq n\text{ or }\alpha\neq -\beta\\ (\ref{cond-res})\text{ satisfied}}} \mathcal{M}(m,n,p)\frac{C}{\sqrt{s|\partial^{2}_{\eta}\phi(\xi,\eta_{0})|}}\frac{\widetilde{f}_{\alpha,m}(\lambda^{\alpha,\beta}_{m,n}\xi)}{\cro{\lambda^{\alpha,\beta}_{m,n}\xi}_{m}} \frac{\widetilde{f}_{\beta,n}((1-\lambda^{\alpha,\beta}_{m,n})\xi)}{\cro{(1-\lambda^{\alpha,\beta}_{m,n})\xi}_{n}} 
}{s}.
\end{equation}
\begin{rmq}
This equation is different and much simpler to deal with than the Duhamel formula $(\ref{duh})$ for the original equation $(\ref{equation})$. In fact, it does not involve any integration in $\xi$: then estimating the integral term will be based on Young's inequality instead of Hölder-like inequalities.
\end{rmq}
We are going to give a few useful bounds for $\lambda=\lambda^{\alpha,\beta}_{m,n}:=\frac{1}{1+\alpha\beta\sqrt{\frac{n+1}{m+1}}}$.
\begin{enumerate}
\item If $\alpha\beta=1$, then $\lambda$ and $1-\lambda$ are bounded by $1$. Otherwise, since $m\neq n$, the maximum of $\lambda$ is reached for $|n-m|=1$. This leads to the bound
\begin{equation}\label{boundlambda}
|\lambda|\lesssim \min(m,n).
\end{equation}

\item Since we are in the case $m\leq n$, we also have the bound
\begin{equation}\label{bound1-lambda}
\frac{1}{1-\lambda}=1+\alpha\beta\sqrt{\frac{m+1}{n+1}}\leq 2.
\end{equation}
\end{enumerate}

\subsection{A preliminary bilinear estimate}

\begin{prop}\label{prop-fix-pt-res}
Let $a$, $b$, $c$, and $d$ be four integers, $m$ and $n$ two integers with $m<n$. Then if we define $S=S(a,b,c,d)=H^{a\wedge b}(\cro{x}^{A(c,d)})$, we have
\begin{align*}
\norlp {(\lambda\xi)^{a}((1-\lambda)\xi)^{b}\lambda^{e(c,d)}\frac{\partial_{\xi}^{c}\tld{f}_{\alpha,m}(\lambda\xi)}{\cro{\lambda\xi}_{m}}\frac{\partial_{\xi}^{d}\tld{f}_{\beta,n}((1-\lambda)\xi)}{\cro{(1-\lambda)\xi}_{n}}}{2}{\xi}\lesssim \frac{1}{\sqrt{mn}}\nor{f_m}{S}\nor{f_n}{S},
\end{align*}
where 
\begin{align*}
a\wedge b=\max(a,b),~
e(c,d)=\left\{ \begin{array}{rcl}
0&\text{ if }& c\leq d,\\
\frac{1}{2}&\text{ if }& c> d,
\end{array}\right.~and~ A(c,d)=\left\{\begin{array}{rcl}
\max(c,d)&\text{ if }&c\neq d,\\
c+1&\text{ if }&c=d.
\end{array}\right.
\end{align*}
\end{prop}
For example, for $a=1$, $b=N$, $c=d=0$, we have 
\begin{align*}
\norlp {(\lambda\xi)((1-\lambda)\xi)^{N}\frac{\tld{f}_{\alpha,m}(\lambda\xi)}{\cro{\lambda\xi}_{m}}\frac{\tld{f}_{\beta,n}((1-\lambda)\xi)}{\cro{(1-\lambda)\xi}_{n}}}{2}{\xi}\lesssim \frac{1}{\sqrt{mn}}\nor{f_m}{H^N(\cro{x})}\nor{f_n}{H^N(\cro{x})}.
\end{align*}
If $a=3/2$, $b=N$, $c=1$, $d=0$, we have 
\begin{align*}
\norlp {(\lambda\xi)^{3/2}((1-\lambda)\xi)^{N}\lambda^{1/2}\frac{\partial_{\xi}\tld{f}_{\alpha,m}(\lambda\xi)}{\cro{\lambda\xi}_{m}}\frac{\tld{f}_{\beta,n}((1-\lambda)\xi)}{\cro{(1-\lambda)\xi}_{n}}}{2}{\xi}\lesssim \frac{1}{\sqrt{mn}}\nor{f_m}{H^{N}(\cro{x})}\nor{f_n}{H^{N}(\cro{x})}.
\end{align*}
\begin{pr}
Let us only prove two cases: $a>b$, $c>d$ and $a>b$, $c<d$.\\
If $a>b$ and $c>d$, then
\begin{align*}
&\norlp {(\lambda\xi)^{a}((1-\lambda)\xi)^{b}\sqrt{\lambda}\frac{\partial_{\xi}^{c}\tld{f}_{\alpha,m}(\lambda\xi)}{\cro{\lambda\xi}_{m}}\frac{\partial_{\xi}^{d}\tld{f}_{\beta,n}((1-\lambda)\xi)}{\cro{(1-\lambda)\xi}_{n}}}{2}{\xi}\\
&\leq 
\norlp {(\lambda\xi)^{a}\sqrt{\lambda}\frac{\partial_{\xi}^{c}\tld{f}_{\alpha,m}(\lambda\xi)}{\cro{\lambda\xi}_{m}}}{2}{\xi}\norlp{((1-\lambda)\xi)^{b}\frac{\partial_{\xi}^{d}\tld{f}_{\beta,n}((1-\lambda)\xi)}{\cro{(1-\lambda)\xi}_{n}}}{\infty}{\xi}\\
&\leq 
\norlp {\xi^{a}\frac{\partial_{\xi}^{c}\tld{f}_{\alpha,m}(\xi)}{\cro{\xi}_{m}}}{2}{\xi}\norlp{\xi^{b}\frac{\partial_{\xi}^{d}\tld{f}_{\beta,n}(\xi)}{\cro{\xi}_{n}}}{\infty}{\xi},
\end{align*}
by a change of variable (dilation). Then, since the norm of the multiplier $\frac{1}{\cro{D}_{p}}$ is bounded by $\frac{1}{\sqrt{p}}$ (Proposition \ref{propmult}),
\begin{align*}
&\norlp {(\lambda\xi)^{a}((1-\lambda)\xi)^{b}\sqrt{\lambda}\frac{\partial_{\xi}^{c}\tld{f}_{\alpha,m}(\lambda\xi)}{\cro{\lambda\xi}_{m}}\frac{\partial_{\xi}^{d}\tld{f}_{\beta,n}((1-\lambda)\xi)}{\cro{(1-\lambda)\xi}_{n}}}{2}{\xi}\leq 
\frac{1}{\sqrt{mn}}\norlp {\xi^{a}\partial_{\xi}^{c}\tld{f}_{\alpha,m}(\xi)}{2}{\xi}\norlp{\xi^{b}\partial_{\xi}^{d}\tld{f}_{\beta,n}(\xi)}{\infty}{\xi}.
\end{align*}
Then by $L^2$ continuity of the Fourier transform,
\begin{align*}
\norlp {\xi^{a}\partial_{\xi}^{c}\tld{f}_{\alpha,m}(\xi)}{2}{\xi}&\lesssim 
\norlp {D^{a}\left(x^c f_{\alpha,m}\right)}{2}{x}\\
&\lesssim \nor{f_{\alpha,m}}{H^{a}(\cro{x}^{c})}.
\end{align*}
By $L^{1}\rightarrow L^{\infty}$ continuity of the Fourier transform, 
\begin{align*}
\norlp{\xi^{b}\partial_{\xi}^{d}\tld{f}_{\beta,n}(\xi)}{\infty}{\xi}&\lesssim \norlp{D^b x^d f_{\beta,n}}{1}{x}\\
&\lesssim \sqrt{\norlp{D^b x^d f_{\beta,n}}{2}{x}\norlp{x D^b x^d f_{\beta,n}}{2}{x}},
\end{align*}
by Proposition \ref{leml1}. Then, 
\begin{align*}
\norlp{D^b x^d f_{\beta,n}}{2}{x}\lesssim \nor{f_{\beta,n}}{H^b (\cro{x}^d)}\lesssim \nor{f_{\beta,n}}{H^a (\cro{x}^c)}.
\end{align*}
Moreover,
\begin{align*}
xD^b = D^{b}x-bD^{b-1}.
\end{align*}
Hence,
\begin{align*}
\norlp{x D^b x^d f_{\beta,n}}{2}{x}&\lesssim \nor{f_{\beta,n}}{H^b(\cro{x}^{d+1})}\\
&\lesssim\nor{f_{\beta,n}}{H^a(\cro{x}^{c})}.
\end{align*}
This proves the theorem in the case $a>b$ and $c>d$.\\
In the case $a>b$, $c<d$, then
\begin{align*}
\norlp {(\lambda\xi)^{a}((1-\lambda)\xi)^{b}\frac{\partial_{\xi}^{c}\tld{f}_{\alpha,m}(\lambda\xi)}{\cro{\lambda\xi}_{m}}\frac{\partial_{\xi}^{d}\tld{f}_{\beta,n}((1-\lambda)\xi)}{\cro{(1-\lambda)\xi}_{n}}}{2}{\xi}
&\leq 
\norlp {(\lambda\xi)^{a}\frac{\partial_{\xi}^{c}\tld{f}_{\alpha,m}(\lambda\xi)}{\cro{\lambda\xi}_{m}}}{\infty}{\xi}\norlp{((1-\lambda)\xi)^{b}\frac{\partial_{\xi}^{d}\tld{f}_{\beta,n}((1-\lambda)\xi)}{\cro{(1-\lambda)\xi}_{n}}}{2}{\xi}\\
&\leq 
\norlp {\xi^{a}\frac{\partial_{\xi}^{c}\tld{f}_{\alpha,m}(\xi)}{\cro{\xi}_{m}}}{\infty}{\xi}\frac{1}{\sqrt{1-\lambda}}\norlp{\xi^{b}\frac{\partial_{\xi}^{d}\tld{f}_{\beta,n}(\xi)}{\cro{\xi}_{n}}}{2}{\xi},
\end{align*}
by a change of variables. Then using the bound $(\ref{bound1-lambda})$ reduces the problem to the previous case.

Hence gathering (\ref{contr-est-res-syst-noweight}), (\ref{contr-est-res-syst-weight-x}) and (\ref{contr-est-res-syst-weight-x2}) prove Proposition \ref{bilinear-estimate-resonant-system}.
\end{pr}

\subsection{Bilinear estimates for (\ref{half-duh-res}).}

We are going to use Proposition \ref{prop-fix-pt-res} to prove the following one:
\begin{prop}\label{bilinear-estimate-resonant-system}
Let $N\geq \frac{3}{2}$. Then, if 
\begin{align*}
I(\xi):=\frac{1}{\sqrt{|\partial^{2}_{\eta}\phi(\xi,\eta_{0})|}}\frac{\widetilde{f}_{\alpha,m}(\lambda\xi)}{\cro{\lambda\xi}_{m}} \frac{\widetilde{f}_{\beta,n}((1-\lambda)\xi)}{\cro{(1-\lambda)\xi}_{n}},
\end{align*}
we have the following bounds
\begin{align*}
\nor{\mathcal{F}^{-1}(I)}{H^{N}(\cro{x})}\lesssim m^{\frac{3}{2}}\frac{1}{\sqrt{mn}}\nor{f_{\alpha,m}}{H^{N}(\cro{x})}\nor{f_{\beta,n}}{H^{N}(\cro{x})},\\
\nor{\mathcal{F}^{-1}(I)}{H^{N}(\cro{x}^2)}\lesssim m^{\frac{9}{2}}\frac{1}{\sqrt{mn}}\nor{f_{\alpha,m}}{H^{N}(\cro{x}^2)}\nor{f_{\beta,n}}{H^{N}(\cro{x}^2)}.
\end{align*}
\end{prop}
\begin{pr}
We are going to prove the estimate for the weight $\cro{x}^2$, the weight $\cro{x}$ being dealt with similarly. We are simply going to estimate the $L^2$ norms of $|\xi|^N I$, $|\xi|^N \partial_{\xi}I$ and $|\xi|^N\partial_{\xi}^{2}I$.
\begin{enumerate}
\item 
We first recall that

\begin{align*}
\partial_{\eta}^{2}\phi^{\alpha\beta}_{m,n,p}(\xi,\eta)=\frac{2m+2}{(\eta^2+2m+2)^{\frac{3}{2}}}+\eps\frac{2n+2}{\left((\xi-\eta)^2+2n+2\right)^{\frac{3}{2}}}\cdot
\end{align*}
where $\eps=\alpha\beta$. Then, if $\eta=\eta_0:=\frac{\xi}{1+\eps\sqrt{\frac{n+1}{m+1}}}$, a  calculation shows that
\begin{align*}
\partial_{\eta}^{2}\phi^{\alpha\beta}_{m,n,p}(\xi,\eta_{0}(\xi))= \frac{2m+2}{\lambda\left(\lambda^2\xi^2+2m+2\right)^{\frac{3}{2}}}\cdot
\end{align*}
Hence 
\begin{align*}
\frac{1}{\sqrt{|\partial_{\eta}^{2}\phi|}}&\lesssim \frac{\lambda\cro{\lambda\xi}_{m}^{\frac{3}{2}}}{2m+2} \lesssim \cro{\lambda\xi}_{m}^{\frac{3}{2}},
\end{align*}
by the bound $(\ref{boundlambda})$ on $\lambda$. Then
\begin{align*}
|\xi|^N|I|\lesssim |\xi|^N \cro{\lambda\xi}_{m}^{\frac{3}{2}}\frac{\widetilde{f}_{\alpha,m}(\lambda\xi)}{\cro{\lambda\xi}_{m}} \left|\frac{\widetilde{f}_{\beta,n}((1-\lambda)\xi)}{\cro{(1-\lambda)\xi}_{n}}\right|\cdot
\end{align*}
If we write $|\xi|^N=\frac{1}{(1-\lambda)^N}|(1-\lambda)\xi|^N$, by the bound $(\ref{bound1-lambda})$ we have
\begin{align*}
|\xi|^N|I|\lesssim m^{\frac{3}{4}}(\lambda\xi)^{a}((1-\lambda)\xi)^{b}\left|\frac{\partial_{\xi}^{c}\tld{d}_{\alpha,m}(\lambda\xi)}{\cro{\lambda\xi}_{m}}\frac{\partial_{\xi}^{f}\tld{f}_{\beta,n}((1-\lambda)\xi)}{\cro{(1-\lambda)\xi}_{n}}\right|,
\end{align*}
with
\begin{align*}
a=\frac{3}{2},~b=N,~c=0,~d=0.
\end{align*}
Then we are in the framework of Proposition \ref{prop-fix-pt-res}, and we conclude by 
\begin{align}\label{contr-est-res-syst-noweight}
\nor{|\xi|^N I}{L^2}\lesssim \frac{m^{\frac{3}{4}}}{\sqrt{mn}}\nor{f_{\alpha,m}}{H^{N}}\nor{f_{\beta,n}}{H^{N}}.
\end{align}
\item Then we have to do the same for the weighted norms, i.e. for the $\xi$ derivative of $I$. We can write \begin{align*}
\partial_{\xi}\left(\frac{1}{\sqrt{|\partial^{2}_{\eta}\phi(\xi,\eta_{0})|}}\frac{\widetilde{f}_{\alpha,m}(\lambda\xi)}{\cro{\lambda\xi}_{m}} \frac{\widetilde{f}_{\beta,n}((1-\lambda)\xi)}{\cro{(1-\lambda)\xi}_{n}}\right)= I_{1}+I_{2}+I_{3}+I_{4}+I_{5},
\end{align*}
where
\begin{align*}
I_{1}&:=\partial_{\xi}\left({\frac{1}{\sqrt{|\partial^{2}_{\eta}\phi(\xi,\eta_{0})|}}}\right)\frac{\widetilde{f}_{\alpha,m}(\lambda\xi)}{\cro{\lambda\xi}_{m}} \frac{\widetilde{f}_{\beta,n}((1-\lambda)\xi)}{\cro{(1-\lambda)\xi}_{n}}
,\\
I_{2}&:=\frac{1}{\sqrt{|\partial^{2}_{\eta}\phi(\xi,\eta_{0})|}}\lambda\frac{\partial_{\xi}\widetilde{f}_{\alpha,m}(\lambda\xi)}{\cro{\lambda\xi}_{m}} \frac{\widetilde{f}_{\beta,n}((1-\lambda)\xi)}{\cro{(1-\lambda)\xi}_{n}},\\
I_{3}&:=\frac{1}{\sqrt{|\partial^{2}_{\eta}\phi(\xi,\eta_{0})|}}\lambda^2\xi\frac{\widetilde{f}_{\alpha,m}(\lambda\xi)}{\cro{\lambda\xi}_{m}^{3}} \frac{\widetilde{f}_{\beta,n}((1-\lambda)\xi)}{\cro{(1-\lambda)\xi}_{n}},\\
I_{4}&:=\frac{1}{\sqrt{|\partial^{2}_{\eta}\phi(\xi,\eta_{0})|}}\frac{\widetilde{f}_{\alpha,m}(\lambda\xi)}{\cro{\lambda\xi}_{m}} (1-\lambda)\frac{\partial_{\xi}\widetilde{f}_{\beta,n}((1-\lambda)\xi)}{\cro{(1-\lambda)\xi}_{n}},\\
I_{5}&:=\frac{1}{\sqrt{|\partial^{2}_{\eta}\phi(\xi,\eta_{0})|}}\frac{\widetilde{f}_{\alpha,m}(\lambda\xi)}{\cro{\lambda\xi}_{m}} (1-\lambda)^2\xi\frac{\widetilde{f}_{\beta,n}((1-\lambda)\xi)}{\cro{(1-\lambda)\xi}_{n}^{3}}.
\end{align*}
The expression of $\partial_{\xi}\left({\frac{1}{\sqrt{|\partial^{2}_{\eta}\phi(\xi,\eta_{0})|}}}\right)$ is 
\begin{align*}
\partial_{\xi}\left({\frac{1}{\sqrt{|\partial^{2}_{\eta}\phi(\xi,\eta_{0})|}}}\right)=\frac{\lambda^{\frac{3}{2}}(\lambda\xi)}{\sqrt{2m+2}((\lambda\xi)^2+2m+2)^{\frac{1}{4}}}.
\end{align*}
This given, we can write that for every $j$,
\begin{align*}
|\xi|^{N}I_{j}\lesssim m^{\gamma_{j}}(\lambda\xi)^{a_{j}}((1-\lambda)\xi)^{b_{j}}\lambda^{e(c_j,d_j)}\frac{\partial_{\xi}^{c_{j}}\tld{f}_{\alpha,m}(\lambda\xi)}{\cro{\lambda\xi}_{m}}\frac{\partial_{\xi}^{d_j}\tld{f}_{\beta,n}((1-\lambda)\xi)}{\cro{(1-\lambda)\xi}_{n}}.
\end{align*}
The values of the coefficients are summed up in this array.
\begin{align*}
\begin{array}{|c|c|c|c|c|c|}
\hline
j&a_j & b_j & c_j & d_j & \gamma_j \\
\hline
1 &1 &N & 0& 0& 3/4\\
\hline
2 & 3/2&N & 1& 0&3/4\\
\hline
3 &1/2 & N& 0& 0&3/2\\
\hline
4 & 3/2 & N & 0& 1&3/2\\
\hline
5 & 3/2& N & 0& 0& 3/2\\
\hline
\end{array}
\end{align*}
This implies, by Proposition \ref{prop-fix-pt-res},
\begin{align}\label{contr-est-res-syst-weight-x}
\nor{|\xi|^N \partial_{\xi}I}{L^2}\lesssim \frac{m^{\frac{3}{2}}}{\sqrt{mn}}\nor{f_{\alpha,m}}{H^{N}(\cro{x_{1}})}\nor{f_{\beta,n}}{H^{N}(\cro{x_{1}})}.
\end{align}
\item Dealing with the weight $x^2$ is very similar: we have to compute
\begin{align*}
\partial_{\xi}\left(\frac{\lambda^{\frac{3}{2}}(\lambda\xi)}{\sqrt{2m+2}((\lambda\xi)^2+m+1)^{\frac{1}{4}}}\right)=\frac{\lambda^{\frac{9}{2}}\frac{\xi^2}{2}+2\lambda^{\frac{5}{2}}(m+1)}{\sqrt{2m+2}((\lambda\xi)^{2}+2m+2)^{\frac{5}{4}}}.
\end{align*}
This quantity can be bounded as follows:
\begin{align*}
\left|\partial_{\xi}\left(\frac{\lambda^{\frac{3}{2}}(\lambda\xi)}{\sqrt{2m+2}((\lambda\xi)^2+m+1)^{\frac{1}{4}}}\right)\right|\lesssim \lambda^{\frac{15}{4}}+\lambda^{\frac{7}{4}}.
\end{align*}
When applying $\partial_{\xi}$ to $I_{j}$ we will get a sum of five terms $I_{j,1},\dots,I_{j,5}$. We state that for all $1\leq j\leq 5$ and $1\leq k\leq 5$, for all $N\in\N$,
\begin{align*}
|\xi|^{N}I_{j,k}\lesssim m^{\gamma_{j,k}}(\lambda\xi)^{a_{j,k}}((1-\lambda)\xi)^{b_{j,k}}\lambda^{e(c_{j,k},d_{j,k})}\frac{\partial_{\xi}^{c_{j,k}}\tld{f}_{\alpha,m}(\lambda\xi)}{\cro{\lambda\xi}_{m}}\frac{\partial_{\xi}^{d_{j,k}}\tld{f}_{\beta,n}((1-\lambda)\xi)}{\cro{(1-\lambda)\xi}_{n}},
\end{align*}
with the parameters satisfying the following bounds:
\begin{align*}
a_{j,k}&\leq \frac{3}{2},\\
b_{j,k}&=N,\\
c_{j,k}+d_{j,k}&\leq 2,\\
\gamma_{j,k}&\leq \frac{9}{2}.
\end{align*}
Hence we obtain
\begin{align}\label{contr-est-res-syst-weight-x2}
\nor{|\xi|^N \partial_{\xi}I}{L^2}\lesssim \frac{m^{\frac{9}{2}}}{\sqrt{mn}}\nor{f_{\alpha,m}}{H^{N}(\cro{x_{1}})}\nor{f_{\beta,n}}{H^{N}(\cro{x_{1}})}.
\end{align}
\end{enumerate}
Hence gathering (\ref{contr-est-res-syst-noweight}), (\ref{contr-est-res-syst-weight-x}) and (\ref{contr-est-res-syst-weight-x2}) we prove Proposition \ref{bilinear-estimate-resonant-system}.
\end{pr}

\subsection{From the estimate to the theorem}

Now we are going to prove Theorem \ref{thm-existence-resonant-system} (in the case of a weight equal to $\cro{x_{1}}^2$). More precisely we are going to prove the following proposition which leads to the theorem by a contraction argument.

\begin{prop}\label{fixed-point-resonant}
Let $N\geq\frac{3}{2}$, $M>6$. Then if $f_0\in\He^M H^{N}(\cro{x_{1}}^2)$, then for all $t$,
\begin{align*}
\nor{f(t)}{\He^M H^{N}(\cro{x_{1}}^2)}\lesssim \nor{f_0}{\He^M H^{N}(\cro{x_{1}}^2)}+\sqrt{t}\nor{f(t)}{\He^M H^{N}(\cro{x_{1}}^2)}^2.
\end{align*}
\end{prop}

\begin{pr}
The proof only consists in summing over $m$ and $n$ the inequality proven in Proposition \ref{bilinear-estimate-resonant-system}. We first write
\begin{align*}
\nor{\tld{f}_{p,\pm}}{H^{N}(\cro{x_{1}}^2)}
&\lesssim \Bigg\| \inte{0}{t}{\cro{\xi}^N
\sum_{\substack{m,n\in\Z\\ \alpha,\beta\in\{\pm 1\}\\ m\neq n\text{ or }\alpha\neq -\beta\\m\leq n\\ (\ref{cond-res})\text{ satisfied}}} \mathcal{M}(m,n,p)\frac{1}{\sqrt{s|\partial^{2}_{\eta}\phi(\xi,\eta_{0})|}}\frac{\widetilde{f}_{\alpha,m}(\lambda^{\alpha,\beta}_{m,n}\xi)}{\cro{\lambda^{\alpha,\beta}_{m,n}\xi}_{m}} \frac{\widetilde{f}_{\beta,n}((1-\lambda^{\alpha,\beta}_{m,n})\xi)}{\cro{(1-\lambda^{\alpha,\beta}_{m,n})\xi}_{n}} 
}{s}\Bigg\| \\
&\lesssim\inte{0}{t}{\cro{\xi}^N
\sum_{\substack{m,n\in\Z\\ \alpha,\beta\in\{\pm 1\}\\ m\neq n\text{ or }\alpha\neq -\beta\\m\leq n\\ (\ref{cond-res})\text{ satisfied}}} \mathcal{M}(m,n,p)\norlp{\frac{1}{\sqrt{s|\partial^{2}_{\eta}\phi(\xi,\eta_{0})|}}\frac{\widetilde{f}_{\alpha,m}(\lambda^{\alpha,\beta}_{m,n}\xi)}{\cro{\lambda^{\alpha,\beta}_{m,n}\xi}_{m}} \frac{\widetilde{f}_{\beta,n}((1-\lambda^{\alpha,\beta}_{m,n})\xi)}{\cro{(1-\lambda^{\alpha,\beta}_{m,n})\xi}_{n}}}{2}{\xi} 
}{s}.
\end{align*}
Then by Proposition \ref{bilinear-estimate-resonant-system},
\begin{align*}
\nor{\tld{f}_{p,\pm}}{H^{N}(\cro{x_{1}}^2)}&\lesssim
\inte{0}{t}{\!\!\!\!\! \sum_{\substack{m,n\in\Z\\ \alpha,\beta\in\{\pm 1\}\\ m\neq n\text{ or }\alpha\neq -\beta\\m\leq n\\ (\ref{cond-res})\text{ satisfied}}} \!\!\!\!\! \mathcal{M}(m,n,p)\frac{1}{\sqrt{s}}m^{\frac{9}{2}}\frac{1}{\sqrt{mn}}\nor{f_{\alpha,m}}{H^{N}(\cro{x_{1}}^2)}\nor{f_{\beta,n}}{H^{N}(\cro{x_{1}}^2)}}{s}.
\end{align*}
Then, using $\nor{f_{\alpha,m}}{H^{N}(\cro{x_{1}}^2)}\leq m^{-M}\nor{f}{\He^{M}H^{N}(\cro{x_{1}}^2)}$ and integrating lead to
\begin{align*}
p^M \nor{\tld{f}_{p,\pm}}{H^{N}(\cro{x_{1}}^2)}&\lesssim \sqrt{t}\nor{f}{\He^{M}H^{N}(\cro{x_{1}}^2)}^2 p^M\sum_{\substack{m,n\in\Z\\ \alpha,\beta\in\{\pm 1\}\\ m\neq n\text{ or }\alpha\neq -\beta\\m\leq n\\ (\ref{cond-res})\text{ satisfied}}} \mathcal{M}(m,n,p) m^{4-M}n^{\frac{1}{2}-M}.
\end{align*}
Since $M>6$, we are in the framework of the half resummation Theorem \ref{resummation}-(\ref{halfresummation}): there exists a sequence $(u_p(t))_{p\in\N}$ in $\ell^2$ such that
\begin{equation*}
p^M\nor{\tld{f}_{p,\pm}(t)}{H^{N}(\cro{x_{1}}^2)}\lesssim \sqrt{t}\nor{f(t)}{\He^{M}H^{N}(\cro{x_{1}}^2)}^2 u_p(t).
\end{equation*}
This proves Proposition \ref{fixed-point-resonant}.
\end{pr}

\section{Validity of the approximation}\label{section-validity-approx}
Let $f$ be a solution to the initial system with initial data $f_0$. Let $g$ be a solution to the resonant system with the same initial data. Our aim is to estimate the difference $h:=f-g$. We are going to prove the approximation theorem \ref{approximation theorem}.

\subsection{Duhamel formula for $h$}

So as to clarify the notations, let us call $D_1$ and $D_2$ the two following bilinear forms:
\begin{enumerate}
\item the bilinear operator $D_1$ corresponds to the Duhamel formula for the original equation:
\begin{align*}
D_{1}(a,b):=\sum_{m,n}\sum_{\alpha,\beta=\pm 1}\alpha\beta\mathcal{M}(n,m,p)\inte{0}{t}{  \inte{\R}{}{  e^{\mp is\phi^{\alpha,\beta}_{m,n,p}}\frac{\widetilde{a}_{\alpha,m}(\eta)}{\cro{\eta}_{m}} \frac{\widetilde{b}_{\beta,n}(\xi-\eta)}{\cro{\xi-\eta}_{n}}  }{\eta}  }{s},
\end{align*}
with $\mathcal{M}(m,n,p)$ is the interaction term between three Hermite functions (\ref{def-interaction-term}).
\item the bilinear operator $D_2$ corresponds to the resonant equation:
\begin{align*}
D_{2}(a,b):=\inte{0}{t}{
\sum_{\substack{m,n\in\Z\\ \alpha,\beta\in\{\pm 1\}\\ m\neq n\text{ or }\alpha\neq -\beta\\ (\ref{cond-res})\text{ satisfied}}} \frac{C_{sp}}{\sqrt{t|\partial^{2}_{\eta}\phi(\xi,\eta_{0})|}}\frac{\widetilde{a}_{\alpha,m}(\lambda^{\alpha,\beta}_{m,n}\xi)}{\cro{\lambda^{\alpha,\beta}_{m,n}\xi}_{m}} \frac{\widetilde{b}_{\beta,n}((1-\lambda^{\alpha,\beta}_{m,n})\xi)}{\cro{(1-\lambda^{\alpha,\beta}_{m,n})\xi}_{n}} 
}{s},
\end{align*}
with $C_{sp}$ the constant occuring in the Stationary Phase Lemma (Proposition \ref{stat-phase}) and $\lambda^{\alpha,\beta}_{m,n}=\left(1+\alpha\beta\sqrt{\frac{n+1}{m+1}}\right)^{-1}$.
\end{enumerate}

The Duhamel formula for $h$ can be written as follows:
\begin{align*}
h&=f-g\\
&=D_{1}(f,f)-D_{2}(g,g)\\
&=D_{1}(g+h,g+h)-D_{2}(g,g)\\
&=D_{1}(h,h)+2D_{1}(g,h)+(D_{1}-D_{2})(g,g).
\end{align*}
Our aim is to estimate the $S^{M,N}_{t}$ norm of $h$ in terms of $t$ and $\eps$ (the size of $g_0$). So as to do it, we will first establish a differential inequality and then use Gronwall's Lemma.

\begin{lem}\label{lemma for approximation theorem}
Let $\mathcal{N}(t)$ and $\mathcal{M}(t)$ be the $\He^{M_{0}}L^2$ and $S^{M,N}_{t}$ norm of $h(t)$. Then
\begin{align*}
\mathcal{N}(t)\lesssim \inte{0}{t}{\left(s^{\frac{1}{2}+\omega}\mathcal{N}(s)\mathcal{M}(s)^2+s^{-\frac{1}{4}+\omega} \mathcal{N}(s)\mathcal{M}(s)+s^{-\frac{1}{4}+\omega} \eps \mathcal{N}(s)+\cro{s}^{-1}\eps^2+s^{-\frac{1}{2}}\eps^3\right)}{s}.
\end{align*}
\end{lem}
Before going through the proof of the lemma, let us prove the approximation theorem.\\\\
\begin{prf}{Theorem}{\ref{approximation theorem}}
First of all, whenever $t\leq C\eps^{-\frac{4}{3(1+\omega)}}$, we know by Theorems \ref{thm} and \ref{thm-existence-resonant-system} that 
\begin{align*}
\mathcal{M}(t)\leq\eps.
\end{align*}
Then the previous inequality can be rewritten as a Gronwall inequality
\begin{align*}
\mathcal{N}(t)\leq K \left(\inte{0}{t}{s^{\frac{1}{2}+\omega}\eps^{2}\mathcal{N}(s)+s^{-\frac{1}{4}+\omega} \eps \mathcal{N}(s)+s^{-\frac{1}{2}+\omega}\eps^2+s^{\omega}s^{\frac{1}{4}} \eps^3}{s}\right),
\end{align*}
which can be simplified again, by using $\eps\leq Ct^{-\frac{3(1+\omega)}{4}}$:
\begin{align*}
s^{\frac{1}{2}+\omega}\eps^{2}&\leq s^{\frac{1}{2}+\omega}Ct^{-\frac{3(1+\omega)}{4}}\eps\leq Cs^{-\frac{1}{4}}s^{-\frac{\omega}{4}}\\
s^{\omega}s^{\frac{1}{4}} \eps^3&\leq s^{-\frac{1}{2}+\omega}\eps^2,
\end{align*}
which leads to 
\begin{align*}
\mathcal{N}(t)&\leq CK \inte{0}{t}{s^{-\frac{1}{4}+\omega} \eps \mathcal{N}(s)+\cro{s}^{-1}\eps^2}{s}\\
&\leq CK \ln\cro{t}\eps^2+CK \inte{0}{t}{s^{-\frac{1}{4}+\omega} \eps \mathcal{N}(s)}{s}.
\end{align*}
Then, Gronwall's lemma gives
\begin{align*}
\mathcal{N}(t)&\leq CK \ln\cro{t}\eps^2+\inte{0}{t}{CK \ln\cro{s}\eps^2 CKs^{-\frac{1}{4}+\omega} \eps \exp\left(CK\inte{s}{t}{r^{-\frac{1}{4}+\omega} \eps}{r}\right)}{s}\\
&\leq CK \ln\cro{t}\eps^2 +(CK)^2\eps^3\inte{0}{t}{\ln\cro{s}s^{-\frac{1}{4}+2\omega}\exp\left({CK\eps(t^{\frac{3}{4}+\omega}-s^{\frac{3}{4}+\omega})}\right)}{s}.
\end{align*}
Whenever $t\leq C \eps^{-\frac{4}{3+\omega}}$, we have
\begin{align*}
&CK \ln\cro{t}\eps^2\leq \frac{CK}{2}\ln(1+C^2\eps^{-\frac{8}{3+\omega}})\eps^2\lesssim \eps^\alpha,~~\forall\alpha<2,\\
&\exp\left({CK\eps(t^{\frac{3}{4}+\omega}-s^{\frac{3}{4}+\omega})}\right)\leq \exp\left(CK\eps^2\right),\\
&(CK)^2\eps^3\inte{0}{t}{\ln\cro{s}s^{-\frac{1}{4}+2\omega}}{s}\leq (CK)^2\eps^3t^{\frac{3}{4}+2\omega}\leq (CK)^2\eps^3\eps^{-\frac{3+8\omega}{3+\omega}}.
\end{align*}
This proves that for all $\alpha<2$, for all $\omega$ such that $3-\frac{3+8\omega}{3+\omega}>\alpha$, i.e. $\omega<3\frac{2-\alpha}{5+\alpha}$ there exists a $C(\alpha,\omega)$ such that, for $\eps$ small enough, for all $t\leq C(\alpha,\omega)\eps^{-\frac{4}{3+\omega}}$, we have
\begin{align*}
\mathcal{N}(t)\leq \eps^{\alpha}.
\end{align*}
This proves the theorem.
\end{prf}
\\
\\
\begin{prf}{Lemma}{\ref{lemma for approximation theorem}}
First of all, it has been proven in Theorem \ref{prop-contraction} that the operator $D_1$ satisfies the following inequality:
\begin{equation}\label{validity-approx-ineg-D1}
\begin{aligned}
\nor{D_{1}(a,b)}{S^{M,N}_{t}}\lesssim &\inte{0}{t}{s^{\omega}\cro{s}^{-\frac{1}{4}}\nor{a}{L^2}\nor{b}{S^{M,N}_{s}}+s^{\omega}\cro{s}^{\frac{1}{2}}\nor{a}{L^2}\nor{b}{S^{M,N}_{s}}^2}{s} \\
&+ t^{\omega+\frac{1}{4}}\left(\nor{a(t)}{L^2}\nor{b(t)}{S^{M,N}_{t}}+\nor{a}{L^2}\nor{b}{S^{M,N}_{1}}\right)
\end{aligned}
\end{equation}
This allows to bound the terms $D_{1}(h,h)$ and $2D_{1}(g,h)$ involved in the Duhamel formula for $h$. The rest of the proof is devoted to the bounds for the remainder term, i.e.
\begin{align*}
D_{1}(g,g)-D_{2}(g,g).
\end{align*}

In order to estimate the term $D_{1}(g,g)-D_{2}(g,g)$, we are going to write it according to the heuristics done in Section \ref{derivation-resonant-system}: we are writing
\begin{align*}
D_{1}(g,g)-D_{2}(g,g)=SI_{\pm,p}(g)+NR_{\pm,p}(g)+Osc_{\pm,p}(g),
\end{align*}
where
\begin{enumerate}
\item $SI_{\pm,p}(g)$ is the \emph{self-interacting} term, corresponding to the Hermite modes giving birth to a non space resonant mode for all $\xi\neq 0$:
\begin{align*}
SI_{\pm,p}(g):=\sum_{m}\sum_{\alpha=\pm 1,~\beta=-\alpha}\alpha\beta\mathcal{M}(m,m,p)\inte{0}{t}{  \inte{\R}{}{  e^{\mp is\phi^{\alpha,\beta}_{m,m,p}}\frac{\widetilde{g}_{\alpha,m}(\eta)}{\cro{\eta}_{m}} \frac{\widetilde{g}_{\beta,m}(\xi-\eta)}{\cro{\xi-\eta}_{m}}  }{\eta}  }{s}.
\end{align*}
\item $NR_{\pm,p}(g)$ corresponds to the stationary phase remainder:
\begin{align*}
NR_{\pm,p}(g):= \sum_{\substack{m,n\in\N\\\alpha,\beta\in\{\pm 1\} \\m\neq n\text{ or }\alpha\neq - \beta}}\alpha\beta\mathcal{M}(n,m,p)NR^{\alpha,\beta}_{m,n}(t,\xi),
\end{align*}
with
\begin{align*}
NR^{\alpha,\beta}_{m,n}(t,\xi):=\int_{0}^{t}  &\left[ \int_{\R} e^{\mp is\phi^{\alpha,\beta}_{m,n,p}}\frac{\widetilde{f}_{\alpha,m}(\eta)}{\cro{\eta}_{m}} \frac{\widetilde{f}_{\beta,n}(\xi-\eta)}{\cro{\xi-\eta}_{n}}  d\eta\right.\\
&\left.  -\frac{Ce^{\mp is\phi^{\alpha,\beta}_{m,n,p}(\xi,\eta_{0}(\xi))}}{\sqrt{t|\partial^{2}_{\eta}\phi(\xi,\eta_{0})|}}\frac{\widetilde{f}_{\alpha,m}(\eta_{0}(\xi))}{\cro{\eta_{0}(\xi)}_{m}} \frac{\widetilde{f}_{\beta,n}(\xi-\eta_{0}(\xi))}{\cro{\xi-\eta_{0}(\xi)}_{n}}\right]ds.
\end{align*}
\item $Osc_{\pm,p}(g)$ corresponds to the modes giving birth to time resonances, i.e.
\begin{align*}
Osc_{\pm,p}(g):=\sum_{\substack{m,n\in\N\\\alpha,\beta\in\{\pm 1\} \\m\neq n\text{ or }\alpha\neq - \beta \\(\ref{cond-res})\text{ satisfied}}}\alpha\beta\mathcal{M}(n,m,p)\inte{0}{t}{  \inte{\R}{}{  e^{\mp is\phi^{\alpha,\beta}_{m,n,p}}\frac{\widetilde{f}_{\alpha,m}(\eta)}{\cro{\eta}_{m}} \frac{\widetilde{f}_{\beta,n}(\xi-\eta)}{\cro{\xi-\eta}_{n}}  }{\eta}  }{s}.
\end{align*}
\end{enumerate}
The next three sections will be dedicated to bounding those three terms, and more precisely proving the following inequality
\begin{align}\label{validity-approx-ineg-D1-D2}
\norlp{D_{1}(g,g)-D_{2}(g,g)}{2}{}\lesssim \inte{0}{t}{\left(\frac{1}{\cro{s}}\eps^2+\frac{1}{\sqrt{s}}\eps^3\right)}{s}.
\end{align}
Then combining (\ref{validity-approx-ineg-D1}) and (\ref{validity-approx-ineg-D1-D2}) give Lemma \ref{lemma for approximation theorem}.
\subsection{Estimates for the self-interaction remainder}

We are going to prove the following lemma:
\begin{lem}\label{estimate-self-interaction-remainder}
There exists a sequence $(u_{p}(s))_{p\in\N}$ in the unit ball of $\ell^2$ such that 
\begin{align*}
p^{M_0}\norlp{SI_{\pm,p}(g)}{2}{}\lesssim \inte{0}{t}{u_{p}(s)\frac{\eps^2}{\cro{s}}}{s}.
\end{align*}
\end{lem}

\begin{prf}{Lemma}{\ref{estimate-self-interaction-remainder}}
We want to estimate the $L^2$ norm of
\begin{align*}
\sum_{m}\sum_{\alpha=\pm 1,~\beta=-\alpha}\alpha\beta\mathcal{M}(m,m,p)\inte{0}{t}{  \inte{\R}{}{  e^{\mp is\phi^{\alpha,\beta}_{m,m,p}}\frac{\widetilde{g}_{\alpha,m}(\eta)}{\cro{\eta}_{m}} \frac{\widetilde{g}_{\beta,m}(\xi-\eta)}{\cro{\xi-\eta}_{m}}  }{\eta}  }{s}.
\end{align*}
We know that in this case the quantity $\phi^{\alpha,\beta}_{m,m,p}(\xi,\eta_0(\xi))$ never vanishes except when $\xi=0$. This is the reason why we will handle separately the zones around and outside the origin.

Let $\chi$ be a smooth function, compactly supported, which is equal to $1$ on $\left[-\frac{1}{2},\frac{1}{2}\right]$ and $0$ outside $[-1,1]$. Then we write

\begin{align*}
\sum_{m}\sum_{\alpha=\pm 1,~\beta=-\alpha}\alpha\beta\mathcal{M}(m,m,p)\inte{0}{t}{  \inte{\R}{}{  e^{\mp is\phi^{\alpha,\beta}_{m,m,p}}\frac{\widetilde{g}_{\alpha,m}(\eta)}{\cro{\eta}_{m}} \frac{\widetilde{g}_{\beta,m}(\xi-\eta)}{\cro{\xi-\eta}_{m}}  }{\eta}  }{s} = SI_{s}+SI_{l},
\end{align*}
with $SI_{s}$ corresponding to the small values of $\xi$
\begin{align*}
SI_{s}:=\chi(\xi)\sum_{m}\sum_{\alpha=\pm 1,~\beta=-\alpha}\alpha\beta\mathcal{M}(m,m,p)\inte{0}{t}{  \inte{\R}{}{  e^{\mp is\phi^{\alpha,\beta}_{m,m,p}}\frac{\widetilde{g}_{\alpha,m}(\eta)}{\cro{\eta}_{m}} \frac{\widetilde{g}_{\beta,m}(\xi-\eta)}{\cro{\xi-\eta}_{m}}  }{\eta}  }{s},
\end{align*}
and $SI_{l}$ corresponding to the large ones. We will use two different strategies for these integrals : $SI_{s}$ will be bounded by using the time resonances method, $SI_{l}$ by a stationary phase.

\subsubsection{Study of $SI_{s}$.}

In the zone $|\xi|\leq 1$, the phase $\phi^{\alpha,\beta}_{m,m,p}$ does not vanish and is easily bounded. In fact, since $\beta=-\alpha$,
\begin{align*}
\phi^{\alpha,\beta}_{m,m,p}=\sqrt{\xi^2+2p+2}+\alpha\left(\sqrt{\eta^2+2m+2}-\sqrt{(\xi-\eta)^2+2m+2}\right).
\end{align*}
Since $|\xi|\leq 1$, 
\begin{align*}
\left|\sqrt{\eta^2+2m+2}-\sqrt{(\xi-\eta)^2+2m+2}\right|&\leq |\xi|\sup_{\eta}\frac{|\eta|}{\sqrt{\eta^2+2m+2}}\\
&\leq 1.
\end{align*}
Then 
\begin{align*}
|\phi^{\alpha,\beta}_{m,m,p}|\geq \sqrt{2p+2}-1\geq 1.
\end{align*} 

Then it is possible to use the time resonances method, applied in page \pageref{time-resonances-method} for example. Writing $e^{\mp is\phi^{\alpha,\beta}_{m,m,p}}=\frac{1}{i\phi^{\alpha,\beta}_{m,m,p}}\partial_{s}e^{\mp is\phi^{\alpha,\beta}_{m,m,p}}$, we get, by integrating by parts,

\begin{align*}
SI_{s}= SI_{s}^{1}+SI_{s}^2+SI_{s}^{3},
\end{align*}
with
\begin{itemize}
\item $SI_{s}^1$ is the boundary term
\begin{equation*}
\begin{aligned}
SI_{s}^1:=\chi(\xi)  \sum_{m}\sum_{\substack{\alpha=\pm 1\\ \beta=-\alpha}}\alpha\beta\mathcal{M}(m,m,p)\int_{\R}\bigg( &\frac{e^{\mp it\phi^{\alpha,\beta}_{m,m,p}}}{\phi^{\alpha,\beta}_{m,m,p}}\frac{\widetilde{g}_{\alpha,m}(t,\eta)}{\cro{\eta}_{m}} \frac{\widetilde{g}_{\beta,m}(t,\xi-\eta)}{\cro{\xi-\eta}_{m}}\\
-&\frac{1}{\phi^{\alpha,\beta}_{m,m,p}}\frac{\widetilde{g}_{\alpha,m}(0,\eta)}{\cro{\eta}_{m}} \frac{\widetilde{g}_{\beta,m}(0,\xi-\eta)}{\cro{\xi-\eta}_{m}} \bigg) d\eta,
\end{aligned}
\end{equation*}
\item $SI_{s}^2$ and $SI_{s}^3$ are the two terms involving time derivatives of $g$:
\begin{align*}
SI_{s}^2&:=\chi(\xi)  \sum_{m}\sum_{\substack{\alpha=\pm 1\\ \beta=-\alpha}}\alpha\beta\mathcal{M}(m,m,p)\inte{0}{t}{  \inte{\R}{}{  \frac{e^{\mp is\phi^{\alpha,\beta}_{m,m,p}}}{\phi^{\alpha,\beta}_{m,m,p}}\frac{\partial_{s}\widetilde{g}_{\alpha,m}(\eta)}{\cro{\eta}_{m}} \frac{\widetilde{g}_{\beta,m}(\xi-\eta)}{\cro{\xi-\eta}_{m}}  }{\eta}  }{s},\\
SI_{s}^3&:=\chi(\xi)  \sum_{m}\sum_{\substack{\alpha=\pm 1\\ \beta=-\alpha}}\alpha\beta\mathcal{M}(m,m,p)\inte{0}{t}{  \inte{\R}{}{ \frac{e^{\mp is\phi^{\alpha,\beta}_{m,m,p}}}{\phi^{\alpha,\beta}_{m,m,p}} \frac{\widetilde{g}_{\alpha,m}(\eta)}{\cro{\eta}_{m}} \frac{\partial_{s}\widetilde{g}_{\beta,m}(\xi-\eta)}{\cro{\xi-\eta}_{m}}  }{\eta}  }{s}.
\end{align*}

\end{itemize}

\paragraph{Estimates for $SI_{s}^1$.} Remarking that $\frac{1}{\phi^{\alpha,\beta}_{m,m,p}}$ is a Coifman-Meyer multiplier, we have the following bound:
\begin{align*}
\norlp{SI_{s}^{1}}{2}{\xi}\lesssim \sum_{m}m^{-2M-1}\mathcal{M}(m,m,p)\sum_{\alpha=-\beta=\pm 1}\alpha\beta  &\left(\nor{g_{\alpha,m}(t)}{L^{\infty}}\nor{g_{\beta,n}(t)}{L^{2}}\right.\\
&\left.+\nor{g_{\alpha,m}(0)}{L^{\infty}}\nor{g_{\beta,n}(0)}{L^{2}}\right),
\end{align*}
i.e.
\begin{align*}
\norlp{SI_{s}^{1}}{2}{\xi}\lesssim \sum_{m}m^{-2M-1}\mathcal{M}(m,m,p)\eps^2.
\end{align*}
Resumming in $p$ is then not a problem since for all $\nu>1/8$ and $\varpi<1/24$,
\begin{align*}
\mathcal{M}(m,m,p)\leq C_{K}\frac{m^{\nu}}{p^{\varpi}} \frac{m^K}{p^K}.
\end{align*}
Hence
\begin{align}\label{valapp-SIs-1}
\norlp{SI_{s}^{1}}{2}{\xi}&\lesssim \eps^2 u_{p}(t),
\end{align}
with $(u_{p}(t))_{p\in\Z}$ in the unit ball of $\ell^2$.
\paragraph{Estimates for $SI_{s}^2$.} Here we take advantage of the fact that $\partial_{s}\tld{g}_{\alpha,m}$ is quadratic in $g$: more precisely,
\begin{align*}
\partial_{s}g_{\alpha,m}=e^{is\cro{D}_{m}}\left(e^{-is\cro{D}_m}g_{\alpha,m}\right)^2,
\end{align*}
which will lead to, using the dispersion inequality (\ref{dispequ}),
\begin{align}\label{valapp-SIs-2}
\norlp{SI_{s}^{2}}{2}{\xi}\lesssim \sum_{m}\sum_{\alpha=-\beta=\pm 1}\alpha\beta \mathcal{M}(m,m,p) \inte{0}{t}{\frac{1}{\cro{s}}\eps^3 m^{-3M-1}}{s}.
\end{align}
Resumming is not a problem either. Moreover, the third integral, $SI_{s}^{3}$, can be bounded in the same way.

\subsubsection{Study of $SI_{l}$.}

If $|\xi|\geq 1/2$ then $\partial_{\eta}\phi^{\alpha,\beta}_{m,m,p}$ never vanishes, so the stationary phase Lemma applies.\\
For $\xi$ fixed and nonzero, the minimum of $\partial_{\eta}\phi(\xi,\eta)$ is reached in $\eta=\xi/2$ and equals
\begin{align*}
\frac{\xi/2}{(\xi^2/4+2m+2)^{\frac{1}{2}}}\geq \frac{1}{4\sqrt{2m+2}}.
\end{align*}

In order to be able to apply Proposition \ref{stat-phase} we need to localize in $\eta$: let $(\psi_{j})_{j\in\Z}$ be a family of functions such that each $\psi_j$ is supported in the annulus $2^{j}\leq |\eta|\leq 2^{j+1}$ and such that $\sum_{j}\psi_j=1$. Let us write $SI_{l}^{j,p}(\xi)$ for the following integral:
\begin{align*}
SI_{l}^{j,p}(\xi):=(1-\chi(\xi))\inte{0}{t}{\inte{\R}{}{\psi_{j}(\eta)   e^{\mp is\phi^{\alpha,\beta}_{m,m,p}}\frac{\widetilde{g}_{\alpha,m}(\eta)}{\cro{\eta}_{m}} \frac{\widetilde{g}_{\beta,m}(\xi-\eta)}{\cro{\xi-\eta}_{m}}  }{\eta}}{s}.
\end{align*}
 Then Proposition \ref{stat-phase} applies with $F_{j}(\xi,\eta):=\psi_{j}(\eta)  \frac{\widetilde{g}_{\alpha,m}(\eta)}{\cro{\eta}_{m}} \frac{\widetilde{g}_{\beta,m}(\xi-\eta)}{\cro{\xi-\eta}_{m}}$:

\begin{align*}
SI_{l}^{j,p}(\xi)\lesssim \inte{0}{t}{\frac{1}{\cro{s}}2^{\frac{j}{2}}\sqrt{m}\left(\norlp{F_j}{2}{\eta}(\xi)+\norlp{\partial_{\eta}F_j}{2}{\eta}(\xi)\right)}{s}.
\end{align*}

We want to take the $L^2_{\xi}$ norm of $SI_{l}^{j,p}$. First we know that if $f$ and $g$ are two functions in $L^2$ such that for all $\xi$, $\eta\mapsto f(\eta)g(\xi-\eta)$ is in $L^2$, then 
\begin{align*}
\nor{f(\eta)g(\xi-\eta)}{L^{2}_{\xi,\eta}}=\norlp{f}{2}{}\norlp{g}{2}{}.
\end{align*}
Then, remaking that
\begin{align*}
|F_{j}|\leq \left|\psi_{j}(\eta)  \frac{\widetilde{g}_{\alpha,m}(\eta)}{\cro{\eta}_{m}} \frac{\widetilde{g}_{\beta,m}(\xi-\eta)}{\cro{\xi-\eta}_{m}}\right|+\left|\psi_{j}(\eta) \frac{\widetilde{g}_{\alpha,m}(\eta)}{\cro{\eta}_{m}} \frac{\widetilde{g}_{\beta,m}(\xi-\eta)}{\cro{\xi-\eta}_{m}}\right|
\end{align*}
implies
\begin{align*}
\norlp{F_j}{2}{\xi,\eta}\lesssim \norlp{\psi_{j}(\eta)  \frac{\widetilde{g}_{\alpha,m}(\eta)}{\cro{\eta}_{m}}}{2}{\eta}\norlp{\frac{\widetilde{g}_{\beta,m}(\eta)}{\cro{\eta}_{m}}}{2}{\eta}+\norlp{\psi_{j}(\eta) \frac{\widetilde{g}_{\alpha,m}(\eta)}{\cro{\eta}_{m}}}{2}{\eta}\norlp{ \frac{\widetilde{g}_{\beta,m}(\eta)}{\cro{\eta}_{m}}}{2}{\eta}.
\end{align*}
Finally, since $\left|\psi_{j}(\eta)\frac{1}{\cro{\eta}_{m}}\right|\leq \min(1,2^{-j})$, we obtain

\begin{align*}
\norlp{F_j}{2}{\xi,\eta}&\lesssim \min(1,2^{-j})\nor{g_{\alpha,m}}{L^2}\nor{g_{\beta,m}}{L^2}\\
&\lesssim \min(1,2^{-j})m^{-2M}\eps^2.
\end{align*}

In the same fashion we obtain 
\begin{align*}
\norlp{\norlp{\partial_{\eta}F_j}{2}{\eta}}{2}{\xi} \lesssim \min(1,2^{-j})m^{-2M}\eps^2,
\end{align*}
which leads to
\begin{align*}
\norlp{SI_{l}^{j,p}}{2}{\xi}\lesssim \inte{0}{t}{\frac{1}{\cro{s}}2^{\frac{j}{2}}\sqrt{m}\min(1,2^{-j})m^{-2M}\eps^2}{s}.
\end{align*}

Since $2^{\frac{j}{2}}\sqrt{m}\min(1,2^{-j})$ is summable over $\Z$, we can sum $\norlp{SI_{l}^{j,p}}{2}{\xi}$ over $\Z$ and obtain

\begin{align}\label{valapp-SIl}
\sum_{j\in\Z}\norlp{SI_{l}^{j,p}}{2}{\xi}\lesssim \inte{0}{t}{\frac{u_{p}(s)}{\cro{s}}\eps^2}{s},
\end{align}
with $(u_{p}(s))_{p\in\Z}$ in the unit ball of $\ell^2$.

Finally, Inequalities (\ref{valapp-SIs-1}), (\ref{valapp-SIs-2}) and (\ref{valapp-SIl}) end the proof of Lemma \ref{estimate-self-interaction-remainder}.
\end{prf}

\subsection{Estimates for the non-stationary remainder}

\begin{lem}\label{spr-estimate}
There exists a sequence $(u_p(s))_{p\in\N}$ such that
 \begin{align*}
p^{M_0}\norlp{NR_{\pm,p}(g)}{2}{}\lesssim a_p\frac{\eps^2}{t^{\frac{3}{4}}},
 \end{align*}
 where $(a_p)_{p\in\N}$ is a sequence in the unit ball of $\ell^2$.
\end{lem}

The proof of this lemma relies on the stationary phase Proposition \ref{stat-phase}.\\

\begin{prf}{Lemma}{\ref{spr-estimate}}

The integral term we want to use Proposition \ref{stat-phase} on is the one occuring in Duhamel formula, that is to say:
\begin{align*}
D_{1}(g,g)-SI(g)&:=\sum_{\substack{m,n\in\Z\\ \alpha,\beta\in\{\pm 1\}\\ m\neq n\text{ or }\alpha\neq -\beta\\ (\ref{cond-res})\text{ satisfied}}} \inte{0}{t}{  \inte{\R}{}{  e^{\mp is\phi^{\alpha,\beta}_{m,n,p}}\frac{\widetilde{g}_{\alpha,m}(\eta)}{\cro{\eta}_{m}} \frac{\widetilde{g}_{\beta,n}(\xi-\eta)}{\cro{\xi-\eta}_{n}}  }{\eta}  }{s}\\
&=\sum_{\substack{m,n\in\Z\\ \alpha,\beta\in\{\pm 1\}\\ m\neq n\text{ or }\alpha\neq -\beta\\ (\ref{cond-res})\text{ satisfied}}} I_{m,n}^{\alpha,\beta}(t,\xi).
\end{align*}
Here $\psi(\eta):=\phi^{\alpha,\beta}_{m,n,p}(\xi,\eta)$, the critical point is $\eta_{0}=\lambda^{\alpha,\beta}_{m,n}\xi$ and $F(\xi,\eta):=\frac{\widetilde{g}_{\alpha,m}(\eta)}{\cro{\eta}_{m}} \frac{\widetilde{g}_{\beta,n}(\xi-\eta)}{\cro{\xi-\eta}_{n}}$. Let $\chi_{\rho}\in\mathcal{C}^{\infty}_{0}$ equal to zero on $B(0,\rho)^{c}$.\\
In order to apply Proposition \ref{stat-phase}, we have to find 
\begin{itemize}
\item either an lower bound for $|\psi''|$,
\item or an upper bound for $\frac{\sqrt{|\psi(\eta)|}}{|\psi'(\eta)|}$.
\end{itemize}   Since in some cases (when $\alpha\beta=-1$) $\psi''(\eta)=\partial_{\eta}^{2}\phi$ can vanish, it is better to try to bound directly 
\begin{equation*}
\frac{\sqrt{|\psi(\eta)|}}{|\psi'(\eta)|},
\end{equation*} or rather 
\begin{equation*}
\frac{|\psi(\eta)-\psi(\eta_{0})|}{\psi'(\eta)^2}=\frac{|\phi^{\alpha,\beta}_{m,n,p}(\xi,\eta)-\phi^{\alpha,\beta}_{m,n,p}(\xi,\lambda^{\alpha,\beta}_{m,n}\xi)|}{\left(\partial_{\eta}\phi^{\alpha,\beta}_{m,n,p}(\xi,\eta)\right)^{2}}\cdot
\end{equation*}

The denominator vanishes at infinity or at $\lambda^{\alpha,\beta}_{m,n}\xi$. \\
Since $m\neq n$, 
\begin{equation*}
\partial_{\eta}^{2}\phi^{\alpha,\beta}_{m,n}(\xi,\lambda^{\alpha,\beta}_{m,n}\xi)\neq 0.
\end{equation*}
Then $\frac{|\psi(\eta)-\psi(\eta_{0})|}{\psi'(\eta)^2}$ is well-defined at the point $\lambda^{\alpha,\beta}_{m,n}\xi$ and
\begin{equation*}
\frac{|\psi(\lambda^{\alpha,\beta}_{m,n}\xi)-\psi(\eta_{0})|}{\psi'(\lambda^{\alpha,\beta}_{m,n}\xi)^2}=\frac{1}{|\partial_{\eta}^ {2}\phi^{\alpha,\beta}_{m,n,p}(\xi,\lambda^{\alpha,\beta}_{m,n}\xi)|}\cdot
\end{equation*}\\
Then, understanding the asymptotic behavior of $\frac{|\psi(\eta)-\psi(\eta_{0})|}{\psi'(\eta)^2}$ will allow us to bound it (for sufficiently large values of $\rho$).

\begin{enumerate}
\item if $\alpha\beta=1$, then 
\begin{equation*}
|\phi^{\alpha,\beta}_{m,n,p}(\xi,\eta)-\phi^{\alpha,\beta}_{m,n,p}(\xi,\lambda^{\alpha,\beta}_{m,n})|\sim_{\eta\rightarrow\infty} \ 2|\eta|
\end{equation*} and 
\begin{equation*}
|\partial_{\eta}\phi^{\alpha,\beta}_{m,n,p}(\xi,\eta)|\rightarrow_{\eta\rightarrow\infty} 2,
\end{equation*} hence
\begin{align*}
\frac{|\phi^{\alpha,\beta}_{m,n,p}(\xi,\eta)-\phi^{\alpha,\beta}_{m,n,p}(\xi,\lambda^{\alpha,\beta}_{m,n}\xi)|}{\left(\partial_{\eta}\phi^{\alpha,\beta}_{m,n,p}(\xi,\eta)\right)^{2}}\sim_{\eta\rightarrow\infty} \frac{\eta}{2}\cdot
\end{align*}
\item if $\alpha\beta=-1$, then
\begin{align*}
\lim_{\eta\rightarrow\infty}|\phi^{\alpha,\beta}_{m,n,p}(\xi,\eta)-\phi^{\alpha,\beta}_{m,n,p}(\xi,\lambda^{\alpha,\beta}_{m,n}\xi)|&=  \left|\alpha\sqrt{\left(\lambda^{\alpha,\beta}_{m,n}\xi\right)^{2}+2m+2}+\beta\sqrt{\left((1-\lambda^{\alpha,\beta}_{m,n})\xi\right)^{2}+2n+2}\right|\\
&=\frac{1}{|\lambda_{m,n}^{\alpha,\beta}|}\sqrt{\left(\lambda^{\alpha,\beta}_{m,n}\xi\right)^{2}+2m+2}\cdot
\end{align*}
and
\begin{align*}
|\partial_{\eta}\phi^{\alpha,\beta}_{m,n,p}(\xi,\eta)|\sim_{\eta\rightarrow\infty} \frac{2|n-m|}{\eta^2}.
\end{align*}
Hence
\begin{align*}
\frac{|\phi^{\alpha,\beta}_{m,n,p}(\xi,\eta)-\phi^{\alpha,\beta}_{m,n,p}(\xi,\lambda^{\alpha,\beta}_{m,n}\xi)|}{\left(\partial_{\eta}\phi^{\alpha,\beta}_{m,n,p}(\xi,\eta)\right)^{2}}&\sim_{\eta\rightarrow\infty} \frac{\eta^4\sqrt{\left(\lambda^{\alpha,\beta}_{m,n}\xi\right)^{2}+2m+2}}{2|(m-n)^2\lambda_{m,n}^{\alpha,\beta}|}\\
&\lesssim \eta^4\sqrt{\xi^2+2(\sqrt{m+1}-\sqrt{n+1})^2}.
\end{align*}
\end{enumerate}
These bounds being established, it remains now to apply the stationary phase Proposition \ref{stat-phase}, with
\begin{enumerate}
\item $\psi(\eta):=\phi_{m,n,p}^ {\alpha,\beta}(\xi,\eta)$,
\item $F(\eta):=\frac{\widetilde{g}_{\alpha,m}(\eta)}{\cro{\eta}_{m}} \frac{\widetilde{g}_{\beta,n}(\xi-\eta)}{\cro{\xi-\eta}_{n}}$,
\item the critical point is $\eta_0=\lambda^{\alpha,\beta}_{m,n}\xi$,
\item $M=2$ since $\partial_{\eta}^{3}\phi=\frac{2(m+1)\eta}{(\eta^2+2m+2)^{5/2}}+\alpha\beta\frac{2(n+1)(\xi-\eta)}{((\xi-\eta)^2+2n+2)^{5/2}}$,
\item the bound $m$ is given by
\begin{align*}
\left\{\begin{array}{ccl}
\frac{1}{\rho}&\text{ if }&\alpha\beta=1,\\
\frac{1}{\rho^4\sqrt{\xi^2+2(\sqrt{m+1}-\sqrt{n+1})^2}}&\text{ if }&\alpha\beta=-1.
\end{array}\right.
\end{align*}
\end{enumerate}
This leads to 

\begin{align*}
I_{m,n}^{\alpha,\beta}(t,\xi)=\frac{Ce^{it\phi(\lambda^{\alpha,\beta}_{m,n} \xi)}}{\partial_{\eta}^2\phi^{\alpha,\beta}_{m,n,p}(\xi,\lambda^{\alpha,\beta}_{m,n}\xi)\sqrt{t}}\chi(\lambda^{\alpha,\beta}_{m,n}\xi)\frac{\widetilde{g}_{\alpha,m}(\lambda^{\alpha,\beta}_{m,n}\xi)}{\cro{\lambda^{\alpha,\beta}_{m,n}\xi}_{m}} \frac{\widetilde{g}_{\beta,n}((1-\lambda^{\alpha,\beta}_{m,n})\xi)}{\cro{(1-\lambda^{\alpha,\beta}_{m,n})\xi}_{n}} +NR_{m,n}^{\alpha,\beta}(t,\xi),
\end{align*}
with $NR_{m,n}^{\alpha,\beta}(t,\xi)$ bounded as follows:
\begin{enumerate}
\item if $\alpha\beta=1$,
\begin{align*}
NR\lesssim \frac{1}{t^{\frac{3}{4}}}\bigg(&\rho^{\frac{7}{4}}\norlp{\frac{\widetilde{g}_{\alpha,m}(\eta)}{\cro{\eta}_{m}} \frac{\widetilde{g}_{\beta,n}(\xi-\eta)}{\cro{\xi-\eta}_{n}}}{2}{\eta}+\rho^{\frac{3}{4}}\norlp{\partial_{\eta}\left(\frac{\widetilde{g}_{\alpha,m}(\eta)}{\cro{\eta}_{m}} \frac{\widetilde{g}_{\beta,n}(\xi-\eta)}{\cro{\xi-\eta}_{n}}\right)}{2}{\eta}\bigg),
\end{align*}
\item if $\alpha\beta=-1$,
\begin{align*}
NR\lesssim &\frac{1}{t^{\frac{3}{4}}}\Bigg(\rho^{7}\left(\xi^2+2(\sqrt{m+1}-\sqrt{n+1})^2\right)^{\frac{7}{8}}\norlp{\frac{\widetilde{g}_{\alpha,m}(\eta)}{\cro{\eta}_{m}} \frac{\widetilde{g}_{\beta,n}(\xi-\eta)}{\cro{\xi-\eta}_{n}}}{2}{\eta}
\\ &+\rho^{3}\left(\xi^2+2(\sqrt{m+1}-\sqrt{n+1})^2\right)^{\frac{3}{8}}\norlp{\partial_{\eta}\left(\frac{\widetilde{g}_{\alpha,m}(\eta)}{\cro{\eta}_{m}} \frac{\widetilde{g}_{\beta,n}(\xi-\eta)}{\cro{\xi-\eta}_{n}}\right)}{2}{\eta}\Bigg),
\end{align*}
\end{enumerate}
Remark that the bound found in the case $\alpha\beta=-1$ is bigger than the one in the case $\alpha\beta=1$.\\ Moreover if we define $\chi$ on an annulus instead of a ball we have 
\begin{align*}
\rho^{k}\norlp{\partial_{\eta}\left(\frac{\widetilde{g}_{\alpha,m}(\eta)}{\cro{\eta}_{m}} \frac{\widetilde{g}_{\beta,n}(\xi-\eta)}{\cro{\xi-\eta}_{n}}\right)}{2}{\eta}\sim\norlp{|\eta|^{k}\partial_{\eta}\left(\frac{\widetilde{g}_{\alpha,m}(\eta)}{\cro{\eta}_{m}} \frac{\widetilde{g}_{\beta,n}(\xi-\eta)}{\cro{\xi-\eta}_{n}}\right)}{2}{\eta}.
\end{align*}
Finally we obtain, in the case $\alpha\beta=-1$:

\begin{align*}
NR_{m,n}^{\alpha,\beta}(t,\xi)\lesssim &\frac{1}{t^{\frac{3}{4}}}\left(\left(\xi^2+2(\sqrt{m+1}-\sqrt{n+1})^2\right)^{\frac{7}{8}}\norlp{|\eta|^{7}\frac{\widetilde{g}_{\alpha,m}(\eta)}{\cro{\eta}_{m}} \frac{\widetilde{g}_{\beta,n}(\xi-\eta)}{\cro{\xi-\eta}_{n}}}{2}{\eta}\right.
\\ &+\left.\left(\xi^2+2(\sqrt{m+1}-\sqrt{n+1})^2\right)^{\frac{3}{8}}\norlp{|\eta|^{3}\partial_{\eta}\left(\frac{\widetilde{g}_{\alpha,m}(\eta)}{\cro{\eta}_{m}} \frac{\widetilde{g}_{\beta,n}(\xi-\eta)}{\cro{\xi-\eta}_{n}}\right)}{2}{\eta}\right).
\end{align*}
This asymptotical bound is also valid for the case $\alpha\beta=1$.\\
We are focusing on the first term of the sum, the first one being even easiest. We can wrtie
\begin{align*}
&\norlp{ \left(\xi^2+2(\sqrt{m+1}-\sqrt{n+1})^2\right)^{\frac{7}{8}}\norlp{|\eta|^{7}\partial_{\eta}\left(\frac{\widetilde{g}_{\alpha,m}(\eta)}{\cro{\eta}_{m}} \frac{\widetilde{g}_{\beta,n}(\xi-\eta)}{\cro{\xi-\eta}_{n}}\right)}{2}{\eta}}{2}{\xi}\\
&\lesssim
\nor{\left(\eta^2+2(\sqrt{m+1}-\sqrt{n+1})^2\right)^{\frac{7}{8}}|\eta|^{7}\partial_{\eta}\left(\frac{\widetilde{g}_{\alpha,m}(\eta)}{\cro{\eta}_{m}} \frac{\widetilde{g}_{\beta,n}(\xi-\eta)}{\cro{\xi-\eta}_{n}}\right)}{L^{2}_{\eta,\xi}(\R^2)} \\
&+\nor{\left((\xi-\eta)^2+2(\sqrt{m+1}-\sqrt{n+1})^2\right)^{\frac{7}{8}}|\eta|^{7}\partial_{\eta}\left(\frac{\widetilde{g}_{\alpha,m}(\eta)}{\cro{\eta}_{m}} \frac{\widetilde{g}_{\beta,n}(\xi-\eta)}{\cro{\xi-\eta}_{n}}\right)}{L^{2}_{\eta,\xi}(\R^2)}.
\end{align*}
We only focus on the first term, which is the ``worst'' in terms of cost of derivatives. By sub-linearity, we are reduced to bound
\begin{align*}
\nor{|\eta|^{7+\frac{7}{4}}\frac{\partial_{\eta}\widetilde{g}_{\alpha,m}(\eta)}{\cro{\eta}_{m}} \frac{\widetilde{g}_{\beta,n}(\xi-\eta)}{\cro{\xi-\eta}_{n}}}{L^{2}_{\eta,\xi}(\R^2)}
+
\nor{(\sqrt{m+1}-\sqrt{n+1})^{\frac{7}{4}}|\eta|^{7}\frac{\partial_{\eta}\widetilde{g}_{\alpha,m}(\eta)}{\cro{\eta}_{m}} \frac{\widetilde{g}_{\beta,n}(\xi-\eta)}{\cro{\xi-\eta}_{n}}}{L^{2}_{\eta,\xi}(\R^2)}.
\end{align*}
First,
\begin{align*}
\nor{|\eta|^{7+\frac{7}{4}}\frac{\partial_{\eta}\widetilde{g}_{\alpha,m}(\eta)}{\cro{\eta}_{m}} \frac{\widetilde{g}_{\beta,n}(\xi-\eta)}{\cro{\xi-\eta}_{n}}}{L^{2}_{\eta,\xi}(\R^2)} &\lesssim
\norlp{|\eta|^{8+\frac{3}{4}}\frac{\partial_{\eta}\widetilde{g}_{\alpha,m}(\eta)}{\cro{\eta}_{m}} }{2}{\eta}\norlp{\frac{\widetilde{g}_{\beta,n}(\eta)}{\cro{\eta}_{n}}}{2}{\eta}\\
&\lesssim \frac{1}{\sqrt{mn}}\norlp{|\eta|^{8+\frac{3}{4}}\partial_{\eta}\widetilde{g}_{\alpha,m}(\eta) }{2}{\eta}\norlp{\widetilde{g}_{\beta,n}(\eta)}{2}{\eta}.
\end{align*}
Then, since $g$ is in $\Sigma_{T}^{M,N}$, with $N>9-1/4$,
\begin{align*}
\norlp{|\eta|^{N+3+\frac{3}{4}}\partial_{\eta}\widetilde{g}_{\alpha,m}(t,\eta) }{2}{\eta} = \eps m^{-M}a_m(t),
\end{align*}
with $(a_m(t))_{m\in\N}$ in the unit ball of $\ell^2$. There exists also $(b_n(t))_{n\in\N}$ in the unit ball of $\ell^2$ such that
\begin{align*}
\nor{|\eta|^{8+\frac{3}{4}}\frac{\partial_{\eta}\widetilde{g}_{\alpha,m}(\eta)}{\cro{\eta}_{m}} \frac{\widetilde{g}_{\beta,n}(\xi-\eta)}{\cro{\xi-\eta}_{n}}}{L^{2}_{\eta,\xi}(\R^2)} \lesssim \frac{1}{\sqrt{mn}}m^{-M}n^{-M}a_m(t) b_n(t) \eps^ 2.
\end{align*}
The summation Theorem \ref{resummation} ends the proof of Lemma \ref{spr-estimate}.\\
Similarly, we have 
\begin{equation*}
\nor{(\sqrt{m+1}-\sqrt{n+1})^{\frac{7}{4}}|\eta|^{7}\frac{\partial_{\eta}\widetilde{g}_{\alpha,m}(\eta)}{\cro{\eta}_{m}} \frac{\widetilde{g}_{\beta,n}(\xi-\eta)}{\cro{\xi-\eta}_{n}}}{L^{2}_{\eta,\xi}(\R^2)} \lesssim \frac{\max(m,n)^{\frac{7}{8}}}{\sqrt{mn}}m^{-M}n^{-M}a_m(t) b_n(t) \eps^ 2,
\end{equation*}
which also fits in the hypotheses of Theorem \ref{resummation}, since we assumed $M>M_{0}+\frac{1}{8}$. This ends the proof of Lemma \ref{spr-estimate}.
\end{prf}

\subsection{Estimates for the oscillating term}
The oscillating term is 

\begin{align*}
\mathrm{Osc}_{\pm,p}(g)(\xi):=\inte{0}{t}{
\sum_{\substack{m,n\in\Z\\ \alpha,\beta\in\{\pm 1\}\\ m\neq n\text{ or }\alpha\neq -\beta\\ (\ref{cond-res})\text{ not satisfied}}} \frac{C_{sp}}{\sqrt{t|\partial^{2}_{\eta}\phi(\xi,\eta_{0})|}}e^{\mp is\phi^{\alpha,\beta}_{m,n,p}(\xi,\lambda^{\alpha,\beta}_{m,n}\xi)}\frac{\widetilde{g}_{\alpha,m}(\lambda^{\alpha,\beta}_{m,n}\xi)}{\cro{\lambda^{\alpha,\beta}_{m,n}\xi}_{m}} \frac{\widetilde{g}_{\beta,n}((1-\lambda^{\alpha,\beta}_{m,n})\xi)}{\cro{(1-\lambda^{\alpha,\beta}_{m,n})\xi}_{n}} 
}{s}.
\end{align*}

\begin{lem}\label{validity-approx-osc}
The $L^2$ norm of oscillating term can be bounded as follows: there exists $(u_{p}(s))_{p\in\N}$ inthe unit ball of $\ell^2$ such that
\begin{align*}
p^{M_{0}}\norlp{\mathrm{Osc}_{\pm,p}(g)(\xi)}{2}{\xi}\lesssim \inte{0}{t}{u_p(s)\frac{1}{s^{3/2}}\eps^2+\frac{1}{\sqrt{s}}\eps^3}{s}.
\end{align*}
\end{lem}

We shall not write the full proof of this proposition. We only recall that if $(C)$ is not satisfied, then it is proven in Appendix \ref{resonances-asympt} that 
\begin{align*}
|\phi^{\alpha,\beta}_{m,n,p}(\xi,\lambda^{\alpha,\beta}_{m,n}\xi)|\geq \frac{1}{2\sqrt{n+1}(\sqrt{n+1}+\sqrt{m+1})\sqrt{(\lambda^{\alpha,\beta}_{m,n}\xi)^2+2m+2}}.
\end{align*}
Then a integration by parts in time is feasible, and leads to Lemma \ref{validity-approx-osc}.\\
Lemmas \ref{estimate-self-interaction-remainder}, \ref{spr-estimate} and \ref{validity-approx-osc} finally give Lemma \ref{lemma for approximation theorem}.

\end{prf}

Theorem \ref{approximation theorem} is then proved.

\part*{Appendix}
\appendix

\section{Resonant sets, asymptotics for $\phi$ and its derivatives}\label{resonances-asympt}
For now on, fix $m$, $n$, $p$ three integers, $\alpha$ and $\beta$ equal to $\pm 1$ and consider the phase 
\begin{equation*}
\phi(\xi,\eta)=\sqrt{\xi^2+2p+2}+\alpha\sqrt{\eta^2+2m+2}+\beta\sqrt{(\xi-\eta)^2+2n+2}.
\end{equation*}
In this article we are only stating the results. In particular we are not going to prove Theorem \ref{thm-phase}: for the detailed proofs, see \cite{These}.
For now on, assume that $\alpha=\beta=-1$, the other cases being dealt with similarly.

\subsection{Asymptotics for $\partial_{\eta}\phi$ and $\partial_{\xi}\phi$}\label{section-asymptotique-dphi}
Thanks to the central symmetry, let us focus only on the level lines under the space-time resonant set. This corresponds to negative values of $\partial_{\eta}\phi$ and positive ones of $\partial_{\xi}\phi$.\\
First notice that we have an explicit expression for level lines. The level line $\partial_{\eta}\phi:=-2^{-j}$ is 
\begin{equation*}
\xi=\eta+\sqrt{\frac{(2n+2)}{1-([\eta]_{m}-2^{-j})^2}}([\eta]_{m}-2^{-j}),
\end{equation*}
where $[\eta]_{m}:=\frac{\eta}{\sqrt{\eta^2+2m+2}}$. Similarly, the explicit expression for the level line $\partial_{\xi}\phi=2^{-j}$ is 
\begin{equation*}
\eta=\xi-\sqrt{\frac{2n+1}{1-([\xi]_{p}-2^{-j})}}([\xi]_{p}-2^{-j}).
\end{equation*}
We are going to rewrite these formulas in a more suitable way, adapted to $4$ different asymptotic regimes 
\begin{itemize}
\item the first one will be the asmptotics $|\eta|<<\sqrt{m}$: in this zone, the level lines are almost straight lines.
\item the second to fourth ones correspond to different order of magnitude of the asymptotic parameter
\begin{equation*}
\varrho(m,j,\eta):=\frac{\eta^2}{2^{j}m}.
\end{equation*}
\begin{itemize}
\item when this parameter is very small, we can compute the deviation with respect to the straight line of slope $\Lambda_{m,n}:=1+\alpha\beta\sqrt{\frac{n+1}{m+1}}$.
\item when the parameter is very large, level lines are like straight lines of slope $1$.
\item finally, when it is close to one, level lines are vertical lines.
\end{itemize}
\end{itemize}
It has to be noticed that the formulas written in the following sections are exact formulas: the only thing depending on the asymptotics will be the smallness of the remainder terms. The Lemma we state and which is proven in \cite{These} is the following one.

\begin{lem}
The following asymptotics for the phase $\phi$ occur for all $m,n,j,\eta$.
\begin{enumerate}
\item (low-frequency asymptotics) 
\begin{enumerate}
\item the equation of the level line $\partial_{\eta}\phi=2^{-j}$ can be rewritten as
\begin{equation}\label{deta-lf}
\xi=\left(1+\sqrt{\frac{n+1}{m+1}}\frac{1}{(1-2^{-2j})^{\frac{3}{2}}}\right)\eta\left(1+r(j,m,n,\eta)\right) -\frac{2^{-j}\sqrt{2n+2}}{\sqrt{1-2^{-2j}}},
\end{equation}
where $r(j,m,n,\eta)\leq  \frac{1}{\sqrt{m}} (1-2^{-2j})\eta$.
\item the equation of the level line $\partial_{\xi}\phi=2^j$ can be rewritten as
\begin{equation}\label{dxi-lf}
\eta=\left(1-\sqrt{\frac{n+1}{p+1}}\frac{1}{(1-2^{-2j})^{\frac{3}{2}}}\right)\xi\left(1+r_2\right) +\frac{2^{-j}\sqrt{2n+2}}{\sqrt{1-2^{-2j}}},
\end{equation}
where $r_2\leq  \frac{1}{\sqrt{p}} (1-2^{-2j})\xi$.
\item in the asymptotic zone $|\eta|<<\sqrt{m}$, the width of the band ${-2^{-j}\leq\partial_{\eta}\phi\leq -2^{-(j+1)}}$ is bounded by
\begin{equation}\label{width-eta-small}
C 2^{-j}\min(\sqrt{m},\sqrt{n})
\end{equation}
where $C$ is a constant independent of $j$, $n$ and $m$.
\end{enumerate} 
\item (high-frequency asymptotics with $\varrho$ small)
\begin{enumerate}
\item the level line $\partial_{\eta}\phi:=-2^{-j}$ can be rewritten as
\begin{equation}\label{deta-mf}
\xi=\eta+\sqrt{\frac{n+1}{m+1}}\eta-\sqrt{\frac{n+1}{m+1}}\sqrt{\eta^2+2m+2}\frac{2^{-j}\eta^2}{2m+2}+\sqrt{\frac{n+1}{m+1}}\sqrt{\eta^2+2m+2}r,
\end{equation}
where $r\lesssim \varrho^2$.
\item if $|\eta|\gtrsim\sqrt{m}$ and $\varrho\leq \varrho_0$ for $\varrho_0<1$ well-chosen, the width of the zone 
\begin{align*}
\left\{-2^{-j}\leq\partial_{\eta}\phi\leq -2^{-(j+1)},~2^{k}\leq |\eta|\leq 2^{k+1}\right\},
\end{align*}written $w_{\varrho<<1}(m,n,j,k)$ is bounded as follows.
\begin{equation}\label{width-rho-small}
w_{\varrho<<1}(m,n,j,k)\lesssim \frac{2^{3k}2^{-j}}{2m+2}.
\end{equation}
\end{enumerate}
\item (high-frequency asymptotics with $\varrho$ large)
\begin{enumerate}
\item the equation for the level line $\partial_{\eta}\phi=-2^ {j}$ can be rewritten
\begin{equation}\label{etainf}
\xi=\eta+\frac{\sqrt{2n+2}}{\sqrt{2^{1-j}-2^{-2j}}}(1-2^{-j}+r),
\end{equation}
with $r\leq c\frac{1}{2^{1-j}-2^{-2j}}\frac{m+1}{\eta^2}\lesssim \frac{1}{2-2^{-j}}\frac{1}{\varrho}$.
\item under the hypothesis $\varrho\geq \varrho_{0}$, $\varrho_0$ being chosen in (\ref{width-rho-small}), the width of the band ${-2^{-j}\leq\partial_{\eta}\phi\leq -2^{-(j+1)}}$ is bounded by
\begin{align}\label{width-rho-large}
w_{\varrho\gtrsim 1}(m,n,j)\lesssim 2^{\frac{j}{2}}\sqrt{2n+2}.
\end{align}
\end{enumerate}
\end{enumerate}
\end{lem}

\subsection{Comparison between $\partial_{\eta}\phi$ and $\partial_{\xi}\phi$}\label{section-comparison}

\begin{lem}\label{dxi<deta}
For all $\xi,\eta$ real numbers, we have
\begin{equation*}
\left|\partial_{\xi}\phi\right|\leq\left|\partial_{\eta}\phi\right|.
\end{equation*}
\end{lem}
\cacher{\begin{pr}
Prove the inequality in the zone $\xi\leq\left(1+\sqrt{\frac{n+1}{m+1}}\right)\eta$ (i.e. below the space-time resonant set). In this zone, $\partial_{\eta}\phi$ is negative and $\partial_{\xi}\phi$ is positive, which means that $\left|\partial_{\xi}\phi\right|=\frac{\xi}{\cro{\xi}_{p}}-\frac{\xi-\eta}{\cro{\xi-\eta}_{n}}$ and $\left|\partial_{\eta}\phi\right|=\frac{\eta}{\cro{\eta}_{m}}-\frac{\xi-\eta}{\cro{\xi-\eta}_{n}}$. Then the inequality to prove is 
\begin{equation*}
\frac{\xi}{\sqrt{\xi^{2}+2p+2}}\leq\frac{\eta}{\sqrt{\eta^{2}+2m+2}}.
\end{equation*}
Since $\xi\mapsto \frac{\xi}{\sqrt{\xi^{2}+2p+2}}$ is increasing, in the zone $\xi\leq\left(1+\sqrt{\frac{n+1}{m+1}}\right)\eta$ we have 
\begin{align*}
\frac{\xi}{\sqrt{\xi^{2}+2p+2}}&\leq\frac{\left(1+\sqrt{\frac{n+1}{m+1}}\right)\eta}{\sqrt{\left(1+\sqrt{\frac{n+1}{m+1}}\right)^{2}\eta^{2}+2p+2}}\\
&\leq \frac{\eta}{\sqrt{\eta^{2}+\frac{2p+2}{\left(1+\sqrt{\frac{n+1}{m+1}}\right)^{2}}}}.
\end{align*}
But thanks to the relation between $m$, $n$ and $p$, $\frac{2p+2}{\left(1+\sqrt{\frac{n+1}{m+1}}\right)^{2}}=2m+2$ and the inequality is proved.
\end{pr}
}

\subsection{Distance between $\Sc$ and $\T$}\label{section-distance}
Here we are in the case where 
\begin{equation*}
p\leq n+m\text{ or }p^2+m^2+n^2-2pm-2mn-2pn-2p-2m-2n-3\neq 0.
\end{equation*} 
In this situation we know that $\Sc$ and $\T$ do not intersect, and wonder in this section if we can determine the distance between those two setsin a given ball.\\
The first step will be to evaluate $\phi(\Lambda_{m,n}\eta,\eta)$ (which is no longer equal to $0$) and then to find the width of a neighborhood of the straight line $\xi-\Lambda_{m,n}\eta$ where $\phi$ remains different from $0$. We have the following proposition.

\begin{prop}\label{timres-width}
Let $R>0$. There exists a constant $c$ such that for all $(\xi,\eta)$ satisfying 

\begin{align*}
\sqrt{|\xi|^2+|\eta|^2}&\leq R,\\
\mathrm{dist}((\xi,\eta),\Sc)&\leq \frac{c}{(\sqrt{n+1}+\sqrt{m+1})^2}\frac{1}{R},
\end{align*}
then the modulus of the phase $|\phi|$ is bounded from below (up to a constant independent of $R$, $m$ and $n$) by
\begin{align*}
\frac{1}{(\sqrt{n+1}+\sqrt{m+1})^2}\frac{1}{R}.
\end{align*}
\end{prop}
\subsection{Asymptotics for $\phi$}\label{section-asymptotique-phi}
Now the remaining asymptotics we have to study is the asymptotics for $\phi$: our goal is to find a lower bound (depending only on $m$, $n$, $p$ and $R$) for $\phi$ on the zone $|\eta|,|\xi|\leq R$.\\

\begin{lem}\label{asymptphi}
In the case where there are no time resonances (i.e. $n>p$ or $m>p$) then we have the following lower bound on $\phi$.
\begin{align*}
|\phi|\gtrsim \frac{1}{(\sqrt{n+1}+\sqrt{m+1})^2}\frac{1}{R}.
\end{align*}
\end{lem}

\cacher{\begin{pr}
So as to study this phase, try to understand it along a given direction. Let $\gamma$ be a real number and try to understand $\phi(\gamma\eta,\eta)$.
\begin{align*}
\phi(\gamma\eta,\eta)=&\sqrt{\gamma^2\eta^2+2p+2}-\sqrt{\eta^2+2m+2}-\sqrt{(\gamma-1)^2\eta^2+2m+2}\\
=&|\gamma|\sqrt{\eta^2+\frac{2p+2}{\gamma^2}}-\sqrt{\eta^2+2m+2}-|\gamma-1|\sqrt{\eta^2+\frac{2m+2}{(\gamma-1)^2}}.
\end{align*}

From this formula comes the following alternative: if $\gamma<1$ then $\phi(\gamma\eta,\eta)$ goes to $-\infty$ when $\eta$ goes to $\pm\infty$ ; if $\gamma\geq 1$ then $\phi(\gamma\eta,\eta)$ goes to $0$ when $\eta$ goes to $\pm\infty$. Hence the asymptotics have to be done in that zone. Remark that in the proof of Theorem \ref{space-time-resonances} we proved that $|\phi|$ was minimal on the space resonant set.\\
Since we are in the case $p<m$ or $p<n$ we have $p<n+m+1$ and the previous lower bound can be used.
\begin{align*}
|\phi|&\gtrsim \frac{1}{2\sqrt{n+1}(\sqrt{n+1}+\sqrt{m+1})\sqrt{\eta^2+2m+2}}\\
&\gtrsim \frac{1}{2\sqrt{n+1}(\sqrt{n+1}+\sqrt{m+1})\sqrt{(K_{m,n}R)^2+2m+2}}.
\end{align*}
\end{pr}
}

\section{Some harmonic analysis tools}\label{Fourier-multipliers}
Some proofs in this section are skipped, they can be found in \cite{These}.\\
Before really studying Fourier multipliers, we state this useful lemma on Sobolev spaces.
\begin{lem}\label{leml1}
Let $f$ be in $H^{k}$, such that $xf$ is in $H^{k}$. Then $f$ is in $W^{k,1}$ and
\begin{equation*}
\nor{f}{W^{k,1}}\leq 4\sqrt{\norh{f}{k}{}\norh{\xc f}{k}{}}.
\end{equation*}
\end{lem}
\cacher{\begin{pr}
Let $A$ be a positive real number. By Cauchy-Schwarz inequality,
\begin{align*}
\inte{x\in\R}{}{|D^{k}f(x)|}{x}&=\inte{|x|\leq A}{}{|D^{k}f(x)|}{x}+\inte{|x|\geq A}{}{|D^{k}f(x)|}{x}\\
&=\inte{x\in\R}{}{\mathbf{1}_{|x|\leq A}|D^{k}f(x)|}{x}+\inte{x\in\R}{}{\mathbf{1}_{|x|\geq A}\frac{1}{|x|}|x||D^{k}f(x)|}{x}\\
&\leq 2A\norh{f}{k}{}+\frac{2}{A}\norlp{xD^{k}f}{2}{}\\
&\leq 2A\norh{f}{k}{}+\frac{2}{A}\norh{\xc f}{k}{}.
\end{align*}
Then optimize over $A$ to find $A=\sqrt{\frac{\norh{xf}{k}{}}{\norh{f}{k}{}}}$. This ends the proof.
\end{pr}
}
\subsection{Linear Fourier multiplier estimates}

\subsubsection{Dispersive estimates}
First write down the dispersive estimate for Klein-Gordon equation in dimension $1$ which can be found in \cite{HN12}.
\begin{prop}\label{disp}
 Let $k$ be a positive real number, $u_0$ be a function in the Sobolev space $W^{3/2,1}$. Then for all $t$, $e^{-it\cro{D}}u_0$ and $e^{-it\sqrt{-D^2+k}}u_0$ is in $L^{\infty}$ and the following inequalities hold for all $t>0$.
\begin{align}
 \norlp{e^{-it\cro{D}}u_0}{\infty}{}&\lesssim\frac{1}{\sqrt{\cro{t}}}\nor{u_0}{W^{3/2,1}}, \label{dispequ}\\
\norlp{e^{-it\sqrt{-\Delta+k}}u_0}{\infty}{}&\lesssim\frac{k^{\frac{1}{4}}}{\sqrt{\cro{t}}}\nor{u_0}{W^{3/2,1}}. \label{dispequpoids}
\end{align}
(Global dispersion) Let $f$ in $\Sigma^{N,M}_{t}$ for $M,N$ satisfying (\ref{cond-M})-(\ref{cond-N}). Then for all $t$, $e^{\pm it\sqrt{-\Delta+x_{2}^2+1}}f$ is in $L^{\infty}(\R^2)$ and 
\begin{equation}\label{disp-global}
\nor{e^{\pm it\sqrt{-\Delta+x_{2}^2+1}}f(t)}{L^{\infty}(\R^2)}\lesssim \cro{t}^{-\frac{1}{4}}\nor{f(t)}{S^{N,M}_{t}}.
\end{equation}
\end{prop}

\subsubsection{Other Fourier multiplier estimates.} Given the nature of the Duhamel formula, we will have to deal with mutlipliers of the form $\frac{1}{\cro{\eta}_{n}}$.

\begin{prop}\label{propmultgeneral}
\begin{itemize}
\item There exists a constant $C$ such that the following inequalities hold for all $1\leq q\leq \infty$, all $f$ in $L^q$ and all positive real number $\lambda$.
\begin{align}\tag{\ref{propmultgeneral}-a}\label{propmult}
\norlp{\frac{f}{\sqrt{D^2+\lambda}}}{q}{}\leq \frac{C}{\sqrt{\lambda}}\norlp{f}{q}{}.
\end{align}
\item There exists a constant $C$ such that the following inequalities hold for all integers $p$, all $f$ in $\He^N$,
\begin{align}\tag{\ref{propmultgeneral}-b}\label{propmult-global}
\nor{\frac{f}{\sqrt{-\Delta^2+x_{2}^2+1}}}{\He^N}\leq C\nor{f}{\He^N}.
\end{align}
\item There exists a constant $C$ such that the following inequalities hold for all $1<q<\infty$, all $f$ in $L^q$ and all positive real number $\lambda$:
\begin{equation}\tag{\ref{propmultgeneral}-c}\label{propmult2}
\norlp{\frac{D}{\sqrt{D^2+\lambda}}f}{q}{}\leq C\norlp{f}{q}{}.
\end{equation}
\item More generally, for all $a>1$, $1<q<\infty$, all $f$ in $L^q$ and all positive real number $\lambda$, we have the following estimate.\begin{equation}\tag{\ref{propmultgeneral}-d}\label{propmult3}
\norlp{\left(\frac{|D|}{\sqrt{D^2+\lambda}}\right)^{a}f}{q}{}\leq C_{a}\norlp{f}{q}{}.
\end{equation}
\item For all $M$, $N$ integers, for all $f$ in $H^N$, 
 \begin{equation}\tag{\ref{propmultgeneral}-e}\label{hf}
 \norlp{|D|^{M}(1-\theta_{R}(D))f}{2}{} \leq\frac{1}{R^{N-M}}\norh{f}{N}{},
\end{equation}
where $\theta_{R}(\xi)$ is localizing in the zone $|\xi|\leq R$.
\end{itemize}
\end{prop}

\subsubsection{Combination of dispersion and multiplier estimates}

\begin{prop}\label{disp+mult}
Let $n$ be an integer, $f$ in $B$ and $s>0$. Then we have the following inequality.
\begin{align}\label{disp+mult_equ}
\norlp{e^{is\cro{D}_{n}}\frac{f(s)}{\cro{D}_{n}}}{\infty}{}\lesssim s^{-\frac{1}{4}}n^{-\frac{1}{4}}\sqrt{\norh{f(s)}{\frac{3}{2}}{}\nor{f}{B_{s}}}.
\end{align}
\end{prop}

\cacher{\begin{pr}
First, by Proposition \ref{propmult}, 
\begin{align*}
\norlp{e^{is\cro{D}_{n}}\frac{f(s)}{\cro{D}_{n}}}{\infty}{}\lesssim\frac{1}{\sqrt{n}}\norlp{e^{is\cro{D}_{n}}f(s)}{\infty}{}.
\end{align*}
Then we can apply the dispersion estimate (\ref{dispequ}):
\begin{align*}
\norlp{e^{is\cro{D}_{n}}\frac{f}{\cro{D}_{n}}}{\infty}{}\lesssim\frac{1}{\sqrt{n}}\frac{n^{\frac{1}{4}}}{\sqrt{s}}\norw{f}{\frac{3}{2}}{1}{}.
\end{align*}
Finally, using lemma \ref{leml1} leads to
\begin{align*}
\norlp{e^{is\cro{D}_{n}}\frac{f(s)}{\cro{D}_{n}}}{\infty}{}\lesssim\frac{n^{-\frac{1}{4}}}{\sqrt{s}}\sqrt{\norh{f(s)}{\frac{3}{2}}{}\norh{xf}{\frac{3}{2}}{}}.
\end{align*}
Since $\norh{xf(s)}{\frac{3}{2}}{}\leq \sqrt{s}\nor{f}{B_{s}}$, we have proved the inequality.
\end{pr}
}

\subsection{Behavior with dilation operators}

\begin{defi}
Let $\lambda$ be a real number. Define the \emph{Fourier dilation operator} of parameter $\lambda$, written $\E_{\lambda}$ by $\E_{\lambda}f:=\F^{-1}\left(\xi\mapsto \hat{f}(\lambda \xi)\right)$.
\end{defi}

It is well known that $\left(\E_{\lambda}f\right)(x)=\frac{1}{\lambda}f\left(\frac{x}{\lambda}\right)$. The following lemma generalizes in the case where there is a Fourier multiplier.

\begin{lem}\label{multdilat}
Let $\lambda$ be a real number, $p$ an integer and $g$ the symbol of a Fourier multiplier. Then for all $f$ we have
\begin{equation*}
\norlp{\E_{\lambda}(g(D)f)}{p}{}=\lambda^{\frac{1}{p}-1}\norlp{g(D)f}{p}{}.
\end{equation*}
\end{lem}

\subsection{Bilinear multiplier estimates}\label{bilin}

When dealing with the cutoff functions as the ones defined in Definition \ref{cutoffs}, a question arises: how can we keep "Hölder-like" estimates ? More precisely, define, for $m(\eta,\xi)$ a smooth function of $\eta$ and $\xi$,
\begin{equation*}
 A_{m}(f,g)=\mathcal{F}^{-1}\left(\inte{\R}{}{m(\xi,\eta)\hat{f}(\eta)\hat{f}_{m}(\xi-\eta)}{\eta}\right).
\end{equation*}
The natural question is: is it possible to have an inequality like $\norlp{A_{m}(f,g)}{r}{}\leq C\norlp{f}{p}{}\norlp{f_{m}}{q}{}$, with $p,q,r$ satisfying some conditions (for example the Hölder condition $\frac{1}{r}=\frac{1}{p}+\frac{1}{q}$) ? Answering for a general $m$ is a very hard question, but some useful results are already known and work for the multipliers used in the three kinds of cut-offs.\\
This theory has been studied in the 60's by Coifman and Meyer \cite{CM}, and the now known as Coifman-Meyer estimates will be very useful in our situation (see \cite{MPTT05} for a proof).
\begin{thm}\label{CM}
 (Coifman-Meyer) Suppose that $m\in L^{\infty}(\R^2)$ is smooth away from the origin and satisfies
 
\begin{equation}\label{cond}
 |\partial_{\xi}^{\alpha}\partial_{\eta}^{\beta}m|\leq\frac{C}{(|\xi|+|\eta|)^{\alpha+\beta}},
\end{equation}
for all $\alpha, \beta \leq 3$ (we say that $m$ is a Coifman-Meyer symbol).\\
Then $A_{m}$ is bounded from $L^p\times L^q$ to $L^r$ where $\frac{1}{r}=\frac{1}{p}+\frac{1}{q}$, $1\leq p,q\leq \infty$, $1\leq r<\infty$.
\end{thm}
\begin{rmq}
The condition $(\ref{cond})$ is satisfies if $m$ is $\mathcal{C}^{\infty}$ and homogeneous of degree $0$ on a $(\xi,\eta)$ sphere.
\end{rmq}
For example, the symbol used in the \emph{high-low frequencies}  cut-off will satisfy the condition $(\ref{cond})$. But some other symbols with fail to satisfy the smoothness hypothesis of Coifman-Meyer's theorem: for example $\chi(\xi,\eta)\left|\frac{\xi}{\eta}\right|^{a}$ where $\chi(\xi,\eta)$ localizes around $|\xi-\eta|\leq 2|\eta|$:
\begin{prop}\label{CMex}
Let $a$ be a positive real number. Then the symbols
\begin{align*}
\chi(\xi,\eta)\left|\frac{\xi}{\eta}\right|^{a}~ \text{and}~\left(1-\chi(\xi,\eta)\right)\left|\frac{\xi}{\xi-\eta}\right|^{a}
\end{align*}
satisfy Hölder-like estimates.
\end{prop}

However the symbols used for the space resonant set cutoff will not fit in the conditions above. Another estimate, proven by Bernicot and Germain in \cite{BG2011} will be needed. First define the class $\mathcal{M}^{\Gamma}_{\eps}$.

\begin{defi}
The scalar-valued symbol $m_{\eps}$ belongs to the class $\mathcal{M}_{\eps}^{\Gamma}$ if 
\begin{itemize}
\item $\Gamma$ is a smooth curve in $\R^{2}$.
\item $m_{\eps}$ is supported in $B(0,1)$, as well as in a neighborhood of size $\eps$ of $\Gamma$.
\item The following inequality holds for sufficiently many indices $\alpha$ and $\beta$.
\begin{equation*}
\left|\partial_{\xi}^{\alpha}\partial_{\eta}^{\beta}m_{\eps}(\xi,\eta)\right|\lesssim \eps^{-\alpha-\beta}.
\end{equation*}
\end{itemize}
\end{defi}

\begin{thm}\label{thmger}
Consider $\Gamma$ a compact and smooth curve. Let $p,q,r\in[2,+\infty)$ be exponents satisfying $\frac{1}{p}+\frac{1}{q}+\frac{1}{r}-1\geq 0$. Then there exists a constant $C=C(p,q,r)$ such that for every $\eps>0$ and symbols $m_{\eps}\in\mathcal{M}^{\Gamma}_{\eps}$, then 
\begin{equation*}
\norlp{T_{m_{\eps}}(f,g)}{r'}{}\leq C\eps^{\frac{1}{p}+\frac{1}{q}+\frac{1}{r}-1}\norlp{f}{p}{}\norlp{g}{q}{}.
\end{equation*}
\end{thm}
As a consequence we have the following useful proposition.
\begin{defi}
The scalar-valued symbol $m_{\eps}$ belongs to the class $\mathcal{M}_{\eps,M}^{\Gamma}$ if 
\begin{itemize}
\item $\Gamma$ is a smooth curve in $\R^{2}$.
\item $m_{\eps}$ is supported in $B(0,1)$, as well as in a neighborhood of size $M\eps$ of $\Gamma$.
\item The following inequality holds for sufficiently many indices $\alpha$ and $\beta$.
\begin{equation*}
\left|\partial_{\xi}^{\alpha}\partial_{\eta}^{\beta}m_{\eps}(\xi,\eta)\right|\lesssim \eps^{-\alpha-\beta}.
\end{equation*}
\end{itemize}
\end{defi}
\begin{prop}\label{corger}
Consider $\Gamma$ a compact and smooth curve. Let $p,q,r\in[2,+\infty)$ be exponents satisfying $\frac{1}{p}+\frac{1}{q}+\frac{1}{r}-1\geq 0$. Let $M$ be a real number greater than $1$. Then there exists a constant $C=C(p,q,r)$ such that for every $\eps>0$ and symbols $m_{\eps}\in\mathcal{M}^{\Gamma}_{M,\eps}$, then 
\begin{equation*}
\norlp{T_{m_{\eps}}(f,g)}{r'}{}\leq C M\eps^{\frac{1}{p}+\frac{1}{q}+\frac{1}{r}-1}\norlp{f}{p}{}\norlp{g}{q}{}.
\end{equation*}
\end{prop}
\cacher{\begin{pr}
Consider $[M]$ functions $\chi_1,...\chi_{[M]}$ such that $\sum_{j}\chi_{j}=1$ and $m_{\eps}\chi_{j}$ is supported in a band of width $\leq\eps$. Then $m_{\eps}\chi_{j}$ belongs to $\mathcal{M}^{\Gamma}_{\eps}$. Then it satisfies the Bernicot-Germain inequality and by Minkowski inequality we have the proposition.
\end{pr}
}

In the case of a straight line, this version of the theorem is very useful.

\begin{thm}\label{estbilinband}
Let $\rho,\omega,\mu$ be real numbers, $\Gamma$ a straight line and $S$ be a symbol which satisfies the following properties.
\begin{itemize}
\item it is supported on a ball of radius $\rho$.
\item it is supported in a band of width $\omega$, around $\Gamma$.
\item it satisfies  the following estimate: for sufficiently many indices $a,b$, $|\partial_{\xi}^{a}\partial_{\eta}^{b}\tld{S}|\lesssim \left(\frac{\rho}{\mu}\right)^{a+b}$.
\end{itemize}
Then we have the following bilinear estimate.
\begin{align*}
\norlp{T_{S}(f,g)}{r'}{}\lesssim  \max\left(1,\frac{\omega}{\mu}\right)(\rho\omega)^{\frac{1}{p}+\frac{1}{q}+\frac{1}{r}-1}\norlp{f}{p}{}\norlp{g}{q}{}.
\end{align*}
\end{thm}

\subsection{A $L^2$ stationary phase lemma}
We state a stationary phase Proposition, first proven by Germain, Pusateri and Rousset (to appear) and adapted to our problem.
\begin{prop}\label{stat-phase}
Consider $\rho>0$, $\chi\in\mathcal{C}^{\infty}_{0}$ equal to zero on $B(0,\rho)^{c}$, such that $|\chi'|$ is bounded by $1/\rho$, and $\psi$ in $\mathcal{C}^{\infty}$. Let
\begin{align*}
I=\inte{\R}{}{e^{it\psi(x)}F(x)\chi}{x}.
\end{align*}
\begin{enumerate}
\item (non-stationary phase) Let $m$ be a positive real number such that for all $x\in\supp (\chi)$, $|\psi'(x)|\geq m$. Then
\begin{equation*}
|I|\leq \frac{\sqrt{\rho}}{t m}\left(\norlp{F}{2}{}+\norlp{F'}{2}{}\right).
\end{equation*}
\item (stationary phase)  Let $x_0$ be the only point where $\psi'(x_0)=0$. Let $m$ and $M$ two positive real numbers,  such that for all $x$ in $\mathrm{supp}(\chi)$,
\begin{align*}
&\psi''(x)\geq m,~|\psi'''(x)|\leq M,\\
&\text{or}~ \left|\frac{\sqrt{\psi-\psi(x_0)}}{\psi'}\right|\lesssim \frac{1}{\sqrt{m}},~|\psi'''(x)|\leq M.
\end{align*}
Then
\begin{align*}
I=\frac{Ce^{it\psi(x_0)}}{\psi''(x_0)\sqrt{t}}{\chi}(x_0){F}(x_0)+O\left(\frac{1}{t^\frac{3}{4}}\left(\left(\frac{1}{m^{\frac{1}{4}}\rho}+m^{\frac{1}{4}}\sqrt{M}\right)\norlp{\tld{F}}{2}{}+\frac{1}{m^{\frac{3}{4}}}\norlp{\tld{F}'}{2}{}\right)\right).
\end{align*}
the constants being independent of $\chi$ and $\psi$.
\end{enumerate}
\end{prop}

The proof of Proposition \ref{stat-phase} relies on a change of variable and will not be detailed here (see \cite{These}).

\subsection{Interaction between Hermite functions}

The integral
\begin{align*}
\mathcal{M}(m,n,p)=\inte{\R}{}{\psi_m(x_2)\psi_n(x_2)\psi_p(x_2)}{x_2}
\end{align*}
 can be computed explicitly but the exact formula is not really helpful to get estimates. 

\begin{prop}\label{interact_harm}
Let $\nu>1/8$ and $0\leq\beta<1/24$, $\eps>0$ and $0\leq\theta\leq 1$. Then for all $m\leq n\leq p$ and $K$ integers, there exists $C_K$, $C_{\eps}$ and $C_{\eps,\theta,K}$ three positive constants such that
\begin{align}
&\left|\mathcal{M}(m,n,p)\right| \leq C_K \frac{m^{\nu}}{p^{\beta}}\left(\frac{\sqrt{mn}}{\sqrt{mn}+p-n}\right)^K,\label{interact_harm-gip}\\
&\left|\mathcal{M}(m,n,p)\right| \leq C_{\eps} \max(m,n,p)^{-\frac{1}{4}+\eps},\label{interact_harm-faible}\\
&\left|\mathcal{M}(m,n,p)\right| \leq C_{\eps,\theta,K} \frac{m^{\theta\nu}}{p^{\theta\beta}}\left(\frac{\sqrt{mn}}{\sqrt{mn}+p-n}\right)^{\theta K}p^{-\frac{1}{4}+\frac{\theta}{4}+\theta\eps}.\label{interact_harm_interpol}
\end{align}
\end{prop}

\section{Paraproduct for the Hermite expansions}\label{Hermite-paraproduct}
In this appendix we are going to prove the following theorem.

\begin{thm}\label{resummation}
Let $a>0$, $\gamma>0$, $\left(a_{m}\right)_{m\in\N}$ and $\left(b_{n}\right)_{n\in\N}$ two sequences in $\ell^2$, $M>a+2$. Then there exists and a constant $C_{\gamma}$ and a sequence $(u_{p})_{p\in\N}$ in the unit ball of $\ell^2$ such that for all integer $p$,
\begin{enumerate}
\item \label{classicalresummation}
\begin{equation*}
p^{M}\sum_{m,n}\mathcal{M}(m,n,p)\frac{\max(m,n)^{a}}{\sqrt{mn}}m^{-M}a_m n^{-M}b_n\leq C_{\gamma}p^{a-\frac{3}{4}+\gamma}u_{p},
\end{equation*}
This inequality has two consequences:
\item (bounded sums theorem)\label{bddresummation}
for all $R>0$,
\begin{equation*}
p^{M}\sum_{m,n\leq R}\mathcal{M}(m,n,p)\frac{\max(m,n)^{a}}{\sqrt{mn}}m^{-M}a_m n^{-M}b_n\leq C_{\gamma}R^{a-\frac{3}{4}+\gamma}u_{p},
\end{equation*}
\item \label{halfresummation}(half sums theorem) for all $M_0>2$,
\begin{equation*}
p^{M}\sum_{m<n}\mathcal{M}(m,n,p)n^{a-\frac{1}{2}}m^{-M_0}a_m n^{-M}b_n\leq C_{\gamma}p^{a-\frac{3}{4}+\gamma}u_{p}.
\end{equation*}
\end{enumerate}
\end{thm}

\begin{pr}
We are going to deal with three different cases, corresponding to different orders of magnitude of the input frequencies $m$, $n$ and the output $p$ and use Proposition \ref{interact_harm}. 
\begin{enumerate}
\item Section \ref{p-great}. If $p> Cm$ and $p> Cm$, $C$ large enough chosen later, we will use the fact that the interaction term $\mathcal{M}(m,n,p)$ becomes very small. 
\item Section \ref{p-small}. If $p\leq Cn$ and $p\leq Cm$, we will simply use that $m^{a}\lesssim p^a$ if $a<0$. However three cases will have to be dealt with
\begin{enumerate}
 \item the case $p\leq m\leq n$.
 \item the case $m\leq p\leq n$.
\item the case $m\leq n\leq p$.
\end{enumerate}
\item Section \ref{p-inter}. If $Cm\leq p\leq Cn$ (the case $Cn\leq p\leq Cm$ being dealt with similarly given the symmetry of the situation), then we will try to use the interaction term as a convolution one.
\end{enumerate}

\subsection{If $p>>m$ and $p>>n$.}\label{p-great}

Let $C$ be a constant greater than $1$, and consider the following "low-low->high" term $S_{llh}$

\begin{equation*}
S_{llh}(m,n,p):=p^{M}\sum_{m\leq n\leq \frac{p}{C}}\mathcal{M}(m,n,p)\frac{\max(m,n)^{a}}{\sqrt{mn}}m^{-M}a_{m}n^{-M}b_{n},
\end{equation*}

the term 
\begin{equation*}
\tld{S}_{llh}:=p^{M}\sum_{n\leq m\leq \frac{p}{C}}\mathcal{M}(m,n,p)\frac{\max(m,n)^{a}}{\sqrt{mn}}m^{-M}a_{m}n^{-M}b_{n},
\end{equation*}
being dealt with similarly. First, by Proposition \ref{interact_harm}, if $K\in\N$, $\nu>1/8$ and $\beta<1/24$, we can write
\begin{equation*}
\mathcal{M}(m,n,p)\leq C_{K}\frac{m^{\nu}}{p^{\beta}}\left(\frac{\sqrt{mn}}{\sqrt{mn}+p-n}\right)^K.
\end{equation*}
Then, since $p>Cn$, we can bound $\frac{\sqrt{mn}}{\sqrt{mn}+p-n}$ by $\frac{\sqrt{mn}}{\left(1-\frac{1}{C}\right)p}$ and write
\begin{equation*}
\mathcal{M}(m,n,p)\leq C_{K}\frac{m^{\nu}}{p^{\beta}}m^{\frac{K}{2}}n^{\frac{K}{2}}\left(\frac{C}{C-1}\right)^{K}p^{-K}.
\end{equation*}
Then collect the terms in $m$, $n$ and $p$ in the original sum to get
\begin{align*}
S_{llh}(p)\leq C_{K}\left(\frac{C}{C-1}\right)^{K}p^{M-K-\beta}\left(\sum_{m}m^{\nu+\frac{K}{2}-\frac{1}{2}-M}a_m\right)\left(\sum_{n}n^{a+\frac{K}{2}-\frac{1}{2}-M}b_n\right).
\end{align*}
If 
\begin{equation*}
M>\frac{K}{2}+\max(\nu,a),
\end{equation*}i.e. 
\begin{equation*}
\left(m^{\nu+\frac{K}{2}-\frac{1}{2}-M}\right)_{m\in\N}\text{ and }\left(n^{a+\frac{K}{2}-\frac{1}{2}-M}\right)_{n\in\N}\text{ are  in }\ell^2,
\end{equation*}
then both series in $m$ and $n$ converge. Moreover, if $M<K-\frac{1}{2}$, $\left(S_{llh}(p)\right)_{p\in\N}$ is in $\ell^2$.

\subsection{If $p\leq Cm$ and $p\leq Cn$.}\label{p-small}

The zone $p\leq Cm$ and $p\leq Cn$ corresponds to a "low-low$\rightarrow$low", "high-high$\rightarrow$high" or a "high-high$\rightarrow$low" interaction. We will deal with the term $S_{1}$ defined by 
\begin{equation*}
S_{1}(p):=p^{M}\sum_{\frac{p}{C}\leq m\leq n}\mathcal{M}(m,n,p)\frac{\max(m,n)^{a}}{\sqrt{mn}}m^{-M}a_{m}n^{-M}b_{n},
\end{equation*}

the term 
\begin{equation*}
\tld{S}_{1}:=p^{M}\sum_{\frac{p}{C}\leq n\leq m}\mathcal{M}(m,n,p)\frac{\max(m,n)^{a}}{\sqrt{mn}}m^{-M}a_{m}n^{-M}b_{n},
\end{equation*}
being dealt with similarly. \\
Here, we do not need to find a fine bound for the interaction term $\mathcal{M}(m,n,p)$, we will simply bound it by a constant $\M_{0}$. Hence, collecting the terms in $m$, $n$ and $p$ we get
\begin{equation*}
S_{1}(p)\lesssim p^{M}\left(\sum_{\frac{p}{C}\leq m}m^{-M-\frac{1}{2}}\right)\left(\sum_{\frac{p}{C}\leq n}n^{a-M-\frac{1}{2}}\right).
\end{equation*}
Then, Hölder's inequality and a comparison series-integral give 
\begin{align*}
\sum_{\frac{p}{C}\leq m}m^{-M-\frac{1}{2}}a_m&\lesssim \left(\sum_{\frac{p}{C}\leq m}m^{-2M-1}\right)^\frac{1}{2}\lesssim C^{-M}p^{-M}\\
\sum_{\frac{p}{C}\leq n}n^{a-M-\frac{1}{2}}b_n&\lesssim \left(\sum_{\frac{p}{C}\leq n}n^{2a-2M-1}\right)^\frac{1}{2}\lesssim C^{a-M}p^{a-M}.
\end{align*}
Finally we get
\begin{equation*}
S_{1}(p)\lesssim p^{a-M},
\end{equation*}
which is in $\ell^2$ for $M$ large enough.

\subsection{If $p\geq Cm$ and $p\leq Cn$ or $p\geq Cn$ and $p\leq Cm$.\label{p-inter}}

\subsubsection{Classical paraproduct [Theorem \ref{resummation}-(\ref{classicalresummation})]}
First assume that $Cm\leq p\leq Cn$, the case $Cn\leq p\leq Cm$ being symmetric. Denote bu $S_{lhh}$ the term
\begin{equation*}
S_{lhh}:=p^{M}\sum_{Cm\leq p\leq n}\mathcal{M}(m,n,p)\frac{\max(m,n)^{a}}{\sqrt{mn}}m^{-M}a_{m}n^{-M}b_{n}.
\end{equation*}
Then assume that $Cm\leq p\leq n$: the case $Cm\leq n\leq p\leq Cn$ is dealt with similarly, simply by multiplying by powers of $C$.\\

Let $\theta,\eps>0$, $\nu>1/8$, $0\leq\beta<1/24$ and $K$ integer, there exists a constant $C_{\eps,\theta,K}$ and a sequence $\left(u_{p}\right)_{p\in\N}$ in $\ell^2$ such that 
\begin{align*}
\M(m,n,p) &\leq C_{\eps,\theta,K} \frac{m^{\theta\nu}}{n^{\theta\beta}}\left(\frac{\sqrt{mp}}{\sqrt{mp}+n-p}\right)^{\theta K}n^{-\frac{1}{4}+\frac{\theta}{4}+\theta\eps}u_{p}^{1-\theta}\\
&\leq C_{\eps,\theta,K}  \frac{m^{\theta\nu}}{n^{\theta\beta}}\left(\frac{\sqrt{mp}}{1+n-p}\right)^{\theta K}n^{-\frac{1}{4}+\frac{\theta}{4}+\theta\eps}u_{p}^{1-\theta}.
\end{align*}
Then, when gathering the terms in $m$, $n$ and $p$ we obtain
\begin{equation*}
S_{lhh}\leq  C_{\eps,\theta,K} u_{p}^{1-\theta}p^{M+\frac{\theta K}{2}}\left(\sum_{Cm\leq p}a_{m}m^{\theta\nu+\frac{\theta K}{2}-M-\frac{1}{2}}\right)\left(\sum_{p\leq n}\left(\frac{1}{1+n-p}\right)^{\theta K}b_n n^{a-M-\frac{1}{2}-\theta\beta-\frac{1}{4}+\frac{\theta}{4}+\theta\eps}\right).
\end{equation*}

Then the sum
\begin{equation*}
\sum_{Cm\leq p}m^{\theta\nu+\frac{\theta K}{2}-M-\frac{1}{2}}
\end{equation*}
is finite whenever $M$ is large enough.\\
In order to bound 
\begin{equation*}
\sum_{p\leq n}\left(\frac{1}{1+n-p}\right)^{\theta K}n^{a-M-\frac{1}{2}-\theta\beta-\frac{1}{4}+\frac{\theta}{4}+\theta\eps}b_{n},
\end{equation*}
first we bound $n^{-1}$ by $p^{-1}$ and write
\begin{equation*}
\sum_{p\leq n}\left(\frac{1}{1+n-p}\right)^{\theta K}b_n n^{a-M-\frac{1}{2}-\theta\beta-\frac{1}{4}+\frac{\theta}{4}+\theta\eps}\leq p^{a-M-\frac{1}{2}-\theta\beta-\frac{1}{4}+\frac{\theta}{4}+\theta\eps} \sum_{n\geq p}\left(\frac{1}{1+n-p}\right)^{\theta K}b_n.
\end{equation*}
Then $b_n$ is in $\ell^2$ and so is $\left(\frac{1}{1+n}\right)_{n\in\N}$, whenever $2K\theta=1+\delta$, $\delta>0$. Finally the following bound holds for $S_{lhh}$.
\begin{align*}
S_{lhh}(p)&\lesssim u_{p}^{1-\theta}p^{M+\frac{\theta K}{2}}p^{a-M-\frac{1}{2}-\theta\beta-\frac{1}{4}+\frac{\theta}{4}+\theta\eps}\\
&\lesssim u_{p}^{1-\theta}p^{M+\frac{1}{4}+\frac{\delta}{2}}p^{a-M-\frac{1}{2}-\theta\beta-\frac{1}{4}+\frac{\theta}{4}+\theta\eps}\\
&\lesssim u_{p}^{1-\theta} p^{a-\frac{1}{2}+\theta(\eps+\frac{1}{4}-\beta)+\frac{\delta}{2}}.
\end{align*}
The remaining problem is that $\left(u_{p}^{1-\theta}\right)_{p\in\N}$ is not in $\ell^2$. However, remark that for all $\alpha>0$, $\left(u_{p}^{1-\theta}p^{-\frac{\theta}{2}(1+\alpha)}\right)_{p\in\N}$ is in $\ell^2$: it can be checked by writing
\begin{align*}
\sum_{p} u_{p}^{2-2\theta}p^{-2\theta(1+\alpha)} \lesssim \left(\sum_{p}\left(u_{p}^{2-2\theta}\right)^{\frac{2}{2-2\theta}}\right)^{\frac{2-\theta}{2}}\left(\sum_{p} \left(p^{-2\frac{\theta}{2}(1+\alpha)}\right)\right)^{\frac{1}{\theta}},
\end{align*}
which is finite. Finally, writing $v_{p}:=u_{p}^{1-\theta}p^{-\frac{\theta}{2}(1+\alpha)}$, the following bound holds:
\begin{equation*}
S_{lhh}(p)\lesssim v_{p}p^{a-\frac{1}{2}+\gamma},
\end{equation*}
with $\left(v_p\right)_{p\in\N}\in\ell^2$ and $\gamma=\theta(\frac{1+\alpha}{2}+\eps+\frac{1}{4}-\beta)+\frac{\delta}{2}$.

\subsubsection{Case of a bounded sum [Theorem \ref{resummation}-(\ref{bddresummation})] }
In this case, we are summing over all $m$ and $n$ less than or equal to $R$. Since we are in the situation ``$p\geq Cm$ and $p\leq Cn$ or $p\geq Cn$ and $p\leq Cm$'', $p$ can be bounded by $R$ and the bounded resummation theorem follows.\\
\\
Theorem \ref{resummation}-(\ref{halfresummation}) is skipped: this ends the proof.
\end{pr}

\bibliographystyle{plain}
\bibliography{bibliothese}
\end{document}